\documentclass[a4paper,10pt]{amsart}
\usepackage[english]{babel}
\usepackage[T1]{fontenc} 
\usepackage[utf8]{inputenc}
\allowdisplaybreaks

\usepackage{xcolor}

\usepackage[all]{xy}
\usepackage{amsfonts}
\usepackage{amsmath}
\usepackage{amsthm}
\usepackage{amssymb}
\usepackage{hyperref}
\usepackage{stmaryrd}

\numberwithin{equation}{section}

\newtheorem{theorem}[equation]{Theorem}

\newtheorem{lemma}[equation]{Lemma}
\newtheorem{proposition}[equation]{Proposition}

\theoremstyle{definition}
\newtheorem{definition}[equation]{Definition}

\theoremstyle{remark}
\newtheorem{remark}[equation]{Remark}
\newcommand{\Fq}{\mathbb{F}_{q}}
\renewcommand{\P}{\mathbb{P}}

\renewcommand{\k}{\mathbb{K}}

\newcommand\M{\mathfrak M}
\newcommand\Pic{\mathfrak Pic}
\newcommand\LS{\mathfrak LS}
\newcommand{\gs}{\mathfrak{g}}
\renewcommand{\O}{\mathcal{O}}
\newcommand{\Goppa}{\operatorname{Goppa}}

\DeclareMathOperator{\Diag}{Diag}
\DeclareMathOperator{\Gr}{Gr}
\DeclareMathOperator{\Res}{Res}
\DeclareMathOperator{\Spec}{Spec}
\DeclareMathOperator{\Sch}{Sch}

\newcommand{\Ht}{\mathbb{G}_m^{n-1}}

\DeclareMathOperator{\PGL}{PGL}
\DeclareMathOperator{\pr}{pr}
\DeclareMathOperator{\ev}{ev}
\DeclareMathOperator{\im}{Im}
\DeclareMathOperator{\fix}{Fix}
\DeclareMathOperator{\Supp}{Supp}
\DeclareMathOperator{\Isom}{Isom}

\DeclareMathOperator{\Hom}{Hom}
\DeclareMathOperator{\Mat}{Mat}
\DeclareMathOperator{\PicG}{Pic}
\DeclareMathOperator{\GL}{GL}
\DeclareMathOperator{\Proj}{Proj}

\DeclareMathOperator{\Sets}{Set}
\DeclareMathOperator{\Hilb}{Hilb}

\usepackage{tikz}

\definecolor{xicurve}{RGB}{70,60,170}
\definecolor{xipos}{RGB}{20,120,90}
\definecolor{xineg}{RGB}{130,130,125}
\definecolor{xiband}{gray}{0.90}
\definecolor{xinone}{RGB}{190,215,240}
\definecolor{xipart}{RGB}{247,200,120}
\definecolor{xifull}{RGB}{150,215,190}

\newcommand{\xiPW}{3.55}
\newcommand{\xiPH}{2.05}

\newcommand{\xipanel}[3]{%
\begin{scope}[shift={#3}]
  \pgfmathsetmacro{\xiA}{#2-2+2*(#1)}%
  \pgfmathsetmacro{\xiC}{\xiA*((#1)+1)-((#1)-1)^2}%
  \pgfmathsetmacro{\xiD}{(#2-4)^2-16*(#1)}%
  \pgfmathsetmacro{\xiLo}{2*(#1)-1}%
  \pgfmathsetmacro{\xiRm}{(\xiA-sqrt(max(\xiD,0)))/2}%
  \pgfmathsetmacro{\xiRp}{(\xiA+sqrt(max(\xiD,0)))/2}%
  \pgfmathsetmacro{\xiPd}{max(1,0.16*(#2-\xiLo))}%
  \pgfmathsetmacro{\xiXa}{\xiD>0 ? min(\xiLo-\xiPd,\xiRm-0.5) : \xiLo-\xiPd}%
  \pgfmathsetmacro{\xiXb}{\xiD>0 ? max(#2+\xiPd,\xiRp+0.5) : #2+\xiPd}%
  \pgfmathsetmacro{\xiFa}{-\xiXa*\xiXa+\xiA*\xiXa-\xiC}%
  \pgfmathsetmacro{\xiFb}{-\xiXb*\xiXb+\xiA*\xiXb-\xiC}%
  \pgfmathsetmacro{\xiMn}{min(min(\xiFa,\xiFb),min(\xiD/4,0))}%
  \pgfmathsetmacro{\xiMx}{max(max(\xiFa,\xiFb),max(\xiD/4,0))}%
  \pgfmathsetmacro{\xiYa}{\xiMn-0.10*(\xiMx-\xiMn)}%
  \pgfmathsetmacro{\xiYb}{\xiMx+0.15*(\xiMx-\xiMn)}%
  \pgfmathsetmacro{\xiSx}{\xiPW/(\xiXb-\xiXa)}%
  \pgfmathsetmacro{\xiSy}{\xiPH/(\xiYb-\xiYa)}%
  \pgfmathsetmacro{\xiZr}{(0-\xiYa)*\xiSy}%
  \pgfmathsetmacro{\xiAx}{(\xiA/2-\xiXa)*\xiSx}%
  \pgfmathsetmacro{\xiBa}{(\xiLo-\xiXa)*\xiSx}%
  \pgfmathsetmacro{\xiBb}{(#2-\xiXa)*\xiSx}%
  \fill[xiband] (\xiBa,0) rectangle (\xiBb,\xiPH);
  \draw[xineg,thin] (0,\xiZr) -- (\xiPW,\xiZr);
  \draw[xineg,thin,dashed] (\xiAx,0) -- (\xiAx,\xiPH);
  \draw[xicurve,thick,smooth]
    plot[domain=\xiXa:\xiXb,samples=60]
    ({(\x-\xiXa)*\xiSx},{(-\x*\x+\xiA*\x-\xiC-\xiYa)*\xiSy});
  \ifdim\xiD pt>-0.001pt
    \draw[xicurve,thin,fill=white]
      ({(\xiRm-\xiXa)*\xiSx},\xiZr) circle[radius=1.5pt];
    \draw[xicurve,thin,fill=white]
      ({(\xiRp-\xiXa)*\xiSx},\xiZr) circle[radius=1.5pt];
  \fi
  \pgfmathtruncatemacro{\xiIa}{2*(#1)}%
  \pgfmathtruncatemacro{\xiIb}{#2-1}%
  \foreach \xiI in {\xiIa,...,\xiIb}{%
    \pgfmathsetmacro{\xiV}{-\xiI*\xiI+\xiA*\xiI-\xiC}%
    \ifdim\xiV pt>0pt
      \fill[xipos]
        ({(\xiI-\xiXa)*\xiSx},{(\xiV-\xiYa)*\xiSy})
        circle[radius=1.4pt];
    \else
      \fill[xineg]
        ({(\xiI-\xiXa)*\xiSx},{(\xiV-\xiYa)*\xiSy})
        circle[radius=1.4pt];
    \fi}%
  \node[font=\scriptsize,anchor=south]
    at (\xiPW/2,\xiPH+0.06) {$g=#1$, $n=#2$};
  \node[font=\tiny,anchor=north]
    at (\xiBa,-0.04) {$2g-1$};
  \node[font=\tiny,anchor=north]
    at (\xiBb,-0.04) {$n$};
\end{scope}}

\newcommand{\xilattice}{%
  \foreach \laN in {3,...,32}{%
    \foreach \laG in {0,...,14}{%
      \pgfmathtruncatemacro{\laAdm}{\laN>=2*\laG+1 ? 1 : 0}%
      \ifnum\laAdm=1
        \pgfmathsetmacro{\laD}{(\laN-4)^2-16*\laG}%
        \pgfmathsetmacro{\laT}{mod(\laN,2)==0 ? 0 : 1}%
        \pgfmathsetmacro{\laQ}{(\laG-2)*\laN-\laG*\laG-4*\laG+4}%
        \pgfmathtruncatemacro{\laK}
          {\laD<=\laT ? 0 : (\laQ>0 ? 2 : 1)}%
        \ifcase\laK
          \fill[xinone]
        \or
          \fill[xipart]
        \else
          \fill[xifull]
        \fi
          ({0.335*(\laN-3)},{0.215*\laG})
          rectangle ++(0.295,0.17);
      \fi}}%
  \draw[xicurve,thick,dashed,smooth]
    plot[domain=4:19,samples=60]
    ({0.335*(\x-3)+0.148},{0.215*((\x-4)^2/16)+0.085});
  \foreach \laN in {4,8,12,16,20,24,28,32}{%
    \node[font=\tiny,anchor=north]
      at ({0.335*(\laN-3)+0.148},-0.05) {$\laN$};}
  \foreach \laG in {0,2,...,14}{%
    \node[font=\tiny,anchor=east]
      at (-0.08,{0.215*\laG+0.085}) {$\laG$};}
  \node[font=\tiny,anchor=east]
    at (-0.08,{0.215*15+0.24}) {$g$};
  \node[font=\tiny,anchor=west]
    at (10.20,0.085) {$n$};}

\begin{document}
\emergencystretch 3em

\title{The Goppa morphism}

\author{\'Angel Luis Mu\~noz Casta\~neda}
\address{Departamento de Matem\'aticas, Universidad de Le\'on}
\email{amunc@unileon.es}
\thanks{This work was initiated under the research project TED2021-131158A-I00, funded by MCIN/AEI/10.13039/501100011033 and by the European Union “NextGenerationEU”/PRTR. Subsequent improvements and extensions have been supported by the project PID2023-150787NB-I00, funded by MICIU/AEI /10.13039/501100011033 /FEDER, UE}

\subjclass[2020]{14H10; 14H60; 94A60; 94B27}

\date{}

\keywords{Level structures, moduli spaces, Goppa codes}

\begin{abstract}
We develop a moduli-theoretic framework, over an arbitrary base field, for the classical Goppa construction of algebraic-geometric codes.
We consider the
moduli stack $\LS_{g,n,d}$ of line bundles of degree $d$ on $n$-pointed smooth projective curves of genus $g$, trivialized along the marked points, as the space of Goppa structures and show that evaluation of global sections defines a natural Goppa morphism $\Goppa_{g,n,d}\colon\LS_{g,n,d}\to\Gr(k,n)$, $k=1-g+d$, together with an
extended morphism to $\Gr(k,n)\times\M_{g,n}$.

Under explicit numerical hypotheses, we prove that the extended Goppa
morphism is an immersion. The pullbacks of the Pl\"ucker coordinates are identified with sections of the theta line bundle on the relative Picard scheme, and the locus of non-degenerate codes for which a subset $I$ fails to be an information set is precisely the corresponding theta locus. We identify the fibers of the Goppa morphism over non-degenerate codes $C$ with moduli spaces of non-degenerate $n$-pointed genus-$g$ curves of degree $d$ in
the projective space $\P_C$ canonically determined by $C$. This yields a
Hilbert-incidence description of the Goppa construction and, in a natural range, implies that $\Goppa_{g,n,d}$ is schematic and quasi-projective and
that $\LS_{g,n,d}$ is a quasi-projective scheme. 

In genus zero, we further obtain, for $2\leq d\leq n-3$, a canonical $\mathbb G_m^{n-1}$-family of immersions $\M_{0,n}\hookrightarrow\Gr(d+1,n)$, which connects the Goppa morphism with the classical geometry of $\M_{0,n}$.
\end{abstract}

\maketitle

\tableofcontents

\section{Introduction}

The interaction between algebraic geometry and coding theory goes back to
Goppa's work in the 1970s. In~\cite{Gop77}, Goppa associated linear codes to divisors on algebraic function fields and used the Riemann--Roch theorem to
control their parameters. The construction was developed further in~\cite{Gop82}. Its significance became especially clear with the theorem of
Tsfasman, Vl\u{a}du\c{t} and Zink~\cite{tsfasman-vladut-zink}. They used modular and Shimura curves to construct sequences of algebraic-geometric codes over fields of square cardinality $q\geq 49$ whose asymptotic parameters
exceed the Gilbert--Varshamov bound. This result established algebraic
geometry as a source of coding-theoretic phenomena that were not visible from the classical constructions alone.
 
The role of algebraic geometry in coding theory expanded considerably in the following years. Divisors and differentials entered the study of duality and self-duality~\cite{stichtenoth-selfdual,driencourt-stichtenoth}. Algebraic methods became part of decoding, as in the work of Feng and Rao \cite{feng-rao}. The question of when two divisors give equal geometric Goppa
codes was studied by Munuera and Pellikaan \cite{munuera-pellikaan}. They also
introduced the generalized Picard group $\PicG(X,D)$ in this context to 
parametrize equivalence classes of strongly algebraic-geometric
codes associated with a given pointed curve \cite[\S~7]{munuera-pellikaan}. Pellikaan, Shen and van Wee \cite{pellikaan-which} addressed the basic characterization problem, and delineated the class of codes that the evaluation construction can produce. More recently, projective geometry and the equations of embedded
curves have been used to recover geometric data from a code \cite{marquez4}.
This reconstruction problem also became relevant in code-based cryptography~\cite{couvreur-mcc-ag}. Standard accounts of the subject and of many of these developments can be found in \cite{Tsfasman,stichtenoth}.
 
Much of this history proceeds by applying algebraic geometry to the
construction and analysis of codes. The present article adds a complementary direction. The geometric data of a fixed numerical type form a moduli stack, and the codes they produce are organized by an algebraic morphism from this stack to
a Grassmannian. Once this is made
explicit, questions about representations of codes become questions about a
morphism of moduli spaces.
 
We recall the construction in the form used throughout the article. Let $X$ be a geometrically irreducible smooth projective curve of genus $g$ over a field $\k$. Let $p_1,\dots,p_n\in X(\k)$ be pairwise distinct. Let $L$ be a
line bundle of degree $d$ and let $\gamma_i:L|_{p_i}\simeq\k$ be trivializations. If $n>d>2g-2$, evaluation
defines an injective map
\[
\begin{split}
H^0(X,L)&\longrightarrow \k^n,\\
s&\longmapsto
\bigl(\gamma_1(s(p_1)),\dots,\gamma_n(s(p_n))\bigr),
\end{split}
\]
and $h^0(X,L)=1-g+d=k$. Its image is a $k$-dimensional subspace of $\k^n$,
that is, a point of $\Gr(k,n)$ and, in coding-theoretic language, a linear
code.
 
The datum
$(X,p_1,\dots,p_n,L,\gamma_1,\dots,\gamma_n)$ is a rank-one level
structure in the sense of the Tsfasman--Vl\u{a}du\c{t}
$H$-construction~\cite[\S~3.1.1, p.~272]{Tsfasman}. The trivializations
are intrinsic to the line-bundle formulation of the evaluation construction. For a divisor $G$ whose support is disjoint from the marked points, the classical realization $L=\O_X(G)$ carries canonical identifications of the
fibers $L|_{p_i}$ with $\k$. For an arbitrary line bundle $L$, however, evaluation naturally takes values in
$\bigoplus_i L|_{p_i}$, and identifying this vector space with $\k^n$
requires trivializations $\gamma_i:L|_{p_i}\simeq\k$. Changing these
trivializations acts by coordinate scaling. The common scalar acts trivially
on the Grassmannian. Thus the level structure rigidifies precisely the ambiguity that appears when the construction is formulated intrinsically in terms of line bundles. The same data also behave naturally in families. Curves, marked points,
line bundles and trivializations may all vary over a base, and relative
evaluation then varies functorially.
 
For fixed $(g,n,d)$, these data form a moduli stack $\LS_{g,n,d}$. Relative
evaluation defines the morphism
\[
\Goppa_{g,n,d}\colon\LS_{g,n,d}\longrightarrow\Gr(k,n),
\]
which we call the \emph{Goppa morphism}. This passage does not replace the
classical Goppa construction by a different construction. It organizes the
same evaluation procedure over its natural moduli space. The construction
itself becomes an object of algebraic geometry.
 
Fixing the numerical type is essential, and not merely a convenient normalization. Pellikaan, Shen and van Wee~\cite{pellikaan-which} proved that
every linear code over $\Fq$ admits a weak algebraic-geometric representation
once the degree is allowed to vary without the usual numerical restrictions. The unrestricted evaluation construction is therefore too broad to single out
a class of linear codes, and the meaningful questions concern representations subject to additional geometric and numerical conditions. In the range
considered here the type $(g,n,d)$ is fixed and $H^0(X,L)$ has dimension $k=1-g+d$. The corresponding family is then governed by the geometry
of $\Goppa_{g,n,d}$.
 
This viewpoint gives a common formulation for several questions that have
appeared separately in the coding-theory literature. The image of $\Goppa_{g,n,d}$ is the locus of codes admitting a representation of the
prescribed numerical type. A fiber parametrizes the geometric representations
of a fixed code. Injectivity expresses uniqueness of such a representation. The behaviour of the morphism in families records whether these properties persist after extension of the base field or, more generally, after base change. Self-duality can also be expressed on the source as a fixed-point
condition for a natural involution. These are different questions, but they
are invariants of the same morphism.
 
There is also a direct projective interpretation, connected to the classical
projective-system viewpoint on a linear code. See, for example, \cite{Tsfasman}. Let $C\subset\k^n$ be a non-degenerate code. The nonzero functionals obtained by restricting the coordinate projections to $C$ define
labelled points $q_1,\dots,q_n$ in the projective space $\P_C$. Note that these points need not be pairwise distinct. After choosing a basis of $C$, they are precisely the projective points represented by the columns of a generator matrix. When $C$ is obtained from a sufficiently positive level structure, the
complete linear series of $L$ embeds the underlying curve in $\P_C$ and sends $p_i$ to $q_i$. In this case the $q_i$ are pairwise distinct. The inverse
problem for a code is then converted into a classical incidence problem. One seeks curves of prescribed degree and genus in a fixed projective space
through prescribed marked points. This observation is the geometric core of
our (functorial) description of the fibers.

Theta functions have already appeared in this subject. For instance, Henocq and Rotillon used the theta divisor of the Jacobian in decoding algorithms for geometric Goppa codes~\cite{henocq-rotillon}. Here theta functions play a different role, describing not the decoding algorithm but the Pl\"ucker coordinates that cut out the image of the Goppa morphism. We show that each Pl\"ucker coordinate, pulled back along the Goppa morphism, is a section of the theta line bundle on the relative Picard scheme of the curve, vanishing exactly where the corresponding coordinate 
$k$-subset fails to be an information set.
 
This geometry is also relevant to coding theory and cryptography. Universal injectivity of the Goppa morphism is a moduli-theoretic form of the \emph{unique representation}
problem. It asks whether the code determines the geometric data that produce it. For very strong algebraic-geometric codes over finite fields, this problem
was studied by M\'arquez-Corbella, Mart\'inez-Moro and Pellikaan \cite{marquez4}. Reconstruction of the hidden geometric structure is also the
basis of structural attacks against McEliece-type cryptosystems built from
structured codes. For algebraic-geometric codes this direction includes the
proposal of Janwa and Moreno~\cite{janwa-moreno} and the structural attack of Couvreur, M\'arquez-Corbella and Pellikaan~\cite{couvreur-mcc-ag}. The genus-zero case is closely related to the
classical attack of Sidelnikov and Shestakov on generalized Reed--Solomon codes~\cite{sidelnikov-shestakov}.
 
The description of a fiber addresses a related inverse question. It determines all geometric data of the prescribed type that produce a given code. This is close to the characterization problems studied for algebraic-geometric codes. The function-field version is also relevant to convolutional codes.
Convolutional codes of Goppa type were introduced by Dom\'{\i}nguez P\'erez, Mu\~noz Porras and Serrano Sotelo using families of algebraic curves
\cite{conv-goppa-type}. They were subsequently reformulated as Goppa codes on curves over $\Fq(z)$ in \cite{convolutional-goppa}. Further developments
include the characterization theorem in the convolutional setting \cite{curto-convolutional}. Working over an arbitrary base field lets the
present formalism treat finite fields and fields such as $\Fq(z)$ in the same language. This generality matters because $\Fq(z)$ is not perfect, so purely
inseparable phenomena cannot be excluded a priori.
 
The article establishes this moduli-theoretic framework and determines several
basic geometric properties of the resulting morphism. The contents are
described in \S~\ref{subsec:organization}.

\subsection{Main results}

We now state the main results. The first theorem gives the global form of the
classical evaluation construction.

\begin{theorem}[Theorem~\ref{prop:goppa-morphism-functorial}]\label{thm:intro:goppa-morphism}
If $n>d>2g-2$ and $k=1-g+d$, there exists a morphism of stacks
$\Goppa_{g,n,d}:\LS_{g,n,d}\longrightarrow \Gr(k,n)$
which assigns to each level structure its associated evaluation code.
\end{theorem}

The source is a smooth Deligne--Mumford stack in the stable range
$2g-2+n>0$ (Theorem~\ref{thm:structure-LS}). More precisely, its geometry is controlled by the rigidified universal Picard stack together with the trivializations along the marked divisor. The
next result shows that the evaluation code retains a substantial part of this
modular information.

\begin{theorem}[Lemma~\ref{alvarez-porras}, Theorems~\ref{th:immersionstacks} and~\ref{thm:universal-injectivity}]\label{thm:intro:mono-immersion}
Assume $d>2g-1$ and $(n-d)(d-g)>d$ (in particular, $n/2>d>2g-1$ with $d\geq1$ suffices). Let $(X\to S,\sigma_1,\dots,\sigma_n)$ be a family of
$n$-pointed smooth projective curves of genus $g$, set $D=\sigma_1+\cdots+\sigma_n$, and let $\LS_{X,D,d}$ be the associated space of level structures. Then the extended relative Goppa morphism
$\Phi_{X,D,d}:\LS_{X,D,d}\longrightarrow\Gr(k,n)\times S$
is a monomorphism. If, in addition, $n\ge d+g+2$, the extended Goppa morphism
$\Phi_{g,n,d}:\LS_{g,n,d}\longrightarrow\Gr(k,n)\times\M_{g,n}$ is an immersion of stacks. Moreover, if $n/2>d>2g+1$, the morphism $\Goppa_{g,n,d}$ is universally injective.
\end{theorem}

For a fixed curve over a field, the corresponding Grassmannian construction was given by \'Alvarez-V\'azquez and Mu\~noz Porras \cite{alvarez2}. The global statement above differs in that the pointed curve itself is allowed to vary over $\M_{g,n}$.

The universal injectivity statement recovers the uniqueness theorem for very
strong algebraic-geometric codes over finite fields
\cite[Theorem~2]{marquez4}. It also extends the uniqueness statement to
arbitrary base fields. The effective reconstruction problem considered there is not addressed here.

The Pl\"ucker coordinates admit a theta interpretation. Let
$(X\to S,\sigma_1,\dots,\sigma_n)$ be a family of $n$-pointed smooth projective curves of genus $g$, put $D=\sigma_1+\cdots+\sigma_n$, and set
$Y:=\LS_{X,D,d}$. For a $k$-subset $I\subset\{1,\dots,n\}$, put $D_I:=\sum_{i\in I}\sigma_i$, let $\theta_I$ be the pullback to $Y$ of the
$I$-th Pl\"ucker section, and let $\tau_I:Y\to\Pic_{X,D,g-1}$ be the morphism sending a level structure
$(L,\gamma_1,\dots,\gamma_n)$ to $L(-D_I)$. Let
$\mathcal T$ be the theta line bundle on $\Pic_{X,D,g-1}$, with canonical
section $\theta$, and put
$\mathcal Q:=\Phi_{X,D,d}^{*}\O_{\Gr(k,n)_S}(1)$.

\begin{theorem}[Theorem~\ref{thm:theta-information-relative}]
\label{thm:intro:theta-information}
Assume $n>d>2g-2$. For every $k$-subset $I\subset\{1,\dots,n\}$ there is a canonical isomorphism $\tau_I^{*}\mathcal T\simeq\mathcal Q$ under which
$\tau_I^{*}\theta=\theta_I$. Hence $Z(\theta_I)=\tau_I^{-1}(Z(\theta))$ as closed subschemes of $Y$.
Moreover, for every geometric point $\xi$ of $Y$, the following are equivalent: $\theta_I(\xi)\neq0$, $I$ is an information set of the associated
code $C_\xi$, $H^{0}(X_\xi,L_\xi(-D_{I\xi}))=0$, and $\tau_I(\xi)\notin Z(\theta)$.
\end{theorem}

The theta divisor has previously appeared in the decoding of geometric Goppa codes \cite{henocq-rotillon}. Its role here is different. The Pl\"ucker sections themselves are identified with pullbacks of the canonical theta
section, and their non-vanishing loci are the information-set charts.

The fibers have a direct modular description. For a non-degenerate code $C$, we write $\P_C$ for the projective space of one-dimensional quotients of $C$ and $q_i=[\pr_i|_C]$ for the points determined by the coordinate projections.

\begin{theorem}[Theorem~\ref{thm:fibers-goppa-projective-moduli}]\label{thm:intro:fibers}
Assume $n>d\geq 2g+1$ and let $C\in\Gr(k,n)(\k)$ be a non-degenerate code.
There is a canonical equivalence between the fiber
$\Goppa_{g,n,d}^{-1}(C)$ and the fibered category whose objects over a
$\k$-scheme $T$ are families of $n$-pointed smooth projective curves of
genus $g$ over $T$, together with a non-degenerate closed immersion
$j:X\hookrightarrow\P_C\times T$ of relative degree $d$, sending the
$i$-th marked section to the distinguished point
$q_i\in\P_C$ canonically determined by $C$.
\end{theorem}

The relative form of the preceding description identifies every base change of the Goppa morphism with a classical Hilbert-incidence problem for embedded pointed curves. Incidence schemes of this type were studied explicitly by Perrin using flag Hilbert schemes \cite[\S~1, Example~1.1]{perrin-points}.
Together with Grothendieck's representability of the relative Hilbert functor
\cite{grothendieck-hilbert,nitsure-hilbert}, this yields the following global consequence.

\begin{theorem}[Theorem~\ref{thm:goppa-representable}]
\label{thm:intro:goppa-representable}
Assume $n>d\geq2g+1$. Then $\Goppa_{g,n,d}\colon\LS_{g,n,d}\to\Gr(k,n)$ is schematic and quasi-projective. In particular, $\LS_{g,n,d}$ is a quasi-projective $\k$-scheme.
\end{theorem}

The duality of codes is also compatible with the moduli picture.

\begin{theorem}[Theorem~\ref{prop:selfdual-torsor}]\label{thm:intro:self-duality}
Assume $d\geq 2g+1$ and set $n=2(1-g+d)$. Given a smooth projective curve of genus $g$ with $n$ marked points $(X,D)$, the locus of level structures
determining self-dual codes is the fixed-point subscheme $\fix(\bot)$ of a natural involution $\bot:\LS_{X,D,d}\longrightarrow\LS_{X,D,d}$. The fixed locus $\fix(\bot)$ is a torsor for the fppf topology under the $2$-torsion subgroup scheme $\LS_{X,D,0}[2]$, finite locally free of rank $2^{2d+1}$ over $\k$, and \'etale if
$\operatorname{char}\k\neq2$. In particular, the choice of a $\k$-rational point of $\fix(\bot)$ identifies it with $\LS_{X,D,0}[2]$.
\end{theorem}

Genus zero gives a particularly explicit form of the theory.

\begin{theorem}[Proposition~\ref{prop:LS0nd-explicit} and Theorem~\ref{thm:goppa-immersion-g0}]\label{thm:intro:genus-zero}
Assume $g=0$, $n\geq 3$, and $n>d\geq 1$. Then
$\LS_{0,n,d}\simeq \Ht\times\M_{0,n}$ as schemes. Moreover, if $2\leq d\leq n-3$, the Goppa morphism $\Goppa_{0,n,d}:\LS_{0,n,d}\longrightarrow \Gr(k,n)$ is an immersion and it determines, for each $\lambda\in\Ht(\k)$, an immersion 
\[\Goppa^{\lambda}_{0,n,d}:\M_{0,n}\hookrightarrow\Gr(k,n).\]
\end{theorem}

Over finite fields this is compatible with the classical projective
description of generalized Reed--Solomon codes by distinct points on a normal
rational curve \cite{beelen-grs}. The statement here is scheme-theoretic and,
in the indicated range, includes the immersion of the full parameter space.

\subsection{Geometric interpretation and classical connections}

Theorem~\ref{thm:intro:fibers} identifies the Goppa fibers with classical incidence loci of the type studied by Perrin through flag Hilbert schemes \cite[\S~1, Example~1.1, and \S~2.a]{perrin-points}. Over $\mathbb C$, this
description also has an interpretation in terms of stable maps. Let $\overline{\M}_{g,n}(\P_C,d)$ denote the Kontsevich moduli stack of $n$-pointed genus-$g$ stable maps to $\P_C$ of degree $d$, and let $\ev=(\ev_1,\dots,\ev_n)$ be its total evaluation morphism to $(\P_C)^n$ \cite{fulton-pandharipande}. The fiber of the
Goppa morphism is the open substack of $\ev^{-1}(q_1,\dots,q_n)$ corresponding to smooth non-degenerate closed
immersions.
Evaluation fibers of this form are among the incidence loci that enter Gromov--Witten theory, where incidence conditions are imposed against
the virtual fundamental class. The Goppa fiber retains the actual open locus of smooth embedded curves rather than its virtual count. From this viewpoint, $\LS_{g,n,d}$ is \emph{the gluing} of a family of classical incidence problems as the code $C$ varies.

The genus-zero case makes this relation especially transparent. At
$d=n-2$, the fiber over a code whose distinguished points are in linear general position is canonically $\M_{0,n}$. Its closure in the corresponding Hilbert scheme recovers Kapranov's model of $\overline{\M}_{0,n}$ by rational normal curves through $n$ points in $\P^{n-2}$ \cite[Theorem~0.1]{kapranov-veronese}. Other projective models of pointed rational curves include the log-canonical embedding studied by Keel and Tevelev \cite{keel-tevelev} and the GIT quotients of configurations on $\P^1$ studied by Howard, Millson, Snowden and Vakil \cite{hmsv}. In the range $2\leq d\leq n-3$, Theorem~\ref{thm:intro:genus-zero} instead gives a family of immersions of $\M_{0,n}$ into $\Gr(d+1,n)$. These immersions arise from a different construction and raise the problem of describing their images by equations in Pl\"ucker coordinates.

These examples illustrate a feature that is useful beyond the particular
computations carried out here. A coding-theoretic condition can define a natural locus in a classical moduli space, while a classical moduli problem can appear as a fiber or an image of the Goppa morphism. Self-duality gives
another instance. The coding-theoretic orthogonality condition becomes the fixed locus of an involution on a moduli space of level structures.

\subsection{A broader perspective}\label{subsec:broader}

The rank-one construction suggests a wider moduli-theoretic program. The first extension replaces line bundles by vector bundles. Higher-rank
algebraic-geometric evaluation codes have already been considered, for example in~\cite{savin-vector-bundles}. If $E$ has rank $r$ and the fibers $E|_{p_i}$
are trivialized, evaluation takes values in $\k^{rn}$. Whenever $h^0(X,E)$ is constant 
in families and the image of evaluation is locally free, one obtains a
Grassmannian-valued morphism whose source is the moduli stack of vector bundles with level structure rather than a rigidified relative Picard stack. Note that, when the pointed curve is fixed, this is precisely the setting studied in \cite{alvarez2}.

Principal bundles admit a related formulation. Fix a reductive group $G$ and a
faithful representation $\rho:G\rightarrow\GL(V)$. A principal 
$G$-bundle
determines an associated vector bundle $E$ together with a reduction of structure group, encoded in the compactification theory of Schmitt and
collaborators by tensorial data $(E,\tau)$
\cite{schmitt-decorated-principal,gomez-langer-schmitt-sols}. A level structure on the principal bundle induces one on $E$, so that
evaluation again yields a Grassmannian-valued morphism under suitable cohomological hypotheses, and one may ask how much of the reduction $\tau$ it retains.

The same viewpoint raises arithmetic questions. Where the image of a fixed-type
Goppa morphism is a locally closed locus, one may ask for its dimension,
components, singularities, cycle classes and defining equations, and over a
finite field for its rational points. In favourable situations $\ell$-adic
cohomology and trace formulas become available. Ghorpade and Lachaud provide a precedent in coding theory. They interpret MDS
linear codes as rational points of open strata of Grassmannians
\cite{ghorpade-lachaud-mds}, and study the resulting point-counting problems
using methods based on \'etale cohomology, Lefschetz theorems and the
Grothendieck--Lefschetz trace formula
\cite{ghorpade-lachaud-etale}.

These questions have a cryptographic side. The density of a structured family inside a Grassmannian bears on distinguishability, the geometry of the fibers measures the ambiguity of the hidden representation, and the equations of the image would furnish structural invariants. Such information is already central to the structural attacks discussed above. In all these cases the pattern is the same. An evaluation construction attached to a moduli problem is studied
through the geometry of the resulting morphism to a Grassmannian. In rank one
this is the Goppa morphism studied below, and in higher rank it points towards moduli of vector bundles and of principal bundles. The formalism thereby organizes existing coding-theoretic constructions and suggests new ones, and it produces new maps from familiar moduli spaces to Grassmannians, together with new incidence problems on their images and fibers.

\subsection{Organization}\label{subsec:organization}

Section~\ref{sec:preliminaries} recalls the basic constructions attached to
Goppa codes and level structures. Section~\ref{sec:constructions} constructs
the fibered category of level structures of fixed genus, length and degree and
proves that evaluation is functorial in families. Section~\ref{sec:geometry}
studies the geometry of the resulting moduli spaces and establishes the
monomorphism and immersion properties of the extended Goppa morphism. It also
contains the dimension statements related to the Goppa locus.
Within that section, the pulled-back Pl\"ucker coordinates are
identified with theta sections and their standard affine charts are interpreted
in terms of information sets and systematic generator matrices, which is what
the immersion statement refines.
Section~\ref{sec:fibers} gives the projective interpretation of the fibers,
identifies its relative form with a classical Hilbert-incidence problem and
thereby deduces that the Goppa morphism is schematic and quasi-projective when
$n>d\geq2g+1$. It also proves universal injectivity over an arbitrary base
field.
Section~\ref{sec:duality} treats self-duality from the same geometric
viewpoint. Section~\ref{sec:rational} specializes the theory to genus zero,
where the moduli space and the Goppa morphism admit explicit coordinates and
the fibers are described by rational normal curves.

\subsection{Conventions and notation}

We fix once and for all an arbitrary base field $\k$. We denote by $\Sch$ the
category of $\k$-schemes and by $\Sets$ the category of sets. Given a scheme
$S$, we denote by $\Sch_S$ the category of $S$-schemes. Products of $\k$-schemes are taken over $\Spec\k$ unless otherwise indicated.

Given morphisms of schemes $f:X\rightarrow S$ and $u:T\rightarrow S$, we denote by $X_T$ the pullback $X\times_ST$ and by $f_T:X_T\rightarrow T$ the corresponding morphism.

We write $\mathbb P(\mathcal E)$ for the projective space of
one-dimensional quotients of a locally free sheaf $\mathcal E$. If $C$ is a
$\k$-vector space, we write $\P_C$ for $\P(C)$.

Given natural numbers $k\leq n$, we denote by $\Gr(k,n)$ the Grassmann functor. For a scheme $S$, $\Gr(k,n)(S)$ is the set of locally free
subsheaves $\mathcal E\subset\O_S^n$ of rank $k$ such that $\O_S^n/\mathcal E$ is locally free. The functor $\Gr(k,n)$ is representable by a smooth projective scheme over $\Spec(\mathbb Z)$.

We use the quotient convention for projective bundles and
the subbundle convention for Grassmannians for convenience.

A \emph{linear series} on a projective curve $X$ is a pair $(L,V)$, where
$L$ is an invertible sheaf and $0\neq V\subseteq H^0(X,L)$. Its associated
linear system is $\P(V^\vee)$. The series is \emph{complete} if
$V=H^0(X,L)$. The corresponding linear system is denoted by
$|L|=\P(H^0(X,L)^\vee)$.

We denote by $\mathbb G_m$ the multiplicative group scheme over
$\mathbb Z$, with
$\mathbb G_m(T)=\Gamma(T,\O_T)^*$ for every scheme $T$. 
For a $\k$-scheme $S$ we write $\mathbb G_{m,S}$ for the base change of
$\mathbb G_m$ to $S$. For a commutative $S$-group scheme $A$ we write
$[2]\colon A\to A$ for the squaring morphism and $A[2]:=\ker([2])$ for its
$2$-torsion subgroup scheme. In particular, $\mu_{2,S}:=\mathbb G_{m,S}[2]$.
For an $S$-group scheme $H$, we write
\(
X^{*}(H):=\Hom_{S\text{-grp}}(H,\mathbb G_{m,S})
\)
for its group of characters and 
\(
X_{*}(H):=\Hom_{S\text{-grp}}(\mathbb G_{m,S},H)
\)
for its group of cocharacters.

For the basic theory of categories fibered in groupoids and stacks, we refer
to \cite{laumon,stacks}. We call a morphism of stacks \emph{schematic} if it
is representable by schemes.

For the basic theory of coding theory, we refer to \cite{Tsfasman,stichtenoth,ag-convolutional-survey}. In this article, a $k$-dimensional linear code of length $n$ over a field $\k$ is simply a point of
$\Gr(k,n)(\k)$. In particular, this includes the usual linear block codes
when $\k=\Fq$ and the rational-space description of convolutional codes
when $\k=\Fq(z)$.

\section{Preliminaries on Goppa codes}\label{sec:preliminaries}

\subsection{Basic definitions}
A Goppa structure of genus $g$, length $n$ and degree $d$ is a tuple $\gs:=(X,p_1,\dots,p_n,G)$, where $X$ is a geometrically irreducible smooth projective curve of genus $g$, the points $p_1,\dots,p_n$ are pairwise distinct rational points of $X$, and $G$ is a divisor of degree $d$ whose support contains no $p_i$. This last condition provides, for each $i$, a Zariski open neighbourhood $U_i$ of $p_i$ on which $G$ vanishes, so that 
$\Gamma(U_i,\O_{X}(G))=\Gamma(U_i,\O_{X})$. 
There are canonical isomorphisms 
\begin{equation}\label{eq:canonicaltrivializations}
\gamma^{\mathrm{can}}_i:\O_{X}(G)|_{p_i}\simeq \k.
\end{equation}
Consequently, a Goppa structure $\gs:=(X,p_1,\dots,p_n,G)$ determines a linear map, the \emph{evaluation map}, given by
\begin{equation*}\label{eq:evaluation}
\ev_{\gs}:H^{0}(X,\O_{X}(G))\rightarrow \bigoplus_{i} \O_{X}(G)|_{p_i}\overset{\oplus_i \gamma^{\mathrm{can}}_i}{\simeq} \k^{n}.
\end{equation*}
A Goppa code of genus $g$, length $n$ and degree $d$ is the image $C_{\gs}$ of the evaluation map attached to a Goppa structure $\gs$ of genus $g$, length $n$ and degree $d$ \cite{Gop82}.

This construction admits a natural generalization, the Tsfasman--Vl\u{a}du\c{t} $H$-construction \cite[\S~3.1.1, p.~272]{Tsfasman}, in which level structures of rank one replace Goppa structures. 
\begin{definition}
A level structure of rank one of genus $g$, length $n$ and degree $d$ is a tuple $\gs:=(X,p_1,\dots,p_n,L,\gamma_1,\dots,\gamma_{n})$, where $X$ and $p_1,\dots,p_n$
 are as above, $L$ is an invertible sheaf on $X$ of degree $d$, and $\gamma_1,\dots,\gamma_{n}$ are trivializations 
 $\gamma_i:L|_{p_i}\simeq \k$. 
 \end{definition}
 As before, every level structure $\gs:=(X,p_1,\dots,p_n,L,\gamma_1,\dots,\gamma_{n})$ determines an evaluation map
\begin{equation}\label{eq:evaluationmapcurve}
\ev_{\gs}:H^{0}(X,L)\rightarrow \bigoplus_{i} L|_{p_i}\overset{\oplus \gamma_i}{\simeq} \k^{n}.
\end{equation}
We also call the image of the evaluation map attached to a level structure a Goppa code.

 \subsection{Parameters}\label{subsec:parameters}

Let us denote by $D$ the effective divisor $p_1+\cdots+p_n$ and let
$\gs=(X,p_1,\dots,p_n,L,\gamma_1,\dots,\gamma_{n})$ be a level structure of rank
one. The exact sequence $0\to L(-D)\to L\to\bigoplus_{i} L|_{p_i}\to 0$ gives in
cohomology
 \begin{equation}\label{eq:exactsequencecohomology}
 \xymatrix{
 0\ar[r] & H^{0}(X,L(-D))\ar[r] &H^{0}(X, L) \ar[r] & \oplus_{i=1}^{n} L|_{p_i}  & \\
 \ar[r]& H^{1}(X,L(-D))\ar[r] & H^{1}(X,L)\ar[r] & 0
 }
 \end{equation}
whose middle map, composed with the trivializations, is exactly the evaluation
map \eqref{eq:evaluationmapcurve} attached to $\gs$. If $n>d$, then
$\deg L(-D)=d-n<0$, so $H^{0}(X,L(-D))=0$ and $\ev_{\gs}$ is injective. If
moreover $d>2g-2$, then $\deg(\omega_{X}\otimes_{\O_X} L^{-1})=2g-2-d<0$, so
$H^{1}(X,L)=0$ by Serre duality and the Riemann--Roch theorem gives
\begin{equation}\label{eq:dimension-code}
\dim_{\k} C_{\gs}=\dim_{\k}H^{0}(X,L)=1-g+d=:k.
\end{equation}
We call $C_{\gs}$ a strong Goppa code, and $\gs$ a strong level structure, if
$n>d>2g-2$. For $g=0$ this condition reduces to $n>d$.

Recall that a linear code $C\subseteq\k^{n}$ is \emph{non-degenerate} if
$\pr_{i}|_{C}\neq0$ for every $i=1,\dots,n$, where $\pr_{i}:\k^{n}\to\k$ is the
projection onto the $i$-th coordinate. If $d>2g-1$, every Goppa code $C_{\gs}$ is
non-degenerate. Indeed, $\deg L(-p_{i})=d-1>2g-2$ for every $i$, so the
Riemann--Roch theorem gives $h^{0}(X,L(-p_{i}))=d-g=h^{0}(X,L)-1$. A section of
$L$ outside $H^{0}(X,L(-p_{i}))$ does not vanish at $p_{i}$, and its image under
$\ev_{\gs}$ has nonzero $i$-th coordinate.

 \subsection{Equivalent level structures}
 The set of level structures carries a natural equivalence relation.
 \begin{definition}\label{def:equivalentlevelstructures}
 Let $g,n,d$ be natural numbers.
Two level structures $\gs$ and $\gs'$ of genus $g$, length $n$ and degree $d$ are equivalent if there exist an isomorphism of curves $f:X\xrightarrow{\sim} X'$ and an isomorphism of sheaves $\phi:L\xrightarrow{\sim} f^{*}L'$ such that 
\begin{enumerate}
\item $f(p_i)=p'_i$ for all $i=1,\dots,n$.
\item The square
 \begin{equation}\label{eq:equivalentlevelstructures}
 \xymatrix{
 L|_{p_i}\ar[d]_{\phi|_{p_i}}^{\cong}\ar[rr]^{\gamma_i}_{\cong}&& \k\\
 f^{*}L'|_{p_i}\ar[rr]^{\cong}_{\textrm{\tiny canonical}}&& L'|_{p'_i}\ar[u]_{\gamma'_i}^{\cong}
 }
 \end{equation}
 commutes for every $i=1,\dots,n$.
 \end{enumerate}

 \end{definition}
 \begin{lemma}\label{lm:equivalentstructures}
 If two level structures $\gs,\gs'$ are equivalent, then $C_{\gs}=C_{\gs'}$.
 \end{lemma}
 \begin{proof}
 The commutativity of the square \eqref{eq:equivalentlevelstructures}, together with the canonical isomorphism $H^{0}(X,f^{*}L')\simeq H^{0}(X',L')$, implies that the diagram
  \begin{equation}
 \xymatrix{
 H^{0}(X,L)\ar[d]^{\cong}_{H^{0}(\phi)}\ar[r] \ar@/_{18mm}/[dd]_{\widehat{H^{0}(\phi)}}& L|_{p_i}\ar[d]_{\phi|_{p_i}}^{\cong}\ar[rr]^{\gamma_i}_{\cong}&& \k\\
  H^{0}(X,f^{*}L') \ar[r] & f^{*}L'|_{p_i}\ar[rr]^{\cong}_{\textrm{\tiny canonical}} \ar[d]_{\textrm{\tiny canonical}}^{\cong} && L'|_{p'_i}\ar[u]_{\gamma'_i}^{\cong}\\
  H^{0}(X',L') \ar[u]_{\cong}^{\textrm{\tiny canonical}} \ar[r] & L'|_{p'_i}\ar@{=}[rru]
 }
 \end{equation}
 commutes for every $i=1,\dots,n$. It follows that the triangle
 $$
 \xymatrix{
  H^{0}(X,L)\ar[d]^{\cong}_{\widehat{H^{0}(\phi)}} \ar@{^(->}[r]&\k^{n}\\
  H^{0}(X',L')   \ar@{^(->}[ru]&
  }
 $$
 commutes, and therefore $C_{\gs}=C_{\gs'}$.
 \end{proof}
 \begin{remark}
 Suppose that $X'=X$ and $f=\mathrm{id}_X$ in Definition~\ref{def:equivalentlevelstructures}. The equivalence relation then simplifies, and two level structures $\gs,\gs'$ over the pointed curve $(X,p_1,\dots,p_n)$ are equivalent if there is an isomorphism $\phi:L\xrightarrow{\sim} L'$ for which the triangle
 $$
 \xymatrix{
 L|_{p_i}\ar[d]_{\phi|_{p_i}}\ar[r]^{\gamma_i} & \k\\
 L'|_{p_i}\ar[ru]_{\gamma'_i}&
 }
 $$
 commutes for every $i=1,\dots,n$. This relation agrees with the classical one when $L=\O_{X}(G)$ and $L'=\O_{X}(G')$, with $\Supp(G)$ and $\Supp(G')$ containing no $p_i$. Indeed, the isomorphism $\phi$ is multiplication by a rational function $h\in\k(X)$ with $G=G'+(h)$. Restricting the triangle above to an open subset $U$ containing $D$ yields $s(p_i)h(p_i)=s(p_i)$ for every $s\in \Gamma(U,\O_{X}(G))$, that is, $h(p_i)=1$ for all $i=1,\dots,n$. The converse follows from the same computation. The level structures induced by $G$ and $G'$ are equivalent if and only if there is a rational function $h\in \k(X)$ with $G=G'+(h)$ and $h(p_i)=1$ for all $i=1,\dots,n$, which is the usual equivalence relation for classical Goppa structures \cite{stichtenoth,xing-equal,munuera-pellikaan}.
 \end{remark}
We use the following notation.
\begin{enumerate}
\item $\LS_{X,D,d}$ is the set of equivalence classes of level structures over a pointed curve $(X,p_1,\dots,p_n)$ of genus $g$, length $n$ and degree $d$. Here $D=p_1+\cdots+p_n$.
\item $\LS_{g,n,d}$ is the set of equivalence classes of level structures of genus $g$, length $n$ and degree $d$.
\end{enumerate}

\subsection{The group of level structures on a fixed pointed curve} \label{sec:aremark}
Raynaud \cite{raynaud} developed the theory of rank-one level structures in detail in order to study the Picard functor of families of schemes parametrized by discrete valuation rings. We do not reproduce the general results here. The results we need appear in the course of the paper. We record, however, that for a fixed pointed curve $(X,p_1,\dots,p_n)$ the set $\LS_{X,D}:=\coprod_{d\in\mathbb{Z}}\LS_{X,D,d}$ of level structures of arbitrary degree carries a group structure. Indeed, let $(L_1,\gamma^{1}_1,\dots,\gamma^{1}_{n})$ and $(L_2,\gamma^{2}_1,\dots,\gamma^{2}_{n})$ be level structures in $\LS_{X,D}$. The invertible sheaf 
$L:=L_1\otimes_{\O_X} L_2$
together with the isomorphisms $\gamma_i$ defined by
$$
(L_{1}\otimes_{\O_X} L_2)|_{p_i}\simeq L_1|_{p_i}\otimes_{\k} L_{2}|_{p_i} \simeq \k\otimes_{\k}\k\simeq\k
$$
determines a new level structure $(L,\gamma_1,\dots,\gamma_n)\in \LS_{X,D}$. The neutral element of $\LS_{X,D}$, denoted $\mathbf 1$, is $\O_{X}$, equipped with the trivializations at $p_1,\dots,p_n$ inverse to the canonical maps $\k\rightarrow \O_{X,p_{i}}/\mathfrak{m}_{p_i}$. Every level structure $(L,\gamma_1,\dots,\gamma_n)\in \LS_{X,D}$ admits an inverse $(N,\eta_1,\dots,\eta_n)$ for this operation, namely 
$N=L^{-1}$ 
with $\eta_i$ induced by $\gamma_{i}^{-1}:\k\simeq L|_{p_i}$ through the chain of isomorphisms
$$L^{-1}|_{p_i}\simeq L^{-1}|_{p_i} \otimes_{\k} \k\simeq L^{-1}|_{p_i}\otimes_{\k} L|_{p_i}\simeq (L^{-1} \otimes_{\O_X} L)|_{p_i}\simeq \k.$$

  \subsection{The dual code}\label{subsec:dual-code}
 
 Let us consider the pairing $\langle{-},{-}\rangle:\k^{n}\times \k^n\rightarrow \k$ defined by $(w_1,w_2)\mapsto w_1 \cdot w_{2}^{T}$. Every code $C\subset\k^n$ has a dual code with respect to this pairing, namely 
 $$
 C^{\bot}:=\{w'\in\k^{n}\mid \langle w',w\rangle=0 \text{ for all }w\in C\}.
 $$
 If $C=C_{\gs}$ for some level structure $\gs\in \LS_{X,D,d}$, then $C^{\bot}$ admits the following algebraic-geometric description. From $\gs$ and the two middle maps of \eqref{eq:exactsequencecohomology} we construct the diagram
 \begin{equation}
 \xymatrix{
 &\k^{n}\ar@{..>}[rd]&\\
H^{0}(X, L) \ar[r]\ar@{..>}[ur]^{\ev_{\gs}} & \oplus_{i=1}^{n} L|_{p_i} \ar[r]  \ar[u]^{\cong}_{\oplus \gamma_{i}}& H^{1}(X,L(-D))
 }
 \end{equation}
 in which the dotted arrows denote the indicated compositions. The composition of the two dotted arrows vanishes by the exactness of \eqref{eq:exactsequencecohomology}. Consequently, $C_{\gs}^{\bot}$ is the column space of the matrix of the second dotted arrow with respect to a basis of $H^{1}(X,L(-D))$ or, equivalently, the row space of its transpose.
 
 We now describe the transpose of $\k^n\rightarrow H^{1}(X,L(-D))$ precisely. Let us consider the invertible sheaf $\omega_{X}(p):=\omega_{X}\otimes_{\O_{X}}\O_{X}(p)$, where $\omega_{X}$ is the sheaf of differentials of $X$ and $p\in X$ is a rational point. For every open subset $U\subset X$ one has 
 $$
 \omega_{X}(p)(U)=\{w\in \Omega^{1}_{\k(X)/\k}| \ (w)|_{U}+p|_{U} \geq 0\}.
 $$
Suppose that $p\in U$, and let $t_{p}$ be a generator of the maximal ideal $\mathfrak{m}_{p}$. After shrinking $U$, every $w\in \omega_{X}(p)(U)$ can be written as $w=f\,dt_{p}$ with $f\in\k(X)$ satisfying $(f)|_{U}+p|_{U}\geq 0$. In particular, $w$ is either regular at $p$ or has a pole of order one at $p$. Hence $\Res_{p}(h w)=h(p)\Res_{p}(w)$ for every $h\in \O_{X}(U)$, that is, the map
$
 \omega_{X}(p)(U)\rightarrow \O_{p}(U)=\k, \ w\mapsto \Res_{p}(w),
$
is a morphism of $\O_{X}(U)$-modules. It induces a morphism of sheaves $ \omega_{X}(p)\rightarrow i_{p*}\O_{p}$, which by adjunction determines an isomorphism
$
\Res_{p}: \omega_{X}(p)|_{p}\simeq \k.
$
Combining this isomorphism with the canonical trivializations $\O_{X}(q)|_{p}\simeq\k$ for $p\neq q$, we obtain canonical trivializations
$\Res_{p_{i}}: \omega_{X}(D)|_{p_{i}}\simeq \k$,
 where $D=p_1+\cdots+p_n$. Let us set $M:=\omega_{X}(D)\otimes_{\O_{X}}L^{-1}$. The sheaf $M$ inherits trivializations at the points $p_i$, namely 
\begin{equation*}\label{eq:residue_infinite}
\eta_{i}:M|_{p_{i}}= \omega_{X}(D)|_{p_{i}}\otimes_{\k} L^{-1}|_{p_{i}}\simeq\k\otimes_{\k}\k\simeq \k
\end{equation*} 
 where the first isomorphism is $\Res_{p_i}\otimes \gamma_{i}^{-1}$. For the notation $\gamma_i^{-1}$, see \S~\ref{sec:aremark}. The long exact sequence in cohomology attached to the short exact sequence
$$
\xymatrix{
0 \ar[r] & M(-D) \ar[r]& M\ar[r] & \oplus_{i=1}^{n}M|_{p_i}\ar[r] & 0
}
$$ 
yields, together with Serre duality, a linear map
$$
\Res_{\gs}:H^{1}(X,L(-D))^{\vee}\simeq H^{0}(X,\omega_{X}(D)\otimes_{\O_{X}}L^{-1})\rightarrow \bigoplus_{i=1}^{n}M|_{p_i}\simeq \k^{n}.
$$
This map is the transpose of $\k^n\rightarrow H^{1}(X,L(-D))$ (see \cite[Chap.~3.1, p.~289]{Tsfasman}), and therefore 
$\im(\Res_{\gs})=C_{\gs}^{\bot}$.

\begin{remark}
Throughout the article, a \emph{parity-check matrix} for a code $C\subset\k^n$ of dimension $k$ is an $n\times(n-k)$ matrix $H$ over $\k$ whose left kernel, that is, the set of row vectors $v\in\k^n$ with $vH=0$, equals $C$. Equivalently, $C=\ker(H^{T})$, where $H^{T}$ is regarded as a map $\k^n\to\k^{n-k}$ acting on column vectors.
\end{remark}

\section{The space of level structures and the Goppa morphism, I: constructions}\label{sec:constructions}

Let $n>d>2g-2$ and set $k:=1-g+d$. By \eqref{eq:dimension-code} and Lemma~\ref{lm:equivalentstructures}, there are Goppa maps
\begin{equation}
\begin{split}
\Goppa_{X,D,d}: \LS_{X,D,d}&\rightarrow \Gr(k,n)(\k)\\
\Goppa_{g,n,d}: \LS_{g,n,d}&\rightarrow \Gr(k,n)(\k)
\end{split}
\end{equation}
assigning to each level structure $\gs$ its Goppa code $C_{\gs}$.

In this section we endow $\LS_{X,D,d}$ and $\LS_{g,n,d}$ with a natural algebraic-geometric structure for which the Goppa maps become morphisms of schemes and of stacks, respectively. We show that the \emph{extended} Goppa morphism $\Phi_{X,D,d}: \LS_{X,D,d}\rightarrow \Gr(k,n)_{S}$  is a monomorphism in the range of Lemma~\ref{alvarez-porras}, and $\Phi_{g,n,d}\colon\LS_{g,n,d}\to\Gr(k,n)\times\M_{g,n}$ is an immersion under the stronger numerical hypothesis $n\ge d+g+2$ (Theorem~\ref{th:immersionstacks}). In the immersion range, and after fixing a family of $n$-pointed curves, the corresponding locus of Goppa codes is then a locally closed subscheme of the relative Grassmannian. Projecting this locally closed locus along $\Gr(k,n)\times\M_{g,n}\to\Gr(k,n)$ gives the image of $\Goppa_{g,n,d}$ in $\Gr(k,n)$, which is constructible in that range. Whether it is locally closed is a separate question, which we do not address here. Its rational points are exactly the Goppa codes of type $(g,n,d)$. This raises the problem of describing that image explicitly and of determining its defining equations in suitable affine charts.

\subsection{Category of level structures \texorpdfstring{$\LS_{g,n,d}$}{LS(g,n,d)} and Goppa functor}

\begin{definition}
Let $S$ be an arbitrary scheme. A family of $n$-pointed smooth projective curves of genus $g$ parametrized by $S$ consists of a proper, flat, finitely presented morphism of schemes $\pi: X\rightarrow S$ whose geometric fibers are smooth curves of genus $g$, together with $n$ sections $\sigma_1,\dots,\sigma_n: S\rightarrow X$ such that $\sigma_{i}(s)\neq \sigma_{j}(s)$ for every $s\in S$ and all $i\neq j$. 
\end{definition}
Every section $\sigma_i$ defines a divisor in $X$, which we denote by $p_i$. We write $D:=\sum_{i=1}^{n}p_i$ and denote by $\O_{D}$ the structure sheaf of $D$. For each $S$-scheme $u:T\rightarrow S$, write $\sigma_{iT}:T\rightarrow X_T$ for the base-changed section and $p_{iT}$ for the corresponding divisor. Then $D_T:=\sum_{i=1}^{n}p_{iT}=D\times_S T$, and, denoting its structure sheaf by $\O_{D_T}$, one has
$$i_{D_{T}*}\O_{D_T}=\bigoplus_{i=1}^{n} \sigma_{iT*}\O_{T},$$
where $i_{D_{T}}$ denotes the inclusion $D_{T}\subset X_{T}$.

\begin{lemma}\label{lm:trivializations-surjection}
Let $\pi: X\rightarrow S$ and $\sigma_1,\dots,\sigma_n: S\rightarrow X$ be as above, let $L$ be an invertible sheaf on $X$, and let $\sigma _{i}:S\rightarrow X$ be the section corresponding to $p_{i}$. Tuples of isomorphisms $(\gamma_i:\sigma_{i}^{*}L\simeq \O_{S})_{i=1}^{n}$ correspond bijectively to surjective morphisms $\gamma:L\rightarrow i_{D*}(\O_{D})$.
\end{lemma}
\begin{proof}
For each $i$, the adjunction between $\sigma_{i}^{*}$ and $\sigma_{i*}$ gives
$
\Hom_{\O_{S}}(\sigma _{i}^{*}L, \O_{S})=\Hom_{\O_{X}}(L,\sigma _{i*}\O_{S}),
$
the morphism attached to $f:\sigma _{i}^{*}L\rightarrow \O_{S}$ being the composition $L\rightarrow \sigma _{i*}\sigma _{i}^{*}L\rightarrow \sigma_{i*}\O_{S}$. Since $\sigma_i$ is a closed immersion, and in particular affine, a morphism $L\rightarrow \sigma _{i*}\O_{S}$ is surjective if and only if $\sigma _{i*}\sigma _{i}^{*}L\rightarrow \sigma _{i*}\O_{S}$ is surjective. Since $\sigma _{i}^{*}\sigma _{i*}\mathcal{F}=\mathcal{F}$ for every locally free sheaf $\mathcal{F}$ on $S$, the latter is surjective if and only if $\sigma_{i}^{*}L\rightarrow \O_{S}$ is surjective. Finally, a surjection between locally free sheaves of the same rank is an isomorphism, so $\sigma _{i}^{*}L\rightarrow \O_{S}$ is surjective if and only if it is an isomorphism. As $i_{D*}(\O_{D})=\bigoplus_{i=1}^{n}\sigma_{i*}\O_{S}$, the assertion follows by taking the components of $\gamma$.
\end{proof}

\begin{definition}
A family of level structures of genus $g$, length $n$ and degree $d$ is a tuple
\begin{equation}\label{def:level}
\gs=(\pi: X\rightarrow S, \sigma_1,\dots,\sigma_n: S\rightarrow X, L,\gamma: L\rightarrow i_{D*}(\O_{D}))
\end{equation}
where $(\pi: X\rightarrow S, \sigma_1,\dots,\sigma_n: S\rightarrow X)$ is a family of $n$-pointed curves of genus $g$, $L$ is an invertible sheaf of degree $d$ and $\gamma: L\rightarrow i_{D*}(\O_{D})$ is a surjective morphism (see Lemma~\ref{lm:trivializations-surjection}).
\end{definition}

\begin{definition}\label{def:categorylevel}
Let 
\begin{equation}
\begin{split}
\gs_{1}&=(\pi_1: X_1\rightarrow S_1, \sigma^1_1,\dots,\sigma^1_n: S_1\rightarrow X_1, L_1,\gamma_1:L_1\rightarrow i_{D_1* }((\O_{D_1})))\\
\gs_{2}&=(\pi_2: X_2\rightarrow S_2, \sigma^2_1,\dots,\sigma^2_n: S_2\rightarrow X_2, L_2,\gamma_2:L_2\rightarrow i_{D_2* }((\O_{D_2})))
\end{split}
\end{equation}
be two level structures. A morphism from $\gs_1$ to $\gs_2$ is a triple $\theta=(u,f,\phi)$ where
the pair $(u,f)$ determines a cartesian diagram
\begin{equation}\label{eq:morphismofcurves}
\xymatrix{
X_1\ar[r]^{f}\ar[d]^{\pi_1} & X_2\ar[d]^{\pi_2} \\
S_1 \ar[r]^{u} & S_2
}
\end{equation}
such that $f\circ\sigma_{i}^{1}=\sigma_{i}^{2}\circ u$ for all $i=1,\dots,n$
and $\phi$ is an isomorphism of sheaves, $\phi: L_{1}\simeq f^{*}L_2$, such that the square
\begin{equation}%
\xymatrix{
\sigma_{i}^{1*}L_1\ar[rr]^{\gamma_{1i}} \ar[d]_{\sigma_{i}^{1*}\phi} & & \O_{S_1}  \\
\sigma_{i}^{1*}f^{*}L_2 \ar[r]^{\simeq}_{\mathrm{canonical}}& u^{*}\sigma_{i}^{2*}L_2\ar[r]^{\simeq}_{u^{*}(\gamma_{2i})}  & u^{*}\O_{S_2}\ar[u]^{\simeq}_{\mathrm{canonical}}
}
\end{equation}
is commutative.
\end{definition}

\begin{definition}
Let $g,n,d$ be natural numbers. The category of level structures of genus $g$, length $n$ and degree $d$, $\LS_{g,n,d}$, is the category fibered in groupoids over $\Sch$ whose objects are families of level structures $\gs$ and where morphisms $\gs_1\rightarrow \gs_2$ are triples $(u,f,\phi)$ as above. The structure morphism $\LS_{g,n,d}\rightarrow \Sch$ is the forgetful map that sends $\gs$ to the base scheme $S$.
\end{definition}

We now show that the Goppa map $\Goppa_{g,n,d}$ of the previous section admits a functorial generalization (Theorem~\ref{prop:goppa-morphism-functorial}). Precisely, we prove that for $n>d>2g-2$ and $k=1-g+d$ there is a morphism of categories fibered in groupoids over $\Sch$,
$\Goppa_{g,n,d}: \LS_{g,n,d}\rightarrow \Gr(k,n)$,  
whose restriction to $\k$-valued points is the Goppa map. Defining this morphism amounts to specifying the object attached to each family of level structures $\gs$ as in \eqref{def:level}. Such a family determines an exact sequence
$$
\xymatrix{
0\ar[r] &L(-D)\simeq \ker(\gamma)\ar[r] & L \ar[r]^-{\gamma} & i_{D*}(\O_{D})\ar[r] & 0.
}
$$
Assume $n>d>2g-2$. The Riemann--Roch theorem and the non-noetherian form of the theorem on cohomology and base change give $R^{1}\pi_{*}L=0$ and show that $\pi_{*}L$ and $R^{1}\pi_{*}\ker(\gamma)$ are locally free of ranks $1-g+d$ and $n-(1-g+d)=n-1+g-d$, respectively. Here the hypothesis $d>2g-2$ yields $R^{1}\pi_*L=0$. On the other hand, $\deg L(-D)=d-n<0$ on every fiber yields $\pi_*L(-D)=\pi_*\ker(\gamma)=0$. Since $i_{D*}(\O_{D})=\bigoplus_{i=1}^{n} \sigma_{i*}\O_{S}$, one has $\pi_{*}(i_{D*}(\O_{D}))=\O_{S}^{n}$, and the long exact sequence obtained by applying $\pi_*$ to $0\to L(-D)\to L\xrightarrow{\gamma} i_{D*}(\O_D)\to 0$ reads
$0\to\pi_*L\xrightarrow{\pi_*\gamma}\O_S^n\to R^1\pi_*L(-D)\to 0$,
so that $\pi_*\gamma$ is a subbundle inclusion with locally free cokernel $R^1\pi_*L(-D)$. Consequently the subsheaf
$$
C_{\gs}:=\im(\pi_{*}L \lhook\joinrel\xrightarrow{\pi_{*}\gamma} \O_{S}^{n})
$$
is a subbundle of $\O_S^n$ of rank $k$ with locally free quotient, that is, an $S$-point $C_{\gs}\in\Gr(k,n)(S)$.
\begin{theorem}\label{prop:goppa-morphism-functorial}
If $n>d>2g-2$,
the maps $\gs\mapsto C_{\gs}=\im(\pi_{*}L \lhook\joinrel\xrightarrow{\pi_{*}\gamma} \O_{S}^{n})$ and $(u,f,\phi)\mapsto u$ define a morphism of categories fibered in groupoids 
\begin{equation}\label{eq:goppamap}
\Goppa_{g,n,d}: \LS_{g,n,d}\rightarrow \Gr(k,n),
\end{equation}
which we call the Goppa morphism.
\end{theorem}
\begin{remark}
Note that every scheme, and in particular $\Gr(k,n)$, determines canonically a category fibered in groupoids over $\Sch$. Its objects are the pairs $(S,v: S\rightarrow \Gr(k,n))$, and its morphisms $(S,v)\rightarrow (S',v')$ are the morphisms of schemes $u:S\rightarrow S'$ with $v=v'\circ u$.
\end{remark}
\begin{proof}
Let $(u,f,\phi):\gs_1\rightarrow \gs_2$ be a morphism, with
$u:S_1\rightarrow S_2$ the underlying morphism of base schemes. We must check
that $C_{\gs_1}=u^{*}C_{\gs_2}$ as subsheaves of
$\O_{S_1}^{n}$. The square \eqref{eq:morphismofcurves} is cartesian, so the
theorem on cohomology and base change \cite[Tag~0D4E]{stacks} applies to it.
Together with the isomorphism $\phi:L_1\simeq f^{*}L_2$ and the identities
$f\circ\sigma^{1}_i=\sigma^{2}_i\circ u$, it identifies the exact sequence
$$
0 \rightarrow \pi_{1*} L_1\rightarrow \O_{S_1}^{n} \rightarrow R^{1} \pi_{1*} L_1(-D_1) \rightarrow 0
$$
with the pullback along $u$ of the corresponding sequence for $\gs_2$,
compatibly with the inclusions into $\O^{n}$. Consequently, the two subbundles
$C_{\gs_1}$ and $u^{*}C_{\gs_2}$ of $\O_{S_1}^{n}$ coincide.
\end{proof}

\subsection{The induced category \texorpdfstring{$\LS_{X,D,d}$}{LS(X,D,d)} and Goppa functor}

Consider the forgetful functor $\Theta:\LS_{g,n,d}\rightarrow \M_{g,n}$. Let $(\pi:X\rightarrow S,\sigma_1,\dots,\sigma_n:S\rightarrow X)$ be a family of $n$-pointed smooth projective curves of genus $g$ parametrized by $S$, and let $D$ be the divisor defined by $\sigma_1,\dots,\sigma_n$. This family determines a functor $\Sch_{S}\rightarrow \M_{g,n}$, and we define the category fibered in groupoids
\begin{equation*}
\LS_{X,D,d}:=\LS_{g,n,d}\times_{\M_{g,n}} \Sch_{S}.
\end{equation*}

\begin{lemma}\label{lm:equiv-cat}
Let $(\pi:X\rightarrow S,\sigma_1,\dots,\sigma_n:S\rightarrow X)$ be a family of $n$-pointed smooth projective curves of genus $g$ parametrized by a scheme $S$ and $D$ the divisor defined by $\sigma_1,\dots,\sigma_n$. Then the category $\LS_{X,D,d}$ is equivalent to the category whose objects are level structures $\gs$ where the underlying family of $n$-pointed smooth projective curves is $(\pi:X_T\rightarrow T,\sigma_{1T},\dots,\sigma_{nT}:T\rightarrow X_T)$, and whose morphisms $\gs_1\rightarrow\gs_2$ are isomorphisms $\phi:L_1\simeq L_2$ such that $(\psi:T_1\rightarrow T_2,\mathrm{id}_X\times \psi,\phi)$ is a morphism in $\LS_{g,n,d}$.
\end{lemma}
\begin{proof}
Recall that objects of $\LS_{X,D,d}$ are triples $(\gs,T\rightarrow S,\alpha)$ where $\gs=(\pi: Y\rightarrow T, \beta_1,\dots,\beta_n: T\rightarrow Y, N,\gamma: N\rightarrow i_{E* }((\O_{E})))$ is an object of $\LS_{g,n,d}$ and $\alpha$ is an isomorphism of pointed curves fitting into the diagram
$$
\xymatrix{
Y\ar[rr]^{\alpha} \ar[rd] & & X_{T}\ar[ld]\\
&T\ar@/_/[ur]_{\alpha_{iT}} \ar@/^/[ul]^{\beta_i}&
}
$$
On the other hand, morphisms $(\gs_1,T_1\rightarrow S,\alpha_1)\rightarrow (\gs_2,T_2\rightarrow S,\alpha_2)$ are pairs $(\theta:\gs_1\rightarrow \gs_2,\psi:T_1\rightarrow T_2)$ such that the exterior square
$$
\xymatrix{
Y_1\ar[rr]^{\alpha_1} \ar[rd] \ar[ddd]_{\Theta(\theta)}& & X_{T_1}\ar[ld]\ar[ddd]^{\mathrm{id}_{X}\times\psi}    \\
&T_1\ar@/^/[ur]^{\alpha^1_{iT}} \ar@/_/[ul]_{\beta^1_i} \ar[d]_{\psi}&  \\
& T_2 \ar@/_/[dr]_{\alpha^2_{iT}} \ar@/^/[dl]^{\beta^2_i} &\\
Y_2\ar[rr]_{\alpha_2} \ar[ru] & & X_{T_2}\ar[lu]
}
$$
commutes. In particular, $\psi=u$, where $u$ is the first component of $\theta$ (see \eqref{eq:morphismofcurves}). To each object $(\gs,T\rightarrow S,\alpha)$ we associate the level structure $\alpha_{*}\gs$ over the pointed curve $(X_{T}\rightarrow T,\alpha_{1T},\dots,\alpha_{nT}: T\rightarrow X_{T})$ determined by the pair
$$
(\alpha^{-1*}L,\alpha^{-1*}\gamma:\alpha^{-1*}L\rightarrow \alpha^{-1*}i_{E* }((\O_{E})))
$$
This assignment is well defined because $\alpha^{-1*}i_{E* }((\O_{E}))$ is canonically isomorphic to $i_{D_{T}*}(\O_{D_{T}})$ and $\alpha^{-1*}\gamma$ remains surjective. Finally, to each morphism $(\theta,\psi): (\gs_1,T_1\rightarrow S,\alpha_1)\rightarrow (\gs_2,T_2\rightarrow S,\alpha_2)$ we associate the morphism $\alpha_{1*}\gs_1\rightarrow \alpha_{2*}\gs_2$ given by $\psi, \mathrm{id}_{X}\times\psi$ and the isomorphism of invertible sheaves
$$
\alpha_{1}^{-1*}L\overset{\alpha_{1}^{-1*}\phi}{\simeq}\alpha_{1}^{-1*}\Theta(\theta)^{*}L_2=(\mathrm{id}_{X}\times\psi)^{*}\alpha_{2}^{-1*}L_2.
$$
These assignments are compatible with composition and with base change, and the
construction $(\gs,T\to S,\alpha)\mapsto\alpha_*\gs$ has the evident inverse
sending a level structure over $(X_T,\sigma_{iT})$ back to the same datum together
with $\alpha=\mathrm{id}$. Hence they define an equivalence of categories.
\end{proof}
From \eqref{eq:goppamap} we get the corresponding Goppa morphism
\begin{equation*}
\Goppa_{X,D,d}:\LS_{X,D,d}\rightarrow \Gr(k,n),
\end{equation*}
valid whenever $n>d>2g-2$.

\begin{remark}
The sets of equivalence classes of level structures introduced in \S~\ref{sec:preliminaries} are canonically identified with the sets of isomorphism classes of $\k$-valued points of the corresponding fibered categories. The sets $\LS_{g,n,d}$ and $\LS_{X,D,d}$ defined there are the sets of isomorphism classes of objects of $\LS_{g,n,d}(\k)$ and of $\LS_{X,D,d}(\k)$, respectively. We retain the lighter notation whenever no confusion can arise.
\end{remark}

\section{The space of level structures and the Goppa morphism, II: geometry}\label{sec:geometry}

\subsection{The category \texorpdfstring{$\LS_{X,D,d}$}{LS(X,D,d)}. Structure}\label{sec:relativestructure}

Let $(\pi:X\rightarrow S,\sigma_1,\dots,\sigma_n:S\rightarrow X)$ be a family of $n$-pointed smooth projective curves of genus $g$ over a scheme $S$, and let $D$ be the divisor defined by the sections $\sigma_1,\dots,\sigma_n$.
We denote by $\Pic_{X/S}$ the relative Picard functor, that is, the fppf sheaf on $\Sch_S$ associated to the presheaf which sends an $S$-scheme $u:T\to S$ to the set of isomorphism classes of invertible sheaves on $X_T$, and by $\Pic_{X,D,d}\subset\Pic_{X/S}$ its part of relative degree $d$.
Likewise, we denote by $\LS^{\#}_{X,D,d}$ the relative moduli functor of level structures of degree $d$, that is, the functor that sends an $S$-scheme $u:T\to S$ to the set of isomorphism classes of tuples
$(X_T\rightarrow T,\sigma_{1T},\dots,\sigma_{nT},L,\gamma)$,
where $L$ is an invertible sheaf on $X_T$ of relative degree $d$ and
$\gamma:L\longrightarrow i_{D_T*}(\O_{D_T})$
is a surjective morphism.

\begin{remark}\label{rmk:LS-sheaf}
Every automorphism in $\LS_{X,D,d}(T)$ is trivial. Indeed, an automorphism of $(L,\gamma)$ is an automorphism $\phi$ of $L$ with $\gamma\circ\phi=\gamma$. Since the fibers of $X_T\to T$ are geometrically connected, $\phi$ is multiplication by a unit $u\in\Gamma(T,\O_{T})^{*}$, and restricting the identity $\gamma\circ\phi=\gamma$
along the section $\sigma_{1T}$ gives $u=1$. Consequently $\LS_{X,D,d}$ is
equivalent to the category fibered in sets associated with $\LS^{\#}_{X,D,d}$. In particular, it is a sheaf of sets on $\Sch_S$. We pass freely between these two viewpoints and use only the notation $\LS_{X,D,d}$.
\end{remark}

Both $\LS_{X,D,d}$ and $\Pic_{X,D,d}$ are representable by finitely
presented, formally smooth algebraic $S$-spaces. See
\cite[Th\'eor\`eme~2.3.1 and Corollaire~2.3.2]{raynaud} and \cite[Theorem~7.3]{Artin}, respectively. Moreover, $\Pic_{X,D,d}$ is a proper $S$-scheme. Indeed, $\Pic_{X,D,0}$ is the relative Jacobian and is proper over $S$ \cite[Proposition~9.5.20]{kleiman-picard}, while tensoring with $\O_X(d\sigma_1)$ identifies $\Pic_{X,D,0}$ with $\Pic_{X,D,d}$.

\begin{lemma}\label{lm:preliminariesraynaud}
The functors $\LS_{X,D,d}$ and $\Pic_{X,D,d}$ satisfy the following properties.
\begin{enumerate}
\item $\LS_{X,D,d}$ is acted on by the group $S$-scheme $H_S:=(\mathbb{G}_{m,S})^{n}/\mathbb{G}_{m,S}\simeq \mathbb{G}_{m,S}^{n-1}$.
\item The forgetful morphism $o:\LS_{X,D,d}\rightarrow \Pic_{X,D,d}$ is an $H_S$-torsor for the \'etale topology.
\item The morphism $o:\LS_{X,D,d}\rightarrow \Pic_{X,D,d}$ is formally smooth.
\item $\LS_{X,D,d}$ has relative dimension $g+n-1$ over $S$.
\end{enumerate}
\end{lemma}
\begin{proof}
(1) For each $\gs:=(L,\gamma)\in \LS_{X,D,d}(T)$ and each tuple $\underline{\lambda}:=(\lambda_{1},\dots,\lambda_{n})\in(\mathbb{G}_{m,S})^{n}(T)$ we construct $\underline{\lambda}\cdot \gs\in \LS_{X,D,d}(T)$ as follows. The invertible sheaf underlying $\underline{\lambda}\cdot \gs$ is $L$ itself. Let
$\gamma_i:\sigma_{iT}^{*} L \simeq \O_{T}, \ i=1,\dots,n$,
be the trivializations determined by $\gamma$. The level structure $\underline{\lambda}\cdot \gs$ is the one determined by the trivializations 
$$\sigma_{iT}^{*}L\overset{\gamma_i}{\simeq} \O_{T}\overset{\lambda_{i}}{\simeq}\O_{T}, \ i=1,\dots,n.$$
This defines an action of $(\mathbb{G}_{m,S})^{n}$ on $\LS_{X,D,d}$. The diagonal subgroup $\mathbb{G}_{m,S}$ acts trivially, since the trivializations $\lambda\circ\gamma_{i}$ and $\gamma_{i}$ are equivalent through the automorphism of $L$ given by multiplication by $\lambda^{-1}$. The action descends to $\eta: H_S \times \LS_{X,D,d}\rightarrow \LS_{X,D,d}$.\\
(2) The forgetful morphism $o:\LS_{X,D,d}\rightarrow \Pic_{X,D,d}$ is an epimorphism for the \'etale topology \cite[Proposition~2.1.2]{raynaud}, so it remains to prove that
the induced morphism 
$\eta^{\#}:H_S\times_{S} \LS_{X,D,d}\rightarrow \LS_{X,D,d}\times_{\Pic_{X,D,d}}\LS_{X,D,d}$
is an isomorphism.
\begin{enumerate}
\item[a)] We first show that $\eta^{\#}$ is a monomorphism, that is, that the elements $t_{i}:=\lambda_{i}^{'-1}\lambda_{i}\in\Gamma(T,\O_{T})^{*}$ coincide for all $i$ whenever $\underline{\lambda},\underline{\lambda}'\in(\mathbb{G}_{m,S})^{n}(T)$ and $\gs:=(L,\gamma)\in \LS_{X,D,d}(T)$ satisfy $\underline{\lambda}\cdot \gs=\underline{\lambda}'\cdot \gs$. The equality $\underline{\lambda}\cdot \gs=\underline{\lambda}'\cdot \gs$ amounts to the existence of an automorphism $\varphi:L\simeq L$ making the diagram
$$
\xymatrix{
\sigma_{iT}^{*}L\ar[r]^{\gamma_i}\ar[d]_{\varphi_{i}:=\sigma_{i}^{*}\varphi} & \O_{T} \ar[r]^{\lambda_{i}} & \O_{T}\ar@{=}[d]\\
\sigma_{iT}^{*}L\ar[r]_{\gamma_i} & \O_{T} \ar[r]^{\lambda_{i}^{'}} & \O_{T}
}
$$
commutative for every $i$. It follows that $t_i=\gamma_i \varphi_i \gamma_i^{-1}$ for each $i=1,\dots,n$. 
Since the fibers of $X_T\to T$ are geometrically connected, $\varphi$ is multiplication by a unit $u\in\Gamma(T,\O_T)^{*}$. Hence $\varphi_i$ is multiplication by $u$ for every $i$, and therefore $t_i=\gamma_i\varphi_i\gamma_i^{-1}=u$ for all $i$.
%Writing $\gamma_{ij}:=\gamma_i \gamma_{j}^{-1}$, one obtains $t_{i}=\gamma_{ij} t_{j}\gamma_{ij}^{-1}\in\Gamma(T,\O_{T})^{*}$, and therefore $t_i=t_j$ for all $i,j$. 
\item[b)] We next show that $\eta^{\#}$ is an epimorphism for the \'etale topology. 
Let $T$ be an $S$-scheme and let $((L_{1},\gamma_{1}),(L_{2},\gamma_2))$ be a $T$-valued point of $\LS_{X,D,d}\times_{\Pic_{X,D,d}}\LS_{X,D,d}$
for which there exists an isomorphism $\Psi:L_1\simeq L_{2}$. Let us set $\Psi_{i}:=\sigma_{iT}^{*}\Psi:\sigma_{iT}^{*}L_1\simeq \sigma_{iT}^{*}L_{2}$ and $\lambda_{i}:=\gamma_{1i}\Psi_{i}^{-1}\gamma_{2i}^{-1}\in\Gamma(T,\O_{T})^{*}$ for $i=1,\dots,n$, and let $\underline{\lambda}$ denote the tuple $(\lambda_{1},\dots,\lambda_{n})$. Since the diagram 
$$
\xymatrix{
\sigma_{iT}^{*}L_{1}\ar[d]_{\Psi_{i}}\ar[r]^{\gamma_{1i}} & \O_{T}\\
\sigma_{iT}^{*}L_{2}\ar[r]_{\gamma_{2i}} & \O_{T}\ar[u]^{\lambda_{i}}
}
$$
commutes for each $i=1,\dots,n$, one has $\underline{\lambda}\cdot(L_{2},\gamma_2)=(L_{1},\gamma_1)$ in $\LS_{X,D,d}(T)$. We now treat the general case. Let us fix $((L_{1},\gamma_{1}),(L_{2},\gamma_2))\in (\LS_{X,D,d}\times_{\Pic_{X,D,d}}\LS_{X,D,d} )(T)$ for which there exists an isomorphism $\Psi:L_1\simeq L_{2}\otimes f_{T}^{*}\mathcal{M}$ with $\mathcal{M}$ an invertible sheaf on $T$, and let $Z\rightarrow T$ be an \'etale surjective morphism trivializing $\mathcal{M}$. Let us denote by $((L'_{1},\gamma'_{1}),(L'_{2},\gamma'_2))\in (\LS_{X,D,d}\times_{\Pic_{X,D,d}}\LS_{X,D,d} )(Z)$ and by $\Psi'$ the pullbacks to $Z$. The preceding case applies over $Z$ and produces $\underline{\lambda}'\in(\mathbb{G}_{m,Z})^{n}$ with $\underline{\lambda}'\cdot (L'_{1},\gamma'_{1})=(L'_{2},\gamma'_{2})$, which proves the assertion.
\item[c)] A morphism of algebraic spaces that is both a monomorphism and an epimorphism for the \'etale topology is an isomorphism. Hence $\eta^{\#}$ is an isomorphism of algebraic spaces (see \cite[\S~3.4.1, Definition~3.6, Corollaire~3.7.1]{laumon}).
\end{enumerate}
(3) The morphism $o:\LS_{X,D,d}\rightarrow \Pic_{X,D,d}$ is an $H_S$-torsor by (2), so it suffices to prove that $H_S$ is formally smooth over $S$. This holds because $\mathbb{G}_{m,S}$ is smooth over $S$.

(4) By (2) the morphism $o$ is an $H_S$-torsor, and $H_S\simeq\mathbb{G}_{m,S}^{n-1}$ has relative dimension $n-1$ over $S$. Since $\Pic_{X,D,d}$ has relative dimension $g$ over $S$, the assertion follows.
\end{proof}

As shown above, $\LS_{X,D,d}$ is an algebraic $S$-space. Let $b\colon\LS_{X,D,d}\rightarrow S$ be its structure morphism. The pair $(\Goppa_{X,D,d},b)$ defines the \emph{extended Goppa morphism} 
$\Phi_{X,D,d}:\LS_{X,D,d}\rightarrow \Gr(k,n)_{S}:=\Gr(k,n)\times S$.

\begin{lemma}\label{alvarez-porras}
Assume $d>2g-1$ and
\begin{equation}\label{eq:mono-numerical}
(n-d)(d-g)>d.
\end{equation}
Then the extended Goppa morphism
$\Phi_{X,D,d}: \LS_{X,D,d}\rightarrow \Gr(k,n)_{S}$
is a monomorphism. Condition \eqref{eq:mono-numerical} holds, in particular, whenever $n/2>d>2g-1$ and $d\geq1$.
\end{lemma}
\begin{remark}
Condition \eqref{eq:mono-numerical} is exactly what the slope estimate of Step~2
below requires. It holds whenever $n/2>d>2g-1$ and $d\geq1$, which is the range in which the
result is established, for $S$ the spectrum of a field, in \cite{alvarez2}.
Since $(n-d)(d-g)>d$ and $d-g>0$ force $n>d$, condition
\eqref{eq:mono-numerical} implies $n>d>2g-1$. Deciding whether the monomorphism
persists on the larger range $n>d>2g-1$ would require a bound on the minimal
degree of a quotient line bundle of the kernel bundle $K$ sharper than the one
provided by semistability. We do not address this question here.
Accordingly, every statement below that invokes this lemma carries
\eqref{eq:mono-numerical} explicitly among its hypotheses. When $n$ is determined by
$g$ and $d$ (as in the self-duality results, where $n=2(1-g+d)$), condition
\eqref{eq:mono-numerical} is equivalent to $d\geq 2g+1$, and this is the form used
there.
\end{remark}
\begin{proof}
By Remark~\ref{rmk:LS-sheaf} the functor $\LS_{X,D,d}$ is a sheaf of sets, so
$\Phi_{X,D,d}$ is a monomorphism if and only if the map
$\LS_{X,D,d}(T)\to\Gr(k,n)_{S}(T)$ is injective for every $S$-scheme $T$. We
verify this condition directly. Let us fix $T$, write $\pi\colon X_T\to T$ for the base-changed family, and let $(L,\gamma)$
and
$(L',\gamma')$ be two level structures over $(X_T,\sigma_{iT})$ with the same
image under $\Phi_{X,D,d}$, that is, inducing the same subbundle
$C\subseteq\O_T^n$. Such a coincidence is the datum of an isomorphism $\lambda$
making the triangle
\begin{equation}\label{eq:iso-com}
\xymatrix{
\pi_{*}L\ar[rr]^{\simeq}_{\lambda} \ar[rd]_{\pi_{*}(\gamma)} & & \pi_{*}L'\ar[dl]^{\pi_{*}(\gamma')}\\
& \O_{T}^{n}&
}
\end{equation}
commute. Since $d>2g-1>2g-2$, cohomology and base change give
$\pi_{*}L(-D)=\pi_{*}L'(-D)= R^{1}\pi_{*}L=R^{1}\pi_{*}L'=0$, and $\pi_{*}L$, $\pi_{*}L'$ are
locally free of rank $1-g+d$, their formation commuting with base change. Moreover,
$d\ge 2g$ implies that $L$ and $L'$ are globally generated on every geometric fiber,
so the evaluation maps $\psi\colon \pi^{*}\pi_{*}L\to L$, $\psi'\colon \pi^{*}\pi_{*}L'\to
L'$ are surjective with kernels $K:=\ker\psi$, $K':=\ker\psi'$ locally free of rank
$d-g$ and fiberwise degree $-d$. All this fits into the diagram
\begin{equation}\label{eq:diagram}
\xymatrix{
&0\ar[r]&L(-D)\ar[rd]&&&\\
0\ar[r] &  K\ar[r]^{i} & \pi^{*}\pi_{*}L\ar[r]^{\psi} \ar[dd]^{\pi^{*}(\lambda)}_{\simeq} & L\ar[r]\ar[rd]^{\gamma} & 0 &\\
&&&& i_{D*}(\O_{D})\ar[r] & 0\\
0\ar[r] & K'\ar[r]^{i'} & \pi^{*}\pi_{*}L'\ar[r]^{\psi'} &L'\ar[r] \ar[ru]_{\gamma'}& 0 & \\
&0\ar[r]&L'(-D)\ar[ru]&& &
}
\end{equation}
whose pentagon commutes by \eqref{eq:iso-com}, that is,
$\gamma'\circ\psi'\circ \pi^{*}\lambda=\gamma\circ\psi$.

\emph{Step 1: reduction to a fiberwise vanishing.} Let us set
$\theta:=\psi'\circ \pi^{*}(\lambda)\circ i\colon K\to L'$. Since $\psi\circ i=0$, the commutativity of the pentagon gives $\gamma'\circ\theta=\gamma\circ\psi\circ i=0$, so $\theta$ factors through
$\ker\gamma'=L'(-D)$. We regard $\theta$ as a morphism $K\to L'(-D)$, that is, as a global section of the locally free sheaf $K^{\vee}\otimes L'(-D)$ on $X_T$.
We claim that $\theta=0$. It suffices to show
$\pi_{*}\bigl(K^{\vee}\otimes L'(-D)\bigr)=0$, because then
$\theta\in H^{0}(X_T,K^{\vee}\otimes L'(-D))=H^{0}(T,\pi_{*}(K^{\vee}\otimes
L'(-D)))=0$. Write $\mathcal{F}:=K^{\vee}\otimes L'(-D)$. It is locally free of finite rank on $X_T$, and so of finite presentation and flat over $T$. Assume that $H^{0}(X_{\bar t},\mathcal{F}_{\bar t})=0$ for every geometric point $\bar t$ of
$T$. By flat base change the same vanishing holds on every fiber $X_{t}$ with
$t\in T$. The assertion $\pi_{*}\mathcal{F}=0$ is local on $T$, so we may assume $T=\Spec A$. Since $\pi$ is proper and of finite presentation, the object
$E:=R\pi_{*}\mathcal{F}$ of $D(A)$ is perfect and its formation commutes with
arbitrary base change \cite[Tag~0B91]{stacks}. In particular, $H^{i}(E\otimes^{\mathbf{L}}_{A}\kappa(t))=H^{i}(X_{t},\mathcal{F}_{t})$ for every
$t\in T$. These groups vanish for $i<0$, for $i=0$ by hypothesis, and for $i>1$ because the fibers of $\pi$ are curves. It follows that $E$ is perfect of tor-amplitude in $[1,1]$ \cite[Tag~068V]{stacks}. In particular, $\pi_{*}\mathcal{F}=H^{0}(E)=0$. Thus it remains to prove the fiberwise statement over an algebraically closed field.

\emph{Step 2: the fiberwise vanishing.} Let us fix a geometric point $\bar t$ and omit it from the notation. Then $X$ is a smooth curve over an
algebraically closed field, the sheaf $K$ is the kernel bundle of the complete linear system $|L|$, of rank $d-g$ and slope
$\mu(K)=-d/(d-g)$, and it is semistable by \cite[Theorem~1.2]{butler} (see
Remark~\ref{rmk:butler-char}). Moreover, $L'(-D)$ is a line bundle of degree $d-n<0$. Suppose that $\Hom(K,L'(-D))\neq0$ and let $N\subseteq L'(-D)$ be the image of a nonzero morphism. Then $N$ has rank $1$ and $\deg N\le\deg L'(-D)=d-n$. On the other hand, $\deg N=\mu(N)\ge\mu(K)=-d/(d-g)$ because $N$
is a quotient of the semistable bundle $K$. Hence $-d/(d-g)\le d-n$, that is,
$(n-d)(d-g)\le d$, which contradicts \eqref{eq:mono-numerical}. Therefore $\Hom(K,L'(-D))=0$, as required.

\emph{Step 3: conclusion.}
By Steps 1 and 2 we have $\theta=0$, so $\pi^{*}(\lambda)(K)\subseteq\ker\psi'=K'$. The same argument applied to $\lambda^{-1}$ gives $\pi^{*}(\lambda)^{-1}(K')\subseteq K$. Hence $\pi^{*}(\lambda)$ restricts to an isomorphism $K\xrightarrow{\sim}K'$ and, passing to the quotients in
\eqref{eq:diagram}, descends to a unique isomorphism
$\bar\lambda\colon L\xrightarrow{\sim}L'$ with $\psi'\circ \pi^{*}\lambda=\bar\lambda
\circ\psi$. Restricting the pentagon along the marked sections and using that
$\psi,\psi'$ are surjective, the identity $\gamma'\circ\psi'\circ
\pi^{*}\lambda=\gamma\circ\psi$ yields $\gamma_i=\gamma_i'\circ\sigma_{iT}^{*}
\bar\lambda$ for every $i$. Thus $\bar\lambda$ is an isomorphism of level structures $(L,\gamma_i)\xrightarrow{\sim}(L',\gamma_i')$. It is unique because $\psi$ is an epimorphism and $\bar\lambda$ is determined by
$\lambda$. Therefore $(L,\gamma)$ and $(L',\gamma')$ define the same element of $\LS_{X,D,d}(T)$, so that $\Phi_{X,D,d}(T)$ is injective. Since $T$ was arbitrary, $\Phi_{X,D,d}$ is a monomorphism.
\end{proof}

\begin{remark}\label{rmk:butler-char}
For $g\geq1$, Step~2 uses Butler's theorem
\cite[Theorem~1.2]{butler}. It gives the semistability of the kernel bundle
$K$ for a globally generated line bundle of degree $d\geq2g$, with no
restriction on the characteristic. For $g=0$, one has
$K\simeq\O_{\P^1}(-1)^{\oplus d}$. Thus Step~2 and
Lemma~\ref{alvarez-porras} are valid in arbitrary characteristic.
\end{remark}

\begin{theorem}\label{th:relativestructure}
The algebraic space $\LS_{X,D,d}$ is representable in the category of $S$-schemes, and the forgetful morphism $o\colon\LS_{X,D,d}\to\Pic_{X,D,d}$ is affine. Moreover, the structure morphism $\LS_{X,D,d}\rightarrow S$ is separated, of finite presentation, and smooth of relative dimension $g+n-1$.
\end{theorem}
\begin{proof}

By Lemma~\ref{lm:preliminariesraynaud}(2) the morphism $o$ is an $H_S$-torsor for
the \'etale topology, and $H_S\simeq\mathbb{G}_{m,S}^{n-1}$ is affine, smooth and of finite
presentation over $S$. Choose an \'etale covering $\{U_i\to\Pic_{X,D,d}\}_{i\in I}$
trivializing it, so that $\LS_{X,D,d}\times_{\Pic_{X,D,d}}U_i\simeq H_S\times_S U_i$
is affine over $U_i$ for every $i$. An \'etale covering is in particular an fpqc
covering, and descent data for affine morphisms are effective
\cite[Tag~0245]{stacks}. Hence $\LS_{X,D,d}$ is representable by a scheme and $o$
is affine.

Being affine, $o$ is separated and quasi-compact. Being an $H_S$-torsor, it is
smooth of relative dimension $n-1$ and of finite presentation, these properties
being local for the \'etale topology on the target. Composing with the proper
morphism $\Pic_{X,D,d}\to S$ we conclude that $\LS_{X,D,d}\to S$ is separated,
of finite presentation and smooth of relative dimension $(n-1)+g$, in accordance
with Lemma~\ref{lm:preliminariesraynaud}(4).
\end{proof}

\subsection{Theta sections on \texorpdfstring{$\LS_{X,D,d}$}{LS(X,D,d)} and information sets}\label{subsec:theta-relative}

Assume $n>d>2g-2$ and keep $k=1-g+d$. Write $Y:=\LS_{X,D,d}$ and let $\pi_Y\colon X_Y:=X\times_S Y\to Y$ be the base-changed family, carrying the universal level structure $(\mathcal L,\gamma_1,\dots,\gamma_n)$. This exists because $\LS_{X,D,d}$ represents its functor of points (Remark~\ref{rmk:LS-sheaf}). Since $d>2g-2$ we have $R^1\pi_{Y*}\mathcal L=0$ and $\pi_{Y*}\mathcal L$ is locally free of rank $k$, and $\Phi_{X,D,d}$ classifies the subbundle $\pi_{Y*}\mathcal L\subset\O_Y^{n}$.

Let $\mathcal U\subset\O^{n}_{\Gr(k,n)_{S}}$ be the tautological subbundle, put $\O(1):=(\det\mathcal U)^{-1}$, and for $I\subset\{1,\dots,n\}$ with $|I|=k$ let $p_I\in H^{0}(\Gr(k,n)_{S},\O(1))$ be the Pl\"ucker section, that is, the determinant of the coordinate projection $\mathcal U\to\O^{I}$. Set
\begin{equation*}\label{eq:theta_sheaf}
\mathcal Q:=\Phi_{X,D,d}^{*}\O(1)
\end{equation*}
and $\theta_I:=\Phi_{X,D,d}^{*}p_I\in H^0(\LS_{X,D,d},\mathcal Q)$. Since
$\Phi_{X,D,d}^{*}\mathcal U=\pi_{Y*}\mathcal L$, we get $\mathcal Q\simeq(\det\pi_{Y*}\mathcal L)^{-1}$ and $\theta_I=\det(\ev_I)$,
where $\ev_I\colon\pi_{Y*}\mathcal L\to\O^{I}_{Y}$ is evaluation along
$D_I:=\sum_{i\in I}\sigma_i$ read through the trivializations $\gamma_i$.

On the Picard side, $\pi$ has the section $\sigma_1$ and $\pi_*\O_X=\O_S$, so by \cite[(1.3)]{raynaud} the relative Picard functor is computed by actual
invertible sheaves and there is a Poincar\'e sheaf $\mathcal P$ on $X\times_S\Pic_{X,D,g-1}$, unique up to isomorphism once normalized by $\sigma_1^{*}\mathcal P\simeq\O$. Let 
$\rho\colon X\times_S\Pic_{X,D,g-1}\to\Pic_{X,D,g-1}$ be the projection. The
relative degree of $\mathcal P$ is $g-1$, so its Euler characteristic on the
fibers of $\rho$ vanishes. Consequently, the determinant of cohomology $\mathcal T:=(\det R\rho_{*}\mathcal P)^{-1}$ is unchanged when $\mathcal P$ is
twisted by the pullback of an invertible sheaf on the base, and it carries a canonical section $\theta$ \cite{deligne-determinant}. We write $W_{X,D}:=Z(\theta)$, the relative theta locus.

Since $\deg\mathcal L(-D_I)=d-k=g-1$, the assignment
$\mathcal L\mapsto\mathcal L(-D_I)$ defines a morphism of $S$-schemes
$\tau_I\colon\LS_{X,D,d}\to\Pic_{X,D,g-1}$. Recall that, for a $k$-dimensional
code $C\subset \k^{n}$ over a field $\k$, a subset $I$ of cardinality $k$ is an
\emph{information set} if the coordinate projection $C\to\k^{I}$ is an
isomorphism.

\begin{theorem}\label{thm:theta-information-relative}
Assume $n>d>2g-2$ and let $I\subset\{1,\dots,n\}$ with $|I|=k$.
\begin{enumerate}
\item There is a canonical isomorphism $\tau_I^{*}\mathcal T\simeq\mathcal Q$ under which $\tau_I^{*}\theta=\theta_I$. In particular $Z(\theta_I)=\tau_I^{-1}(W_{X,D})$ as closed subschemes of $\LS_{X,D,d}$.
\item Let $\xi$ be a geometric point of $\LS_{X,D,d}$, with associated invertible sheaf $L_\xi$ and code $C_\xi$. The following are equivalent: $\theta_I(\xi)\neq0$, $I$ is an information set of $C_\xi$, 
$H^{0}(X_\xi,L_\xi(-D_{I\xi}))=0$, $\tau_I(\xi)\notin W_{X,D}$.
\end{enumerate}
\end{theorem}

\begin{proof}
(1) On $X_Y$ the exact sequence
$0\to\mathcal L(-D_{IY})\to\mathcal L\to\mathcal L|_{D_{IY}}\to0$ together with
the trivializations $\gamma_i$ identifies $\pi_{Y*}(\mathcal L|_{D_{IY}})$ with
$\O^{I}_{Y}$. As $R^{1}\pi_{Y*}\mathcal L=0$, the complex $R\pi_{Y*}\mathcal L(-D_{IY})$ is represented by $[\pi_{Y*}\mathcal L\xrightarrow{\ev_I}\O^{I}_{Y}]$ in degrees $0$ and $1$. Taking determinants,
$(\det R\pi_{Y*}\mathcal L(-D_{IY}))^{-1}\simeq(\det\pi_{Y*}\mathcal L)^{-1}
=\mathcal Q$, and the canonical section of the left-hand side is
$\det(\ev_I)=\theta_I$.

By construction of $\tau_I$ there is an invertible sheaf $N$ on $Y$ with $(\mathrm{id}_{X}\times\tau_I)^{*}\mathcal P\simeq \mathcal L(-D_{IY})\otimes\pi_Y^{*}N$. The formation of the determinant of cohomology commutes with base change \cite{deligne-determinant}, so $\tau_I^{*}\mathcal T\simeq \bigl(\det R\pi_{Y*}((\mathrm{id}_{X}\times\tau_I)^{*}\mathcal P)\bigr)^{-1}$ compatibly with the canonical sections. Since the Euler characteristic vanishes, twisting by $\pi_Y^{*}N$ leaves this determinant unchanged. Combining the two computations gives (1).

(2) At $\xi$ the section $\theta_I$ is $\det(\ev_{I\xi})$, which is nonzero exactly when $\ev_{I\xi}\colon H^{0}(X_\xi,L_\xi)\to\kappa(\xi)^{I}$ is an isomorphism. Under $\Phi_{X,D,d}$ this map is the coordinate projection $C_\xi\to\kappa(\xi)^{I}$, so this happens exactly when $I$ is an information set of $C_\xi$. Its kernel is $H^{0}(X_\xi,L_\xi(-D_{I\xi}))$, which gives the third condition, and the fourth is (1) evaluated at $\xi$.
\end{proof}

\begin{remark}\label{rmk:matroid-erasure}
At a geometric point, the subsets $I$ with $\theta_I\neq0$ are exactly the bases of the matroid represented by $C_\xi$ \cite{oxley-matroid-theory}. Equivalently, the complementary of $I$, $I^{c}$, is a
uniquely correctable erasure pattern of size $n-k$. The equivalence between being an information 
set and the vanishing of $H^{0}(X_\xi,L_\xi(-D_{I\xi}))$
is the standard evaluation criterion for algebraic-geometric codes. What
Theorem~\ref{thm:theta-information-relative} adds is that these loci are the pullbacks of a single theta locus along the morphisms $\tau_I$, hence the zero schemes of sections of one invertible sheaf $\mathcal Q$ on $\LS_{X,D,d}$.
\end{remark}

\begin{remark}\label{rmk-classical-theta}
Assume $\k=\mathbb C$ and let $X$ be a geometric fiber. After choosing a symplectic basis of $H_1(X,\mathbb Z)$ and an identification of $\Pic^{g-1}(X)$ with the Jacobian of $X$, the canonical theta section of $\mathcal T$ is represented analytically, up to the usual translation and a
nowhere-vanishing factor depending on these choices, by the classical Riemann
theta function. Consequently, Theorem~\ref{thm:theta-information-relative} identifies the restriction of $\theta_I$ to the corresponding fiber of
$\LS_{X,D,d}$ with the pullback of a classical Riemann theta function along $L\mapsto L(-D_I)$. Thus, over $\mathbb C$, the non-vanishing condition characterizing information sets may equivalently be expressed as the non-vanishing of a classical theta function at the point of $\Pic^{g-1}(X)$ represented by $L(-D_I)$. The determinant-of-cohomology formulation used above is intrinsic and remains valid over arbitrary base fields and in families.
\end{remark}

Write $\Gr_I(k,n)_{S}:=\{p_I\neq0\}$ and $U_I:=\Phi_{X,D,d}^{-1}(\Gr_I(k,n)_{S})=\{\theta_I\neq0\}$. Since every code
of dimension $k$ admits an information set, the open sets $U_I$ cover $\LS_{X,D,d}$. Writing $I=\{i_1<\cdots<i_k\}$, every point of $\Gr_I(k,n)_{S}$ has a unique systematic generator matrix whose $i_r$-th column
is the $r$-th standard basis vector, and its remaining entries $a_{rj}$, for $j\notin I$, identify $\Gr_I(k,n)_{S}\simeq\mathbb A^{k(n-k)}_{S}$. Each $a_{rj}$ equals the signed ratio $p_{(I\setminus\{i_r\})\cup\{j\}}/p_I$, the
sign being the one induced by writing the replacement subset in increasing
order. Hence on $U_I$ the systematic entries of the universal code are ratios of theta sections, $a_{rj}=\pm\,\theta_{(I\setminus\{i_r\})\cup\{j\}}/\theta_I$.

This raises the question of how much of the level structure these ratios retain. In the range of Lemma~\ref{alvarez-porras}, the monomorphism property
shows that they determine the level structure functorially over $S$. The purpose of the rest of this subsection is to show that, under one further numerical hypothesis, they do more scheme-theoretically, in that on each $U_I$ the systematic generator matrix realizes the space of level structures as a locally closed subscheme 
of the affine space of systematic matrices.

\subsection{Canonical immersion of \texorpdfstring{$\LS_{X,D,d}$}{LS(X,D,d)}}

The immersion statement requires one further numerical hypothesis. The monomorphism of Lemma~\ref{alvarez-porras} does not suffice by itself, since a monomorphism locally of finite type need not be an immersion. We shall use the following criterion.

\begin{lemma}\label{lem:immersion-criterion}
Let $Z$ be a locally noetherian scheme and let $j\colon Y\to Z$ be a monomorphism of finite type. Assume that for every discrete valuation ring $A$ with fraction field $K$, every morphism $\Spec A\to Z$ whose set-theoretic image is contained in that of $j$, and every lift $\Spec K\to Y$ of the generic point, the lift extends to $\Spec A\to Y$. Then $j$ is an immersion.
\end{lemma}

\begin{proof}
The extension is unique because $j$ is a monomorphism, so the hypothesis is condition d) of \cite[Corollaire~15.7.6]{ega4-3}. By condition a) of \emph{loc.\ cit.}, $j$ factors as a proper morphism $g\colon Y\to Z''$ followed
by an immersion $Z''\to Z$. Since $j$ is a monomorphism, so is $g$, and a proper
monomorphism is a closed immersion \cite[Corollaire~18.12.6]{ega4-4}.
\end{proof}

The extension property required by Lemma~\ref{lem:immersion-criterion} for the level parameters is established below.

Let $A$ be a discrete valuation ring over $S$, with uniformizer $t_A$, fraction field $K$ and residue field $\kappa$, and let $v\colon K^{*}\to\mathbb Z$ be its normalized valuation, so that $v(t_A)=1$.  
We write
\[
H_A:=(\mathbb G_{m,A})^{n}/\mathbb G_{m,A}
\]
for the quotient of $(\mathbb G_{m,A})^{n}$ by the diagonal embedding, a split torus over $A$. It will act on the level trivializations introduced below.

\begin{lemma}\label{lem:character-lattice}
Set $\Lambda_0:=\{(a_1,\ldots,a_n)\in\mathbb Z^n\mid a_1+\cdots+a_n=0\}$. Then the groups of characters and cocharacters satisfy
\[
X^{*}(H_A)\simeq\Lambda_0,
\qquad
X_{*}(H_A)\simeq\mathbb Z^{n}/\langle(1,\ldots,1)\rangle,
\]
the class of the $i$-th basis vector of $\mathbb Z^{n}$ corresponding to the $i$-th coordinate of $(\mathbb G_{m,A})^{n}$.
\end{lemma}

\begin{proof}
For a commutative group $\Lambda$ let $\Diag_A(\Lambda):=\Spec A[\Lambda]$ be
the associated diagonalizable $A$-group scheme, where $A[\Lambda]$ is the group
algebra of $\Lambda$ over $A$, so that $\Diag_A(\mathbb Z)=\mathbb G_{m,A}$ and $\Diag_A(\Lambda)(B)=\Hom_{\mathbb Z}(\Lambda,B^{*})$ for every $A$-algebra $B$.
The contravariant functor $\Diag_A$ is exact \cite[Ch.~II, \S1, no.~2.11]{demazure-gabriel}. Applied to the exact sequence
\begin{equation}\label{eq:sum-of-coordinates}
0\to\Lambda_0\to\mathbb Z^{n}\overset{\Sigma}{\longrightarrow}\mathbb Z\to0,
\end{equation}
with $\Sigma$ the sum of the coordinates, it gives
\[
1\to\mathbb G_{m,A}\to(\mathbb G_{m,A})^{n}\to\Diag_A(\Lambda_0)\to1,
\]
the first morphism being the diagonal embedding. Hence $H_A\simeq\Diag_A(\Lambda_0)$. Since $\Spec A$ is connected, $A$ being local, \emph{loc.\ cit.} gives $X^{*}(\Diag_A(\Lambda))\simeq\Lambda$, and the
antiequivalence gives
\[
\Hom_{A\text{-gr}}(\Diag_A(\Lambda'),\Diag_A(\Lambda))\simeq
\Hom_{\mathbb Z}(\Lambda,\Lambda').
\]
The first of these yields $X^{*}(H_A)\simeq\Lambda_0$.

For the cocharacters, a cocharacter is a morphism $\Diag_A(\mathbb Z)\to\Diag_A(\Lambda_0)$, so the second formula with $\Lambda'=\mathbb Z$ and $\Lambda=\Lambda_0$ gives $X_{*}(H_A)\simeq\Hom_{\mathbb Z}(\Lambda_0,\mathbb Z)$. Since $\mathbb Z$ is
free, \eqref{eq:sum-of-coordinates} splits, so applying the contravariant functor $\Hom_{\mathbb Z}(-,\mathbb Z)$ to it yields an exact sequence
\[
0\to\mathbb Z\xrightarrow{\;\Sigma^{\vee}\;}\mathbb Z^{n}\to
\Hom_{\mathbb Z}(\Lambda_0,\mathbb Z)\to0,
\]
where $\mathbb Z^{n}$ is identified with its dual by the standard basis. Here $\Sigma^{\vee}(1)=\Sigma=(1,\ldots,1)$, hence $X_{*}(H_A)\simeq\mathbb Z^{n}/\langle(1,\ldots,1)\rangle$. Under this
identification the class of $e_i$ is the cocharacter obtained by composing the
inclusion of the $i$-th factor $\mathbb G_{m,A}\hookrightarrow(\mathbb G_{m,A})^{n}$ with the projection onto $H_A$.
\end{proof}

Put $X_A:=X\times_S\Spec A$ and let $L$ be an invertible sheaf of relative degree $d$ on $X_A$, let $\sigma_{i,A}\colon\Spec A\to X_A$ be the $i$-th marked section, and put $L_i:=\sigma_{i,A}^{*}L$, an invertible $A$-module. The scheme of level
trivializations of $L$ along the marked sections is $\widetilde P_L:=\prod_{i=1}^{n}\underline{\Isom}_A(L_i,\O_A)$, a right torsor
under $(\mathbb G_{m,A})^{n}$. Since $A$ is local, each $L_i$ is free of rank
one, so $\widetilde P_L\simeq(\mathbb G_{m,A})^{n}$. The diagonal $\mathbb G_{m,A}$ acts on it by simultaneous rescaling $\lambda\cdot(\gamma_1,\ldots,\gamma_n)=(\lambda\gamma_1,\ldots,\lambda\gamma_n)$,
and we set
\begin{equation*}\label{eq:PL}
P_L:=\widetilde P_L/\mathbb G_{m,A},
\end{equation*}
a right torsor under $H_A$. Non-canonically $P_L\simeq(\mathbb G_{m,A})^{n-1}$. In particular $P_L$ is affine over $\Spec A$.

\begin{lemma}\label{lem:valuative-class}
Let $s_\eta\in P_L(K)$ be a rational section. Fix a section $s\in P_L(A)$ and let $h\in H_A(K)$ be the unique element with $s_\eta=s_K\cdot h$. Then the rule $\chi\mapsto v(\chi(h))$ for $\chi\in X^{*}(H_A)$ defines an element
\[
\nu(s_\eta)\in\Hom_{\mathbb Z}(X^{*}(H_A),\mathbb Z)\simeq X_{*}(H_A)
\]
that does not depend on $s$, the identification being that of
Lemma~\ref{lem:character-lattice}. Moreover,
\[
\nu(s_\eta)=0
\quad\Longleftrightarrow\quad
s_\eta\text{ extends to an }A\text{-point of }P_L.
\]
\end{lemma}

\begin{proof}
A section $s\in P_L(A)$ exists because $A$ is local, so that every invertible $A$-module is free and each $\sigma_{i,A}^{*}L$ admits a trivialization. Since $P_L(K)$ is a right $H_A(K)$-torsor, the element $h$ with $s_\eta=s_K\cdot h$ exists and is unique. We first check that $\nu(s_\eta)$ does
not depend on $s$. If $s'=s\cdot a$ is another section of $P_L(A)$, with $a\in H_A(A)$, then $s_\eta=s'_K\cdot h'$ with $h'=a_K^{-1}h$. For every
character $\chi$ one has $\chi(a_K)\in A^{*}$, and so $v(\chi(a_K))=0$. Thus $v(\chi(h'))=v(\chi(h))$ for all $\chi$, as required.

For the last assertion, choose a representative $\widetilde h=(h_1,\ldots,h_n)\in(K^{*})^{n}$ of $h$. By
Lemma~\ref{lem:character-lattice} a character of $H_A$ is an element $\mathbf a=(a_1,\ldots,a_n)\in\Lambda_0$, acting by $\mathbf a(h)=\prod_{i}h_i^{a_i}$, which is well defined because $\sum_i a_i=0$. Hence
\[
\nu(s_\eta)(\mathbf a)=v\Bigl(\prod_{i}h_i^{a_i}\Bigr)=\sum_{i}a_i\,v(h_i),
\]
so that $\nu(s_\eta)$ is the functional on $\Lambda_0$ given by the scalar
product with $(v(h_1),\ldots,v(h_n))$. Under the identification $X_{*}(H_A)\simeq\mathbb Z^{n}/\langle(1,\ldots,1)\rangle$ of
Lemma~\ref{lem:character-lattice}, it is therefore represented by $(v(h_1),\ldots,v(h_n))$, and $\nu(s_\eta)=0$ if and only if all the
valuations $v(h_i)$ are equal. Equivalently, after multiplying $\widetilde h$
by a common scalar in $K^{*}$, all its entries are units of $A$. This is
precisely the condition that $h\in H_A(A)$, in which case $s_K\cdot h$ extends
to an $A$-section of $P_L$. The converse is immediate.
\end{proof}

We now relate the valuative class to the base divisor of the limiting linear series.

Let $u:\Spec A\to\Gr(k,n)_S$ be an $A$-point and let $C_A\subset A^n$ be the associated rank-$k$ direct summand. Assume that the generic fiber $C_K$ is the
Goppa code of a level structure $\gs_\eta=(L_K,\gamma_{1,K},\ldots,\gamma_{n,K})$ on $X_K$, and let $L$ be an extension of $L_K$ to an invertible sheaf on $X_A$. The generic
trivializations $\gamma_{1,K},\ldots,\gamma_{n,K}$ define a rational section $s_\eta\in P_L(K)$, to which Lemma~\ref{lem:valuative-class} applies.

Since $C_K$ is identified with $H^0(X_K,L_K)$, the evaluation map gives a
morphism $\psi:\pi^*C_A\longrightarrow L\otimes_A K$, where $\pi:X_A\to\Spec A$ is the structure morphism. Choose the smallest integer $m$
such that $t_A^m\im(\psi)\subseteq L$ inside $L\otimes_A K$, and set $\phi:=t_A^m\psi:\pi^*C_A\longrightarrow L$. By minimality of $m$, the special
fiber $\phi_0:\pi_0^*C_0\longrightarrow L_0:=L|_{X_0}$ is nonzero. Otherwise $\im(\phi)\subseteq t_A L$, contradicting the minimality of $m$. 
Since $X_0$ is a smooth geometrically connected curve, its image is of the form $\im(\phi_0)=L_0(-B)$ with $B$ an effective divisor on $X_0$. Thus $B$ is the
fixed divisor of the limiting linear system defined by $\phi_0$. We introduce the notation
\begin{equation}\label{notation_limit}
\begin{split}
&L'':=L_0(-B), \qquad d'':=\deg L''=d-\deg B,\\
%\]
%\[
&V:=\im\bigl(C_0\to H^0(X_0,L'')\bigr),
\qquad
R:=\ker(C_0\to V),
\qquad
r:=\dim_\kappa R.
\end{split}
\end{equation}
The vector space $V$ generates $L''$, and defines a non-degenerate morphism $f_0:X_0\longrightarrow\P(V)$ with $f_0^{*}\O(1)=L''$. Let $\iota:\P(V)\hookrightarrow\P(C_0)$ be the closed immersion induced by the
quotient $C_0\twoheadrightarrow V$. For each $i$, let $\lambda_i:C_A\longrightarrow A$ be the projection onto the $i$-th coordinate,
set
\[
K_i^0:=\ker(\lambda_i\otimes_A\kappa:C_0\to\kappa),
\]
and let $q_i\in\P(C_0)$ be the point defined by the one-dimensional quotient
$\lambda_i\otimes_A\kappa:C_0\to\kappa$.

\begin{proposition}\label{prop:base-points-valuative}
Assume $n>d$ and that $C_0$ is non-degenerate. Then the following hold.
\begin{enumerate}
\item If $p_i\notin\Supp B$, then $K_i^0=\ker\bigl(C_0\to L''|_{p_i}\bigr)$. Equivalently, $\iota(f_0(p_i))=q_i$.
\item The valuative class $\nu(s_\eta)$ vanishes if and only if no marked point $p_i$ belongs to $\Supp B$.
\end{enumerate}
\end{proposition}

\begin{proof}
Choose generators $e_i$ of the free $A$-modules $\sigma_{i,A}^{*}L$, and let $e_i^\vee:\sigma_{i,A}^{*}L\to A$ be the dual trivializations. The tuple $(e_1^\vee,\ldots,e_n^\vee)$ defines a section $s\in P_L(A)$. Writing 
$\gamma_{i,K}=\beta_i e_i^\vee$ with $\beta_i\in K^{*}$, we have $s_\eta=s_K\cdot h$, where $h\in H_A(K)$ is the class of $(\beta_1,\ldots,\beta_n)$. Hence, by Lemma~\ref{lem:valuative-class}, $\nu(s_\eta)$ is represented by $(v(\beta_1),\ldots,v(\beta_n))$ in $\mathbb Z^n/\langle(1,\ldots,1)\rangle$. Put $a_i:=v(\beta_i)-m$. Since $(a_1,\ldots,a_n)$ and $(v(\beta_1),\ldots,v(\beta_n))$ differ by the constant
vector $m\cdot (1,\ldots,1)$, the class of $(a_1,\ldots,a_n)$ is again
$\nu(s_\eta)$.

Let $\mu_i:C_A\to A$ be the composition of $\sigma_{i,A}^*\phi$ with $e_i^\vee$.  Since $\phi=t_A^m\psi$ and $\gamma_{i,K}=\beta_i e_{i,K}^\vee$, the $i$-th coordinate map satisfies
\[
\lambda_i\otimes_A K
=
t_A^{-m}\beta_i\,(\mu_i\otimes_A K).
\]
Equivalently, after identifying $A$ with its image in $K$, $\lambda_i=t_A^{-m}\beta_i\mu_i$ in $\Hom_A(C_A,K)$.
Since $C_0$ is non-degenerate and $A$ is local, each coordinate map $\lambda_i:C_A\to A$ is surjective.
Write $\im(\mu_i)=t_A^{b_i}A$ with $b_i\ge0$. The equality $\lambda_i=t_A^{-m}\beta_i\mu_i$ and the surjectivity of $\lambda_i$ give $a_i+b_i=0$. It follows that $a_i\le0$, and
\[
a_i=0
\quad\Longleftrightarrow\quad
b_i=0
\quad\Longleftrightarrow\quad
\mu_i\otimes\kappa\ne0.
\]
Now $\mu_i\otimes\kappa$ is the evaluation at $p_i$ of $\phi_0:\pi_0^*C_0\to L_0$. Since $\im(\phi_0)=L_0(-B)$, this evaluation is nonzero if and only if $p_i\notin\Supp B$.

If $p_i\notin\Supp B$, then $a_i=0$, so $t_A^{-m}\beta_i$ is a unit of $A$ and
the maps $\lambda_i\otimes\kappa$ and $\mu_i\otimes\kappa$, viewed as maps from $C_0$ to a one-dimensional $\kappa$-vector space, differ by a nonzero scalar.
Their kernels therefore coincide. Moreover, $L''|_{p_i}\to L_0|_{p_i}$ is an isomorphism, since $p_i\notin\Supp B$. Hence
\[
K_i^0=\ker\bigl(C_0\to L''|_{p_i}\bigr).
\]
The quotient $C_0\to L''|_{p_i}$ factors through $V$, and the induced quotient $V\to L''|_{p_i}$ defines $f_0(p_i)$. Since it has the same kernel in $C_0$ as $\lambda_i\otimes\kappa$, the corresponding points of $\P(C_0)$ coincide. Thus $\iota(f_0(p_i))=q_i$, proving (1).

For (2), the sheaf $L''=L_0(-B)$ is generated by $V$, so $d''=\deg L''\ge0$ and $\deg B\le d$. 
Since the $p_i$ are distinct $\kappa$-rational points, $\#\{i\mid p_i\in\Supp B\}\leq\deg B$. 
As $\deg B\le d<n$, at least one marked point does not lie in $\Supp B$, and $a_i=0$ for every such index $i$. If no marked point lies in $\Supp B$, then every $a_i$ vanishes and $\nu(s_\eta)=0$. Conversely, if $p_j\in\Supp B$, then $a_j<0=a_i$ for an index $i$ with $p_i\notin\Supp B$, so the class of $(a_1,\ldots,a_n)$ modulo constants is nonzero.
\end{proof}

\begin{proposition}\label{lem:valuative-cancellation}
Assume $d>2g-1$ and $n\ge d+g+2$.
Let $A,K,\kappa,L,\gs_\eta$ and $C_A$ be as above. Assume that $C_0:=C_A\otimes_A\kappa$ is again a Goppa code of type $(g,n,d)$. That is, there
exists a level structure $\gs_0=(L'_0,\gamma'_{1,0},\ldots,\gamma'_{n,0})$ on $X_0$, with the same marked points, such that $C_{\gs_0}=C_0$. Then the generic level structure extends over $\Spec A$.
\end{proposition}

\begin{proof}
Since $d>2g-1$, the code $C_{0}=C_{\gs_{0}}$ is non-degenerate
(\S~\ref{subsec:parameters}). By Proposition~\ref{prop:base-points-valuative}(2), it suffices
to prove that no marked point lies in $\Supp B$. We may base change to an algebraic closure $\bar\kappa$ of $\kappa$ since degrees are preserved, $H^0(X_0,L'')\otimes_\kappa\bar\kappa
\simeq H^0(X_{0,\bar\kappa},L''_{\bar\kappa})$, and $p_i\in\Supp B$ if and only if its base change lies in $\Supp B_{\bar\kappa}$.

Suppose that at least one marked point lies in $\Supp B$, and define
\[
I:=\{i\in\{1,\ldots,n\}\mid p_i\notin\Supp B\},
\quad
D_I:=\sum_{i\in I}p_i.
\]
Since each marked point outside $I$ contributes at least one to $\deg B$, we have
\begin{equation}\label{eq:val-first}
n-|I|\le\deg B=d-d'',
\quad\text{so}\quad
n\le |I|+d-d''.
\end{equation}

\emph{Claim 1.}
We first claim that, under the above assumption, $r\ge1$ (\textrm{see \eqref{notation_limit}}).
If $r=0$, then $C_0\hookrightarrow H^0(X_0,L'')$, so $h^0(X_0,L'')\ge \dim C_0=k=d-g+1$.  Since $B\ne0$, we have $d''<d$, and
Riemann--Roch gives $h^1(X_0,L'')=h^0(X_0,L'')-d''+g-1\ge d-d''>0$. Consequently, $L''$ is special.  If $g=0$,
this is impossible because $L''$ is generated by $V$, so $d''\ge0$. On the other hand, a
line bundle on $\mathbb P^1$ with $h^1>0$ has degree at most $-2$.  If $g\ge1$, Clifford's theorem \cite[Theorem~IV.5.4]{hartshorne} gives
\[
h^0(X_0,L'')\le \frac{d''}{2}+1.
\]
Together with $h^0(X_0,L'')\ge d-g+1$, this implies $d''\ge2(d-g)$.  On the other hand, $h^1(X_0,L'')>0$ gives $d''\le2g-2$.  It follows that $d\le2g-1$, contradicting our hypothesis.  This proves the claim.

\emph{Claim 2.} We claim that, under the above assumption, $n\le d+g+1$.
Since $V$ generates $L''$, the morphism $f_0:X_0\to\mathbb P(V)$ is
non-degenerate, and a non-degenerate curve in $\mathbb P^{\dim V-1}$ has degree
at least $\dim V-1$ \cite[Lecture~1]{harris-ag}. Therefore $d''\ge\dim V-1$. As $\dim V=\dim C_0-r=k-r=d-g+1-r$, we get
\begin{equation}\label{eq:val-second}
d''\ge d-g-r.
\end{equation}
Since $C_0=C_{\gs_0}$, $d>2g-1$, and $n>d$, the evaluation map identifies $H^0(X_0,L'_0)$ with $C_0$. Indeed, $h^0(X_0,L'_0)=d-g+1=k$, and $H^0(X_0,L'_0(-D))=0$ because $\deg L'_0(-D)=d-n<0$.  Under this identification, $K_i^0=H^0(X_0,L'_0(-p_i))$ for every $i$, and consequently \begin{equation}
\bigcap_{i\in I}K_i^0=H^0(X_0,L'_0(-D_I)).
\end{equation}
For every $i\in I$, Proposition~\ref{prop:base-points-valuative}(1) gives $R\subseteq K_i^0$.  Therefore $R\subseteq H^0(X_0,L'_0(-D_I))$.
Since $r\ge1$, the line bundle $L'_0(-D_I)$ has a nonzero section.  It
follows that $d-|I|\ge0$.  Note that for any line bundle $M$ with a nonzero
section on a smooth projective geometrically connected curve, one has the elementary estimate $h^0(M)\le\deg M+1$.  Applying this to $M=L'_0(-D_I)$, we obtain
\begin{equation}\label{eq:val-third}
r\le h^0(X_0,L'_0(-D_I))\le d-|I|+1.
\end{equation}
Combining \eqref{eq:val-first}, \eqref{eq:val-second}, and \eqref{eq:val-third}, we obtain $n \le |I|+d-d'' \le |I|+d-(d-g-r) =|I|+g+r \le |I|+g+d-|I|+1 =d+g+1$.

This is a contradiction, since $n\ge d+g+2$.
Thus no marked point lies in $\Supp B$, and Proposition~\ref{prop:base-points-valuative}(2) gives $\nu(s_\eta)=0$.
\end{proof}

\begin{theorem}\label{th:immersion-relative}
Assume $d>2g-1$, $n\ge d+g+2$, and $(n-d)(d-g)>d$.
Then the extended relative Goppa morphism $\Phi_{X,D,d}\colon\LS_{X,D,d}\longrightarrow \Gr(k,n)_{S}$ is an immersion. Equivalently, $\LS_{X,D,d}$ is isomorphic to a locally closed subscheme of $\Gr(k,n)_{S}$.
\end{theorem}

\begin{proof}
Assume first that $S$ is locally noetherian. By
Theorem~\ref{th:relativestructure}, $\LS_{X,D,d}$ is a scheme of finite
presentation over $S$. Since $\Gr(k,n)_S\to S$ is separated, $\Phi_{X,D,d}$ is of finite type, and by Lemma~\ref{alvarez-porras} it is a
monomorphism. 
We now apply Lemma~\ref{lem:immersion-criterion}. Let $A$ be a discrete valuation ring with fraction field $K$ and residue field $\kappa$, and consider a commutative diagram
\[
\xymatrix{
\Spec K \ar[r]^-{\gs_\eta} \ar[d] &
\LS_{X,D,d} \ar[d]^-{\Phi_{X,D,d}} \\
\Spec A \ar[r]^-{u} &
\Gr(k,n)_S.
}
\]
Assume that the set-theoretic image of $\Spec A$ is contained in that of $\Phi_{X,D,d}$. Let $s\in\Spec A$ be the closed point and set $z_0:=u(s)$.
The fiber of $\Phi_{X,D,d}$ over $z_0$ is non-empty and its structure
morphism to $\Spec\kappa(z_0)$ is a monomorphism. Hence it is $\Spec\kappa(z_0)$ \cite[Tag~03DP]{stacks}. After base change $\kappa(z_0)\to\kappa$, the special fiber $C_0$ of $u$ is therefore a Goppa
code over $\kappa$.

Write $\gs_\eta=(L_K,\gamma_{1,K},\ldots,\gamma_{n,K})$. By the properness of $\Pic_{X,D,d}\to S$ (Theorem~\ref{th:relativestructure}), the class of $L_K$
extends over $\Spec A$. Since $X_A\to\Spec A$ admits the section $\sigma_{1,A}$ and $\PicG(A)=0$, this class is represented by an invertible
sheaf $L$ on $X_A$ extending $L_K$.
Consider, on the other hand, $C_A\subset A^n$, the rank-$k$ direct summand defined by $u$. Its generic
fiber is the Goppa code of $\gs_\eta$, while its special fiber $C_0$ is, as
shown above, again a Goppa code of type $(g,n,d)$. Hence
Proposition~\ref{lem:valuative-cancellation} yields an extension $\gs_A\in\LS_{X,D,d}(A)$ of $\gs_\eta$.

The two $A$-points $\Phi_{X,D,d}(\gs_A)$ and $u$ agree over $\Spec K$. Since $\Spec A$ is integral and $\Gr(k,n)_S$ is separated over $S$, they agree over $\Spec A$. Thus $\gs_A$ completes the valuative diagram, and it is unique because $\Phi_{X,D,d}$ is a monomorphism. Lemma~\ref{lem:immersion-criterion}
therefore shows that $\Phi_{X,D,d}$ is an immersion.

For arbitrary $S$, the result follows by noetherian approximation. The assertion is local on $S$, so we may assume $S=\Spec R$ and write $R$ as a
filtered colimit of noetherian rings $R_\alpha$. By the standard limit
arguments of \cite[\S\,8.8 and \S\,8.10]{ega4-3}, the family, the projective
structure, and the marked sections descend, for $\alpha$ sufficiently large,
to $S_\alpha=\Spec R_\alpha$. Since $\LS_{X,D,d}$, $\Gr(k,n)_S$, and $\Phi_{X,D,d}$ commute with base change, $\Phi_{X,D,d}$ is the base change of the corresponding immersion over $S_\alpha$, and hence is an immersion.
\end{proof}

On the open subset $U_I=\Phi_{X,D,d}^{-1}(\Gr_I(k,n)_S)$ associated with the information set $I$, the identification $\Gr_I(k,n)_S\simeq\mathbb A_S^{k(n-k)}$  shows that the systematic
generator matrix realizes $U_I$ as a locally closed subscheme of $\mathbb A_S^{k(n-k)}$.
The corresponding global projective statement is the following.

\begin{proposition}\label{prop:Q-very-ample}
Under the hypotheses of Theorem~\ref{th:immersion-relative}, put $N:=\binom{n}{k}-1$. The invertible sheaf $\mathcal Q$ is very ample on $\LS_{X,D,d}$ relative to $S$, and the corresponding immersion is
\[
\Theta_{X,D,d}:=
\mathrm{Pl}\circ\Phi_{X,D,d}\colon
\LS_{X,D,d}\longrightarrow\P_S^N,
\]
where $\mathrm{Pl}$ denotes the Pl\"ucker embedding.
It is defined by the theta sections $(\theta_I)_{|I|=k}$. More concretely, for a geometric point $\xi=(L,\gamma_1,\ldots,\gamma_n)$ of $\LS_{X,D,d}$, let
\[
\alpha_{I,\xi}:\mathcal T_{L(-D_I)}
\xrightarrow{\sim}\mathcal Q_\xi
\]
be the isomorphism induced by $\tau_I^*\mathcal T\simeq\mathcal Q$. Then
\[
\Theta_{X,D,d}(\xi)=
\left[
\alpha_{I,\xi}\bigl(\theta(L(-D_I))\bigr)
\right]_{|I|=k}
=
[\theta_I(\xi)]_{|I|=k}.
\]
Moreover, for every $a\geq1$ and every invertible sheaf $\mathcal A$ on $S$, $\mathcal Q^{\otimes a}\otimes b^*\mathcal A$ is very ample relative to $S$.
\end{proposition}

\begin{proof}
By Theorem~\ref{th:immersion-relative}, the composition of $\Phi_{X,D,d}$ with the Pl\"ucker embedding is an immersion. Since $\mathcal Q=\Phi_{X,D,d}^*\O(1)$, this proves that $\mathcal Q$ is very ample
relative to $S$.
On the other hand, by Theorem~\ref{thm:theta-information-relative}, the pullback of the Pl\"ucker
coordinate $p_I$ is $\theta_I$, and under $\tau_I^*\mathcal T\simeq\mathcal Q$ the section $\tau_I^*\theta$ corresponds to $\theta_I$. Hence
\(
\Theta_{X,D,d}=[\theta_I]_{|I|=k},
\)
and the stated pointwise formula follows.

Finally, relative very ampleness is preserved by positive tensor powers and by tensoring with the pullback of an invertible sheaf on $S$.
\end{proof}

\subsection{The category \texorpdfstring{$\LS_{g,n,d}$}{LS(g,n,d)}. Structure and geometry}

Lemma~\ref{lm:equiv-cat} and Theorem~\ref{th:relativestructure} show that the morphism $\Theta: \LS_{g,n,d}\rightarrow \M_{g,n}$ is schematic. Consequently, $\LS_{g,n,d}$ and $\Theta$ inherit the properties established for $\LS_{X,D,d}$ and for the structure morphism $\LS_{X,D,d}\rightarrow S$.

We first introduce the universal $n$-pointed Picard stack of degree $d$, denoted $\mathcal{P}ic_{g,n,d}$. Its objects are tuples $(S,X,\sigma_1,\dots,\sigma_n,L)$, where $(X\to S,\sigma_1,\dots,\sigma_n)\in\M_{g,n}(S)$ and $L$ is an invertible
sheaf of relative degree $d$ on $X$. A morphism from $(S_1,X_1,\sigma^1_1,\dots,\sigma^1_n,L_1)$ to $(S_2,X_2,\sigma^2_1,\dots,\sigma^2_n,L_2)$ is a pair $((f,u),\phi)$, where $(f,u)$ is a morphism of pointed curves and $\phi:L_1\simeq f^*L_2$ is an isomorphism.

The scalar automorphisms of the line bundle define a central subgroup $\mathbb G_m$ of the inertia of $\mathcal{P}ic_{g,n,d}$. We denote by
\[
\Pic_{g,n,d}:=\mathcal{P}ic_{g,n,d}\fatslash\mathbb G_m
\]
the rigidification by $\mathbb G_m$ \cite[Definition~7.4.53]{alper}. Let $\Psi:\Pic_{g,n,d}\to\M_{g,n}$ be the forgetful morphism. If $S\to\M_{g,n}$ corresponds to an $n$-pointed family $(X\to S,\sigma_1,\dots,\sigma_n)$, then there is a canonical identification
\[
\Pic_{g,n,d}\times_{\M_{g,n}}S\simeq\Pic_{X,D,d},
\]
where $D$ is the divisor defined by the sections $\sigma_1,\dots,\sigma_n$.

We recall the fundamental structure results concerning the stacks $\M_{g,n}$ and $\Pic_{g,n,d}$.

\begin{lemma}\label{lm:universalpreliminariesraynaud}
Assume $2g-2+n>0$. Then the following hold.
\begin{enumerate}
\item $\M_{g,n}$ is a smooth Deligne--Mumford stack of dimension $3g-3+n$.
\item $\Pic_{g,n,d}$ is an algebraic stack and ${\Psi}: \Pic_{g,n,d}\rightarrow \M_{g,n}$ is schematic and smooth of relative dimension $g$.
\end{enumerate}
\end{lemma}
\begin{proof}
(1) See \cite[Corollary~6.4.15 and Theorem~6.4.16]{alper}.
For the original constructions, see \cite{deligne-mumford} and, for the $n$-pointed case, \cite{knudsen}.
(2) This follows by the defining property of the rigidification and by Lemma~\ref{lm:preliminariesraynaud}.
\end{proof}

These results yield the following theorem.

\begin{theorem}\label{thm:structure-LS}
Assume $2g-2+n>0$. Then the following hold.
\begin{enumerate}
\item The morphism $\Theta: \LS_{g,n,d}\rightarrow \M_{g,n}$ is schematic.
\item $\LS_{g,n,d}$ is a Deligne--Mumford stack.
\item The morphism $\LS_{g,n,d}\rightarrow\Pic_{g,n,d}$ is an $\Ht$-torsor for the \'etale topology.
\item The morphism $\Theta: \LS_{g,n,d}\rightarrow \M_{g,n}$ is separated and smooth of relative dimension $g-1+n$.
\item $\dim (\LS_{g,n,d})=4g-4+2n$.
\end{enumerate}
\end{theorem}
\begin{proof}
We prove the five statements in turn.

(1) Let $(\pi\colon X\to S,\sigma_1,\dots,\sigma_n)\in\M_{g,n}(S)$ be a family of $n$-pointed curves, corresponding to a morphism $S\to\M_{g,n}$. By
Lemma~\ref{lm:equiv-cat} the fiber product $\LS_{g,n,d}\times_{\M_{g,n}}S$ is
equivalent to $\LS_{X,D,d}$, which by Theorem~\ref{th:relativestructure} is
representable by an $S$-scheme. Consequently, $\Theta$ is schematic.

(2) By Lemma~\ref{lm:universalpreliminariesraynaud}(1), and using $2g-2+n>0$, the stack $\M_{g,n}$ is Deligne--Mumford. The morphism $\Theta$ is schematic by (1), and a stack representable by schemes over a Deligne--Mumford stack is
again Deligne--Mumford \cite[Tag~03YO]{stacks}. Hence $\LS_{g,n,d}$ is a Deligne--Mumford stack.

(3) The forgetful morphism $\LS_{g,n,d}\to\Pic_{g,n,d}$ becomes, after base
change along any $S\to\M_{g,n}$ as in (1), the $H_S$-torsor $o\colon\LS_{X,D,d}\to\Pic_{X,D,d}$ of Lemma~\ref{lm:preliminariesraynaud}(2). Being an $\Ht$-torsor is local on the base for the \'etale topology, so the claim
follows.

(4) Separatedness follows from Theorem~\ref{th:relativestructure} after base change as in (1). For smoothness, $\Psi\colon\Pic_{g,n,d}\to\M_{g,n}$ is smooth
by Lemma~\ref{lm:universalpreliminariesraynaud}(2), and $\LS_{g,n,d}\to\Pic_{g,n,d}$ is an $\Ht$-torsor by (3), hence smooth. Therefore,
the composite $\Theta$ is smooth. Its relative dimension is $\dim \Ht+\dim(\Pic_{g,n,d}/\M_{g,n})=(n-1)+g=g-1+n$.

(5) By (4) and Lemma~\ref{lm:universalpreliminariesraynaud}(1), $\dim\LS_{g,n,d}=\dim\M_{g,n}+(g-1+n)=(3g-3+n)+(g-1+n)=4g-4+2n$.
\end{proof}

\begin{remark}
Theorem~\ref{thm:structure-LS} rests on Theorem~\ref{th:relativestructure} and on the smooth $\Ht$-torsor structure of $\LS_{g,n,d}\to\Pic_{g,n,d}$
over the smooth morphism $\Pic_{g,n,d}\to\M_{g,n}$, and does not use
Lemma~\ref{alvarez-porras}. It requires only $2g-2+n>0$. In particular the dimension formula $\dim\LS_{g,n,d}=4g-4+2n$, used in the density
analysis of \S~\ref{subsec:density}, is available in the whole range $n>d>2g-1$. The numerical conditions $d>2g-1$ and $(n-d)(d-g)>d$ enter only in
the monomorphism and immersion statements for $\Phi_{g,n,d}$.
\end{remark}

Viewing $\Gr(k,n)$ as a category fibered over $\Sch$, we form the fiber product
$\Gr(k,n) \times_{\Sch} \M_{g,n}$. 
Its objects are the pairs
$((S,E\subset \O_{S}^{n}),(S,X\rightarrow S,\sigma_1,\dots,\sigma_n))$, 
and a morphism from $((S,E\subset \O_{S}^{n}),(S,X\rightarrow S,\sigma_1,\dots,\sigma_n))$ to $((T,F\subset \O_{T}^{n}),(T,Y\rightarrow T,\beta_1,\dots,\beta_n))$ is a pair $(u,f)$ in which $u:S\rightarrow T$ is a morphism of schemes with $u^{*}F=E$ and $(u,f)$ is a morphism of pointed curves (see Definition~\ref{def:categorylevel}).

As in the relative case (see \S~\ref{sec:relativestructure}), we construct the extended Goppa morphism
$$
\Phi_{g,n,d}: \LS_{g,n,d} \rightarrow \Gr(k,n) \times_{\Sch} \M_{g,n}.
$$
To a level structure $\gs=(S, X\rightarrow S, \sigma_1,\dots,\sigma_n, L, \gamma_1,\dots,\gamma_n)$ we associate the pair $((S,C_{\gs}),(S, X\rightarrow S, \sigma_1,\dots,\sigma_n))$.
\begin{theorem}\label{th:immersionstacks}
If
$d>2g-1, \ n\ge d+g+2, \ (n-d)(d-g)>d$,
the extended Goppa morphism
$\Phi_{g,n,d}: \LS_{g,n,d} \rightarrow \Gr(k,n) \times_{\Sch} \M_{g,n}$
is an immersion of stacks.
\end{theorem}
\begin{proof}
We must show that for every scheme $S$ and every morphism $h:S\rightarrow \Gr(k,n) \times_{\Sch} \M_{g,n}$ the induced morphism $\LS_{g,n,d}\times_{\Gr(k,n) \times_{\Sch} \M_{g,n} } S\rightarrow S$ is an immersion. The morphism $h$ is a pair $(h_1,h_2)$, where $h_1$ is determined by a locally free subsheaf $(E\subset \O_{S}^{n})\in \Gr(k,n)(S)$ and $h_2$ by a family of smooth projective curves $(X,\sigma_1,\dots,\sigma_n)\in \M_{g,n}(S)$. The fiber product $\LS_{g,n,d}\times_{\Gr(k,n) \times_{\Sch} \M_{g,n} } S$ is representable by the scheme $\Goppa_{X,D,d}^{-1}(E\subset \O_{S}^{n})$, so the assertion follows from Theorem~\ref{th:immersion-relative}.
\end{proof}

\subsection{Dimension bounds and (non-)density of the Goppa locus}\label{subsec:density}

Under the assumption $n>d>2g-1$, every Goppa code is non-degenerate (\S~\ref{subsec:parameters}), and the locus $\Gr(k,n)_{0}$ of non-degenerate linear codes is a Zariski open and
dense subset of the Grassmannian. It is therefore natural to ask whether the image of the Goppa morphism is dense in $\Gr(k,n)$. Related constructions occur in the theory of convolutional Goppa codes
\cite{convolutional-goppa,ag-convolutional-survey,mds-convolutional-goppa}.
A comparison of dimensions answers the question in
the negative for a range of degrees. Indeed, by
Theorem~\ref{thm:structure-LS}(5) one has $\dim\LS_{g,n,d}=4g-4+2n$. Thus, with
$a:=n-2+2g$ and $k=1-g+d$, we have
\begin{equation*}
\Xi(d):=\dim\Gr(k,n)-\dim\LS_{g,n,d}=-d^{2}+ad-\bigl(a(g+1)-(g-1)^{2}\bigr).
\end{equation*}
As a function of $d$, the polynomial $\Xi$ can be written as
\[
\Xi(d)=\frac{\Delta}{4}
-\left(d-\frac{a}{2}\right)^2,
\quad
\Delta:=(n-4)^2-16g.
\]
Thus $\Xi$ is a concave parabola with symmetry axis $d=a/2$. If $\Delta>0$, its roots are
\[
d_0^\pm=\frac{a\pm\sqrt{\Delta}}{2},
\]
and $\Xi(d)>0$ precisely for $d\in R(g,n):=(d_0^-,d_0^+)$. Consequently, for every admissible integer degree
\(
d\in R(g,n)\cap(2g-1,n),
\)
one has
$\dim\LS_{g,n,d}<\dim\Gr(k,n)$, and hence the Goppa locus is not dense in
$\Gr(k,n)$.

One may also ask about the injectivity properties of the Goppa morphism. In the range $n/2>d>2g+1$ the morphism is universally injective (Theorem~\ref{thm:universal-injectivity}). Since the proof rests on the
description of the fibers, we defer both the statement and its proof to Section~\ref{sec:fibers}.

\section{Fibers of the Goppa morphism}\label{sec:fibers}

Throughout this section we fix integers $g\geq 0$, $n\geq 1$, $d\geq 1$
satisfying $n>d>2g-2$, and set $k:=1-g+d$.  

Let $C\in\Gr(k,n)(\k)$ be a non-degenerate code and let $\mathcal{G}_C:=\LS_{g,n,d}\times_{\Gr(k,n)}\Spec\k$
be the fiber category of the Goppa morphism over $C$.  Writing
$\underline{C}_T\subset\O_T^n$ for $C\otimes_{\k}\O_T$, an object of $\mathcal{G}_C(T)$ is a level structure
\[
  (\pi:X\to T,\,\sigma_1,\dots,\sigma_n,\,L,\,\gamma_1,\dots,\gamma_n)
  \;\in\;\LS_{g,n,d}(T)
\]
whose associated Goppa code (see \eqref{eq:goppamap}) is $\underline{C}_T$. Our aim is to describe $\mathcal{G}_C$ in purely geometric terms.

\subsection{The geometric description}\label{subsec:geometric-description}

For every $i=1,\dots,n$, we denote by $\pr_i:\k^n\rightarrow \k$  the projection onto the $i$-th coordinate. Recall that $\P_C$ parametrizes the one-dimensional quotients of $C$.

We first construct the canonical projective data attached to the non-degenerate linear code $C$.

\begin{lemma}\label{lem:canonical-marked-projective}
Let $C\in\Gr(k,n)(\k)$ be a non-degenerate code. For each $i=1,\dots,n$, the coordinate functional $\ell_i:=\pr_i|_C:C\to\k$ determines a $\k$-rational point $q_i:=[\ell_i]\in\P_C$. Let 
$u_C:C\otimes\O_{\P_C} \twoheadrightarrow\O_{\P_C}(1)$ be the tautological quotient.  There is a unique isomorphism $\tau_i:\O_{\P_C}(1)|_{q_i}\xrightarrow{\sim}\k$ such that
$\tau_i\circ q_i^{\ast}u_C=\ell_i$.  In particular, $C$ canonically determines $\P_C$ with the labelled points $q_1,\dots,q_n$, 
together with the trivializations $\tau_i$ of $\O_{\P_C}(1)$ at those points.
\end{lemma}

\begin{proof}
Since $C$ is non-degenerate, $\ell_i\neq 0$ for every $i$, so $\ell_i$ defines the quotient point $q_i=[\ell_i]\in\P_C$. By the quotient convention for
projective space, the fiber at $q_i$ of the tautological quotient $u_C:C\otimes\O_{\P_C}\twoheadrightarrow\O_{\P_C}(1)$ is the one-dimensional
quotient represented by $\ell_i$. There is an isomorphism
$\tau_i:\O_{\P_C}(1)|_{q_i}\xrightarrow{\sim}\k$ with $\tau_i\circ q_i^{\ast}u_C=\ell_i$, and it is unique because $q_i^{\ast}u_C$ is surjective.
\end{proof}

We now define the geometric moduli problem attached to the fiber over $C$. Throughout, the points $q_i$ and the trivializations $\tau_i$ are those of Lemma~\ref{lem:canonical-marked-projective}.

\begin{definition}\label{def:P_C}
We define $\mathcal{P}_C$ to be the category fibered in groupoids over $\Sch$ whose objects over a $\k$-scheme $T$ are tuples $(\pi:X\to T,\sigma_1,\dots,\sigma_n,j)$ satisfying the following conditions.
\begin{enumerate}
\item $(\pi:X\to T,\sigma_1,\dots,\sigma_n)$ is a family of $n$-pointed smooth projective curves of genus $g$.
\item $j:X\hookrightarrow\P_C\times_\k T$ is a closed immersion over $T$ satisfying the following conditions.
\begin{enumerate}
\item[(2.a)] $j\circ\sigma_i=q_i\times\mathrm{id}_T$ for every $i=1,\dots,n$.
\item[(2.b)] the line bundle $L:=j^{\ast}\O_{\P_C\times T}(1)$ has relative degree $d$ over $T$, and for every geometric point $t\to T$ the curve $j_t(X_t)\subset\P_C\times_{\k}\kappa(t)$ is non-degenerate, that is, it is not contained in any hyperplane.
\end{enumerate}
\end{enumerate}
A morphism
$(\pi_1:X_1\to T_1,\sigma_i^1,j_1)\to
(\pi_2:X_2\to T_2,\sigma_i^2,j_2)$ is a morphism of pointed curves
$(u,f)$ (see Definition~\ref{def:categorylevel}), satisfying
$j_2\circ f=(\mathrm{id}_{\P_C}\times u)\circ j_1$.  The pullback along a
morphism $T'\to T$ is obtained by base-changing $X$, the sections $\sigma_i$,
and the closed immersion $j$.
\end{definition}

We can now state the main result of the subsection.

\begin{theorem}\label{thm:fibers-goppa-projective-moduli}
Let $n>d$, let $k:=1-g+d$, and let $C\subset\k^n$ be a non-degenerate $k$-dimensional linear code. Assume $d\geq 2g+1$. Then there is a
canonical equivalence of categories fibered in groupoids $\mathcal{G}_C\simeq\mathcal{P}_C$.  Thus the fiber of the Goppa morphism over $C$ is canonically equivalent to the moduli
stack of $n$-pointed smooth projective curves of genus $g$ over $\k$, equipped with a degree-$d$ non-degenerate closed immersion into $\P_C$ which sends the $i$-th marked point to $q_i$.
\end{theorem}

The proof constructs two functors. The first carries a level structure in the
fiber over $C$ to the morphism defined by its complete linear system, after identifying $\pi_{\ast}L$ with the free sheaf $\underline{C}_T$. The second
pulls back $\O_{\P_C}(1)$ along an embedded family and uses the distinguished
trivializations $\tau_i$ at the marked points.

\medskip
We first describe the two identifications of $\pi_{\ast}L$ with
$\underline{C}_T$.

\begin{lemma}\label{lem:fibers-goppa-restriction-isomorphism}
Let $(\pi:X\to T,\sigma_1,\dots,\sigma_n,j)$ be an object of $\mathcal{P}_C$ and set $L:=j^{\ast}\O_{\P_C\times T}(1)$.  Let
$p_T:\P_C\times T\to T$ and $p_{\P}:\P_C\times T\to\P_C$ be the projections,
and let $u_T:p_T^{\ast}\underline{C}_T\to\O_{\P_C\times T}(1)$ be the pullback of the tautological quotient $u_C$, composed with the canonical identification
$p_T^{\ast}\underline C_T\simeq p_{\P}^{\ast}(C\otimes\O_{\P_C})$.
Define
$$\rho:\underline{C}_T\to\pi_{\ast}L$$ 
as the adjoint of $\pi^{\ast}\underline{C}_T\simeq j^{\ast}p_T^{\ast}\underline{C}_T \xrightarrow{j^{\ast}u_T}L$.  Then $\rho$ is an isomorphism.
\end{lemma}

\begin{proof}
For every geometric point $t\to T$, the theorem on cohomology and base change, applied to the line bundle $L$ of relative degree $d>2g-2$, gives the canonical
identification $(\pi_{\ast}L)\otimes_{\O_T}\kappa(t)\simeq H^0(X_t,L_t)$ and shows that the formation of $\pi_{\ast}L$ is compatible with base change (see \cite[Tag~0D4E]{stacks}).  Under this identification, the fiber of $\rho$ at $t$ is the restriction map
\[
\rho_t:
C\otimes_\k\kappa(t)
\simeq
H^0(\P_C\times_\k\kappa(t),\O(1))
\xrightarrow{\,j_t^*\,}
H^0(X_t,L_t).
\]
We prove that $\rho_t$ is an isomorphism. By Definition~\ref{def:P_C}(2.b), $\deg L_t=d>2g-2$. It follows that $H^1(X_t,L_t)=0$, and Riemann--Roch gives $\dim_{\kappa(t)}H^0(X_t,L_t)=1-g+d=k=\dim_\k C$. The non-degeneracy of $j_t(X_t)$ means that no nonzero hyperplane section of $\O_{\P_C\times\kappa(t)}(1)$ vanishes identically on $j_t(X_t)$. Therefore $\rho_t$ is injective, and so an isomorphism.
Thus $\rho$ is an isomorphism on every geometric fiber. Since both $\underline C_T$ and $\pi_{\ast}L$ are locally free of rank $k$, it follows
that $\rho$ is an isomorphism \cite[Tag~00O0]{stacks}.
\end{proof}

\begin{lemma}\label{lem:auxiliary-identity-goppa}
Let $\gs=(\pi:X\to T,\sigma_1,\dots,\sigma_n,L,\gamma_1,\dots,\gamma_n)$
be an object of $\mathcal{G}_C(T)$, where $\gamma_i:\sigma_i^{\ast}L\to\O_T$ denotes the $i$-th component of the level
structure in adjoint form (see Lemma~\ref{lm:trivializations-surjection}).  Then $\pi_{\ast}\gamma:\pi_{\ast}L\to\O_T^n$ factors through an isomorphism
\begin{equation}
  \alpha_\gs:\pi_{\ast}L\xrightarrow{\sim}\underline{C}_T,
\end{equation}
characterized by $\pr_i\circ\alpha_\gs=(\pi_{\ast}\gamma)_i$ for all $i$, where $(\pi_{\ast}\gamma)_i:=\pr_i\circ\pi_{\ast}\gamma$.  Moreover,
\begin{equation}\label{eq:auxiliary-identity}
  \gamma_i\circ \ev_{\sigma_i} \circ\alpha_\gs^{-1}
  =\ell_i:\underline{C}_T\longrightarrow\O_T
\end{equation}
for every $i=1,\dots,n$,  where $\ev_{\sigma_i}:\pi_{\ast}L\to\sigma_i^{\ast}L$
is evaluation along the section $\sigma_i$.
\end{lemma}

\begin{proof}
Since $\gs$ lies in $\mathcal{G}_C(T)$, the image of $\pi_{\ast}\gamma$ is $\underline{C}_T\subset\O_T^n$. Consequently, $\pi_{\ast}\gamma$ factors through a surjection $\pi_{\ast}L\twoheadrightarrow\underline{C}_T$.  In the range $d>2g-2$,
cohomology and base change show that $\pi_{\ast}L$ is locally free of rank $k=1-g+d$. On the other hand, $\underline{C}_T$ has the same
rank. Therefore the
surjection is an isomorphism. This is $\alpha_\gs$.

Since $\pi\circ\sigma_i=\mathrm{id}_T$, the $i$-th component of $\pi_{\ast}\gamma$ is $(\pi_{\ast}\gamma)_i=\gamma_i\circ
\ev_{\sigma_i}$.  By the defining property of $\alpha_\gs$
and by the identity $\pr_i|_{\underline{C}_T}=\ell_i$, we get $\gamma_i\circ \ev_{\sigma_i} \circ\alpha_\gs^{-1} =(\pi_{\ast}\gamma)_i\circ\alpha_\gs^{-1}=\ell_i$.
\end{proof}

\paragraph{\textbf{Construction of the functor $F$.}}
Let $\gs=(\pi:X\to T,\sigma_1,\dots,\sigma_n,L,\gamma_1,\dots,\gamma_n)$
be an object of $\mathcal{G}_C(T)$.  Let $\alpha_\gs$ be the
isomorphism of Lemma~\ref{lem:auxiliary-identity-goppa}, and let $\varepsilon:\pi^{\ast}\pi_{\ast}L\to L$ be the adjunction map.  Define
\begin{equation}\label{eq:q-g-def}
  q_\gs:\pi^{\ast}\underline{C}_T
  \xrightarrow{\pi^{\ast}\alpha_\gs^{-1}}
  \pi^{\ast}\pi_{\ast}L
  \xrightarrow{\varepsilon}
  L.
\end{equation}
Since $d\geq 2g+1$, this morphism is surjective. 
By the universal property of the projective bundle, the quotient $q_\gs$ determines a unique $T$-morphism $j_\gs:X\to\mathbb P(\underline{C}_T)\simeq\P_C\times T$, together with
an isomorphism $\phi_\gs:j_\gs^{\ast}\O_{\P_C\times T}(1) \xrightarrow{\sim}L$ such that $q_\gs=\phi_\gs\circ q'_\gs$, where $q'_\gs:\pi^{\ast}\underline{C}_T\to j_\gs^{\ast}\O_{\P_C\times T}(1)$ is the tautological quotient
associated with $j_\gs$.  We set
$$F(\gs):=(\pi,\sigma_1,\dots,\sigma_n,j_\gs).$$

\begin{proposition}\label{prop:F-goppa-functor}
The assignment $\gs\mapsto F(\gs)$ defines a functor
$F:\mathcal{G}_C\to\mathcal{P}_C$.
\end{proposition}

\begin{proof}
Let us first check that $F(\gs)$ is an object of $\mathcal{P}_C(T)$.

The geometric fiber of $j_\gs$ at a geometric point $t\in T$ is the
morphism attached to the complete linear system $|L_t|$, and $L_t$ is
very ample because $d\geq2g+1$ \cite[Corollary~IV.3.2]{hartshorne}, so $j_{\gs,t}$ is a closed immersion. Since $X$ is proper over $T$ and $\P_C\times T$ is separated
over $T$, the morphism $j_\gs$ is proper \cite[Tag~01W6]{stacks}. Moreover, if $y\in\P_C\times T$, $\bar y$ is a geometric point above $y$, and $\bar t$ is its image in $T$, then the base change of $j_\gs^{-1}(y)$ to $\overline{\kappa(y)}$ is the fiber of $j_{\gs,\bar t}$ at $\bar y$, hence is either empty or $\Spec\overline{\kappa(y)}$. Thus $j_\gs^{-1}(y)$ is radicial and geometrically reduced over $\kappa(y)$, and $j_\gs$ is a closed immersion by \cite[Corollaire~18.12.6]{ega4-4}.

Note that, by construction, $j_\gs^{\ast}\O_{\P_C\times T}(1)\simeq L$, so every
geometric fiber has degree $\deg L_t=d$. Also, the image is non-degenerate. This is because under
$\alpha_{\gs,t}$ the restriction map
$C\otimes_\k\kappa(t)\to H^0(X_t,L_t)$ is an isomorphism, so no nonzero
hyperplane section restricts to zero on $X_t$.

It remains to check the marked points. Let us fix an index $i$. The $T$-point $j_\gs\circ\sigma_i$ of $\P_C\times T$ is classified by the invertible
quotient $\sigma_i^{\ast}q_\gs =\ev_{\sigma_i} \circ\alpha_\gs^{-1}$, by \eqref{eq:q-g-def} and the functoriality of the adjunction counit. Since $\gamma_i$ is an isomorphism, this quotient has the same kernel as $\gamma_i\circ\sigma_i^{\ast}q_\gs$,
which is $\ell_i$ by \eqref{eq:auxiliary-identity}. Hence $j_\gs\circ\sigma_i=q_i\times\mathrm{id}_T$.

Let us now verify functoriality. Let $(u,f,\varphi):\gs_1\to\gs_2$ be a morphism in $\mathcal{G}_C$, with $u:T_1\to T_2$, $f:X_1\to X_2$ and $\varphi:L_1\xrightarrow{\sim}f^{\ast}L_2$, and let 
$\beta_\varphi:(\pi_1)_{\ast}L_1\xrightarrow{\sim}u^{\ast}(\pi_2)_{\ast}L_2$ be the isomorphism induced by $\varphi$ and by cohomology and base change.
Identifying $u^{\ast}\underline C_{T_2}$ canonically with $\underline C_{T_1}$,
the compatibility of $\varphi$ with the trivializations $\gamma_i$ says exactly
that $u^{\ast}\alpha_{\gs_2}\circ\beta_\varphi=\alpha_{\gs_1}$
after composition with the inclusion $\underline C_{T_1}\hookrightarrow\O_{T_1}^n$. By naturality of the adjunction
map this yields $\varphi\circ q_{\gs_1}=f^{\ast}q_{\gs_2}$ as
invertible quotients of $\pi_1^{\ast}\underline C_{T_1}$, whence $j_{\gs_2}\circ f=(\mathrm{id}_{\P_C}\times u)\circ j_{\gs_1}$
by the universal property of $\mathbb P(\underline C_T)$ \cite[Tag~01OA]{stacks}. Compatibility with identities and composition is immediate.
\end{proof}

\paragraph{\textbf{Construction of the functor $G$.}}
Let $\mathfrak{p}=(\pi:X\to T,\sigma_1,\dots,\sigma_n,j)$ be an object of $\mathcal{P}_C(T)$, and set $L:=j^{\ast}\O_{\P_C\times T}(1)$.  For each $i$,
the identity $j\circ\sigma_i=q_i\times\mathrm{id}_T$ gives the canonical
trivialization
$$\gamma_i:\sigma_i^{\ast}L=(j\circ\sigma_i)^{\ast}\O_{\P_C\times T}(1)
=(q_i\times\mathrm{id}_T)^{\ast}\O_{\P_C\times T}(1)
\xrightarrow{(\tau_i)_T}\O_T.$$  
By Lemma~\ref{lem:fibers-goppa-restriction-isomorphism}, we have a canonical isomorphism $\rho:\underline{C}_T\xrightarrow{\sim}\pi_{\ast}L$.  We set
$$G(\mathfrak{p})=(\pi,\sigma_1,\dots,\sigma_n,L,\gamma_1,\dots,\gamma_n)$$ 
and define $\alpha_{G(\mathfrak p)}:=\rho^{-1}$.

\begin{proposition}\label{prop:G-functor}
The assignment $\mathfrak{p}\mapsto G(\mathfrak{p})$ defines a functor $G:\mathcal{P}_C\to\mathcal{G}_C$.
\end{proposition}

\begin{proof}
We first prove that $G(\mathfrak p)$ lies in $\mathcal{G}_C(T)$, that is, that
its Goppa code is $\underline C_T$.  
Let $\mathrm{incl}\colon\underline C_T\hookrightarrow\O_T^n$ be the inclusion. Since $(\pi_{\ast}\gamma)_i=\gamma_i\circ \ev_{\sigma_i}$, if the identities $\gamma_i\circ \ev_{\sigma_i}\circ\rho=\ell_i$, $i=1,\dots,n$, hold, then $\pi_{\ast}\gamma\circ\rho=\mathrm{incl}$. Since $\rho$ is an isomorphism, it would follow that $\im(\pi_{\ast}\gamma)=\underline C_T$. Thus it suffices to prove these identities.

Let $u_T:p_T^{\ast}\underline C_T\to\O_{\P_C\times T}(1)$ be the tautological quotient.  As $p_T\circ j=\pi$, the pullback $j^{\ast}u_T$ is a quotient
$\pi^{\ast}\underline C_T\to L$, and $\rho$ is by construction its adjoint.
Restriction along a section is compatible with this adjunction, so $\ev_{\sigma_i} \circ\rho=\sigma_i^{\ast}j^{\ast}u_T=(j\circ\sigma_i)^{\ast}u_T
=(q_i\times\mathrm{id}_T)^{\ast}u_T$, where the last equality holds because $\mathfrak p$ is
an object of $\mathcal{P}_C(T)$. It follows that $\gamma_i\circ \ev_{\sigma_i}\circ\rho
=(\tau_i)_T\circ(q_i\times\mathrm{id}_T)^{\ast}u_T$ as morphisms $\underline C_T\to\O_T$.  The right-hand side is the base change along $T\to\Spec\k$ of $\tau_i\circ q_i^{\ast}u_C$, which equals $\ell_i$ by
Lemma~\ref{lem:canonical-marked-projective}.

Now let $(u,f):\mathfrak p_1\to\mathfrak p_2$ be a morphism in $\mathcal{P}_C$
and write $L_a:=j_a^{\ast}\O_{\P_C\times T_a}(1)$ for $a=1,2$.  
Pulling back $\O_{\P_C\times T_2}(1)$ along the morphisms involved in the equality $j_2\circ f=(\mathrm{id}_{\P_C}\times u)\circ j_1$ and using the canonical
identification $(\mathrm{id}_{\P_C}\times u)^{\ast}\O_{\P_C\times T_2}(1)
\simeq\O_{\P_C\times T_1}(1)$
yields a canonical isomorphism $\varphi:L_1\xrightarrow{\sim}f^{\ast}L_2$.
Restricting that chain of
identifications along $\sigma^1_i$, and using $f\circ\sigma^1_i =\sigma^2_i\circ u$ together with $(\tau_i)_{T_1}=u^{\ast}(\tau_i)_{T_2}$, we
obtain $\gamma^1_i=u^{\ast}\gamma^2_i\circ(\sigma^1_i)^{\ast}\varphi$.  This
implies that $(u,f,\varphi)$ is a morphism of level structures $G(\mathfrak p_1)\to G(\mathfrak p_2)$, and the construction is compatible with
identities and composition.
\end{proof}

\paragraph{\textbf{Proof of Theorem~\ref{thm:fibers-goppa-projective-moduli}}}

By Propositions \ref{prop:G-functor} and \ref{prop:F-goppa-functor}, it only remains to prove that \(F\) and \(G\) are quasi-inverse equivalences.

Let $\gs=(\pi:X\to T,\sigma_1,\dots,\sigma_n,L,\gamma_1,\dots,\gamma_n)$
be an object of $\mathcal{G}_C(T)$.  Applying $G$ to $F(\gs)=(\pi,\sigma_1,\dots,\sigma_n,j_\gs)$ produces the
line bundle $L':=j_\gs^{\ast}\O_{\P_C\times T}(1)$.  For each $i$, let $\gamma_i':\sigma_i^{\ast}L'\to\O_T$ be the composite
$$\sigma_i^{\ast}L'=(j_\gs\circ\sigma_i)^{\ast}
\O_{\P_C\times T}(1)=(q_i\times\mathrm{id}_T)^{\ast}
\O_{\P_C\times T}(1)\xrightarrow{(\tau_i)_T}\O_T.$$  
The construction of $F$ also gives $\phi_\gs:L'\xrightarrow{\sim}L$.  We show that $\phi_\gs$ is an isomorphism $G(F(\gs))\xrightarrow{\sim}\gs$ in $\mathcal{G}_C(T)$.

It is enough to prove $\gamma_i=\gamma_i'\circ\sigma_i^{\ast}\phi_\gs^{-1}$ for every $i$.
Note that both sides are morphisms $\sigma_i^{\ast}L\to\O_T$.  
Since $\sigma_i^{\ast}q_\gs:\underline{C}_T\to\sigma_i^{\ast}L$ is
surjective, it suffices to prove the $i$-th equality after precomposing with $\sigma_i^{\ast}q_\gs$. On the one hand, $\gamma_i\circ\sigma_i^{\ast}q_\gs
=\gamma_i\circ \ev_{\sigma_i} \circ\alpha_\gs^{-1}=\ell_i$ by
\eqref{eq:auxiliary-identity}.  On the other hand, since $q_\gs=\phi_\gs\circ q'_\gs$, we have $\sigma_i^{\ast}\phi_\gs^{-1}\circ\sigma_i^{\ast}q_\gs =\sigma_i^{\ast}q'_\gs$, and so 
$(\gamma_i'\circ\sigma_i^{\ast}\phi_\gs^{-1})\circ \sigma_i^{\ast}q_\gs=(\tau_i)_T\circ\sigma_i^{\ast}q'_\gs =\ell_i$.  The last equality is the same computation as in the
proof of Proposition~\ref{prop:F-goppa-functor}.  Thus $\gamma_i=\gamma_i'\circ\sigma_i^{\ast}\phi_\gs^{-1}$, and $G\circ F\simeq\mathrm{id}_{\mathcal{G}_C}$.

Conversely, let $\mathfrak p=(\pi:X\to T,\sigma_1,\dots,\sigma_n,j)$ be an
object of $\mathcal{P}_C(T)$.  Let us set $L:=j^{\ast}\O_{\P_C\times T}(1)$,
and let $q_j:\pi^{\ast}\underline{C}_T\to L$ be the quotient obtained by pulling back
the tautological quotient on $\P_C\times T$.  By definition, $\rho:\underline{C}_T\xrightarrow{\sim}\pi_{\ast}L$ is the adjoint of $q_j$.
Applying $F$ to $G(\mathfrak p)$ gives the quotient $\varepsilon\circ \pi^{\ast}\rho:\pi^{\ast}\underline{C}_T\rightarrow \pi^{\ast}\pi_{\ast}L
\rightarrow L$, where $\varepsilon$ is the adjunction
map.  Since $\rho$ is the adjoint of $q_j$, we have $\varepsilon \circ\pi^{\ast}\rho=q_j$.  Therefore the quotient defining $F(G(\mathfrak p))$ is exactly $q_j$.  By the universal property of the
projective bundle \cite[Tag~01OA]{stacks}, the morphism to $\P_C\times T$
determined by this quotient is $j$.  It follows that $F(G(\mathfrak p))=\mathfrak p$,
functorially in $T$.

The two natural isomorphisms constructed above are compatible with pullback
along morphisms $T'\to T$. Consequently, $F$ and $G$ are quasi-inverse
equivalences of
categories fibered in groupoids, which proves Theorem~\ref{thm:fibers-goppa-projective-moduli}.

\subsection{Representability of the Goppa morphism}\label{subsec:goppa-representable}

The description of
Theorem~\ref{thm:fibers-goppa-projective-moduli} admits a relative form. Let $S$ be a $\k$-scheme and let $\mathcal E\subset\O_S^n$ be a rank-$k$
subbundle such that every coordinate projection $\ell_i\colon\mathcal E\to\O_S$ is surjective. As in Lemma~\ref{lem:canonical-marked-projective}, $\ell_i$ determines a section $q_i\colon S\to\P_S(\mathcal E)$ and a canonical trivialization of $\O(1)$
along $q_i$. Replacing $\underline C_T$ by $\mathcal E_T$ throughout the proof
of Theorem~\ref{thm:fibers-goppa-projective-moduli} gives, for $d\geq2g+1$, a canonical equivalence between $\LS_{g,n,d}\times_{\Gr(k,n)}S$ and the category of families of $n$-pointed
smooth projective curves $X\to T$, for $T\to S$, endowed with a closed immersion $j\colon X\hookrightarrow\P_T(\mathcal E_T)$ of relative degree $d$ satisfying $j\circ\sigma_i=q_{i,T}$ and whose geometric fibers $j_t$ are non-degenerate. The equivalence is compatible with
arbitrary base change on $S$.

In a fixed projective space, incidence schemes of curves and points are classically described by flag Hilbert schemes. Perrin treats this flag Hilbert construction in \cite[\S~1, Example~1.1]{perrin-points}, and notes in the Appendix that the results of \S~1 are characteristic-free. In our relative setting, representability follows from Grothendieck's relative Hilbert scheme \cite{grothendieck-hilbert,nitsure-hilbert}. 
The condition that a family contain the prescribed sections is closed, while smoothness and geometric
integrality are open. Since $d>2g-2$, non-degeneracy is also open, being the non-vanishing locus of the determinant of the restriction morphism $\mathcal E_T\to\pi_*\O_X(1)$.

\begin{theorem}\label{thm:goppa-representable}
Assume $n>d\geq2g+1$. Then $\Goppa_{g,n,d}\colon\LS_{g,n,d}\to\Gr(k,n)$ is schematic and
quasi-projective. In particular, $\LS_{g,n,d}$ is represented by a quasi-projective $\k$-scheme.
\end{theorem}

\begin{proof}
Let $G^\circ\subset\Gr(k,n)$ be the open locus where all coordinate projections of the universal subbundle are surjective and the induced
sections of its projective bundle are pairwise disjoint. Since $d\geq2g+1$, the complete linear series of every level structure is very ample, so the Goppa morphism factors through $G^\circ$.

By the relative description above, over $G^\circ$ the stack $\LS_{g,n,d}$ is identified with the locus in $\Hilb^{dm+1-g}_{\P(\mathcal U)/G^\circ}$ parametrizing smooth geometrically
integral, non-degenerate curves containing the prescribed sections.
By the preceding discussion this is a locally closed subscheme of the relative Hilbert scheme, hence is quasi-projective over $G^\circ$. The result follows.
\end{proof}

\begin{remark}\label{rmk:representable-vs-DM}
The Hilbert-scheme step in Theorem~\ref{thm:goppa-representable} is classical. The Goppa-specific input is the relative projective description above. In the range common to Theorems~\ref{thm:goppa-representable}
and~\ref{thm:structure-LS}, it follows that the Deligne--Mumford stack $\LS_{g,n,d}$ is in fact a quasi-projective scheme.
\end{remark}

\subsection{Universal injectivity of the Goppa morphism}

The description of the fibers yields the universal injectivity of the Goppa
morphism in the range $n/2>d>2g+1$. Its geometric core is the following uniqueness statement for embedded
curves through assigned points, whose proof follows the argument of \cite[Theorem~2]{marquez4}.

\begin{proposition}\label{prop:unique-embedded-curve}
Let $\Omega$ be an algebraically closed field and let $g,n,d$ be integers
with $n>2d$ and $d>2g+1$. Set $k:=1-g+d$. Let $q_{1},\dots,q_{n}\in\P^{k-1}(\Omega)$ be pairwise distinct points. If $Y,Y'\subset\P^{k-1}_{\Omega}$ are smooth, irreducible, non-degenerate
curves of genus $g$ and degree $d$ with $q_{i}\in Y\cap Y'$ for all $i=1,\dots,n$, then $Y=Y'$.
\end{proposition}

\begin{proof}
Since $d>2g-2$, Serre duality and Riemann--Roch give $h^{0}(Y,\O_{Y}(1))=d+1-g=k$, so by non-degeneracy the restriction map $H^{0}(\P^{k-1},\O(1))\to H^{0}(Y,\O_{Y}(1))$ is an isomorphism, so $Y$ is
embedded by the complete linear series of an invertible sheaf of degree $d\geq 2g+2$, and its homogeneous ideal is therefore generated by
quadrics~\cite{saint-donat-equations,fujita-defining} (see also \cite[Proposition~11]{marquez4}). In
particular, $Y$ is the set-theoretic intersection of the quadrics containing
it. Let $Q$ be such a quadric. If $Y'\not\subseteq Q$, then $Q\cap Y'$ is
the zero scheme of a nonzero section of $\O_{Y'}(2)$, an effective divisor of degree $2d$ on the smooth curve $Y'$, so $\#(Y'\cap Q)\leq 2d<n$,
contradicting $q_{1},\dots,q_{n}\in Y'\cap Q$
(cf.~\cite[Proposition~14]{marquez4}). Hence $Y'\subseteq\bigcap_{Q\supseteq Y}Q=Y$, and $Y'=Y$ because both are
integral curves.
\end{proof}

\begin{theorem}\label{thm:universal-injectivity}
If $n/2>d>2g+1$, then, for every field extension $\k\subset K$, two objects $\gs,\gs'\in\LS_{g,n,d}(K)$ with the same Goppa code are isomorphic. In particular, $\Goppa_{g,n,d}:\LS_{g,n,d}\to\Gr(k,n)$ is universally injective.
\end{theorem}

\begin{proof}
By Theorem~\ref{thm:goppa-representable}, $\LS_{g,n,d}$ is represented by a scheme, so its $K$-points correspond to isomorphism classes of objects of $\LS_{g,n,d}(K)$. 
Hence the first assertion is equivalent to universal injectivity of the induced morphism between the representing scheme and $\Gr(k,n)$, by \cite[Tag~01S4]{stacks}. By \cite[Tag~0CPN]{stacks}, it therefore suffices to prove that two objects over an algebraically closed field with the same image are isomorphic.

Let us fix an algebraically closed field $\Omega$ over $\k$ and a code $C\in\Gr(k,n)(\Omega)$ in the image of $\Goppa_{g,n,d}$. It is non-degenerate because $n>d>2g-1$ (see \S~\ref{subsec:density}), so
Theorem~\ref{thm:fibers-goppa-projective-moduli}, applied over the base field $\Omega$, identifies the fiber $\mathcal{G}_{C}(\Omega)$ with $\mathcal{P}_{C}(\Omega)$. It is therefore enough to prove that any two objects of $\mathcal{P}_{C}(\Omega)$ are isomorphic.

Let $(X,\sigma_{1},\dots,\sigma_{n},j)$ and $(X',\sigma'_{1},\dots,\sigma'_{n},j')$ be two such objects, and set $Y:=j(X)$, $Y':=j'(X')$. The points $q_{1},\dots,q_{n}$ of
Lemma~\ref{lem:canonical-marked-projective} are pairwise distinct, since the
sections $\sigma_{i}$ are pairwise disjoint and $j$ is a monomorphism. Moreover, $Y,Y'\subset\P_{C}$ are smooth, irreducible, non-degenerate curves of genus $g$ and degree $d$ containing them. As $n>2d$ and $d>2g+1$,
Proposition~\ref{prop:unique-embedded-curve} gives $Y=Y'$ as closed
subschemes of $\P_{C}$. Since $j$ and $j'$ are isomorphisms onto $Y$, the
morphism $f:=j'^{-1}\circ j\colon X\to X'$ is an isomorphism with $j'\circ f=j$, and $f\circ\sigma_{i}=\sigma'_{i}$ because $j'\circ f\circ\sigma_{i}=q_{i}=j'\circ\sigma'_{i}$ and $j'$ is a
monomorphism. By Definition~\ref{def:P_C}, $f$ is an isomorphism in $\mathcal{P}_{C}(\Omega)$.
\end{proof}

\begin{remark}
The key geometric ingredient of the proof of Theorem~\ref{thm:universal-injectivity} is the method of \cite[Theorem~2]{marquez4}, resting on the classical results on quadrics
cited in Proposition~\ref{prop:unique-embedded-curve}. The level-structure
formalism replaces the divisor-equivalence theorem of Munuera and Pellikaan~\cite[Corollary~4.15]{munuera-pellikaan} invoked there, 
and the first assertion of Theorem~\ref{thm:universal-injectivity} extends the
uniqueness part of \cite[Theorem~2]{marquez4} to arbitrary (in particular non-perfect) base fields, such as $\Fq(z)$, a case not covered in the literature. This might be relevant in (convolutional-) code-based
cryptography, where the geometric structure of Goppa codes is exploited by structural distinguishers and key-recovery attacks~\cite{couvreur-mcc-ag}.
\end{remark}

\subsection{An example in genus one}

Throughout this subsection $C\in\Gr(k,n)(\k)$ is a non-degenerate code,
and $q_1,\dots,q_n\in\P_C(\k)$ are the distinguished points
attached to $C$ by Lemma~\ref{lem:canonical-marked-projective}. By
Theorem~\ref{thm:fibers-goppa-projective-moduli}, whenever $n>d\geq 2g+1$ the
fiber $\mathcal{G}_C$ is canonically equivalent to the groupoid of smooth
$n$-pointed curves of genus $g$ and degree $d$ embedded in $\P_C$, with
non-degenerate geometric fibers, and passing through $q_1,\dots,q_n$.

If $q_i=q_j$ for some $i\neq j$, then $\mathcal{G}_C$ is empty by
Definition~\ref{def:P_C}, since the marked sections of an $n$-pointed smooth
curve are pairwise disjoint and the immersion $j$ is a monomorphism. In every
non-empty case the points $q_1,\dots,q_n$ are therefore pairwise distinct. The
genus-zero case is described in Proposition~\ref{prop:fibers-g0}. We
illustrate the general description in genus one, where the answer is classical.
The classification below is geometric, and may be checked after extension of
scalars to an algebraic closure of $\k$.

Assume $g=1$ and $d=3$.
Here $k=3$, so $\P_C\simeq\P^2$, and the condition $n>d=3$ requires $n\geq 4$.
Theorem~\ref{thm:fibers-goppa-projective-moduli} identifies $\mathcal{G}_C$
with the groupoid of smooth plane cubics passing through $q_1,\dots,q_n$. The
genus formula for smooth plane curves gives $(3-1)(3-2)/2=1$, so these are
precisely the smooth genus-one plane cubics through the distinguished points.

Let $Z:=q_1+\cdots+q_n\subset\P^2$ be the reduced zero-dimensional subscheme
defined by the distinguished points, and set $V_Z:=H^0(\P^2,\mathcal{I}_Z(3))$.
Then $\mathcal{G}_C$ is represented by the open subscheme of $\P(V_Z^\vee)$
parametrizing smooth cubics. In particular, if the points impose $n$ independent
conditions on cubic forms, then $\dim V_Z=10-n$, and the cubics through them
form a projective space $\P^{9-n}$. For $n=8$ this is a pencil, and
$\mathcal{G}_C$ is its open subscheme of smooth members.

For $n=9$ there are two basic possibilities. If
$\dim H^0(\P^2,\mathcal{I}_Z(3))=1$, there is a unique cubic through the nine points,
and the fiber is either empty or $\Spec\k$ according as this cubic is not
smooth or smooth. If $Z$ is the scheme-theoretic complete intersection of
two cubics, then $\dim H^0(\P^2,\mathcal{I}_Z(3))=2$, the cubics through $Z$ form a
pencil, and by the Cayley--Bacharach theorem
\cite[Theorem~CB4]{eisenbud-green-harris} (applied over an algebraic closure of
$\k$) any cubic through eight of the nine base points passes through the ninth.
More special configurations may give larger linear systems or fixed components.
In every case, $\mathcal{G}_C$ is the open subscheme of $\P(V_Z^\vee)$ parametrizing
smooth members.

\section{Involutions and self-duality}\label{sec:duality}

In this section we study self-dual Goppa codes from the perspective of level structures. Recall that a linear code $C\subset\k^n$ of dimension $k$ is self-dual if $C=C^{\bot}$. As a consequence, every self-dual code has length $n=2k$.

Stichtenoth~\cite{stichtenoth-selfdual} gave sufficient conditions for
self-duality in terms of differentials, and Driencourt--Stichtenoth \cite{driencourt-stichtenoth} obtained a converse criterion under additional
hypotheses. Over finite fields, Munuera and Pellikaan gave a
characterization in the range relevant here \cite[Theorem~5.5]{munuera-pellikaan}. The
level-structure viewpoint recovers these differential criteria and realizes the self-dual locus as the fixed-point locus of an involution~$\bot$.

Throughout this section let $(\pi:X\to S,\sigma_1,\dots,\sigma_n)$ be a family of $n$-pointed smooth projective curves of genus $g$, and set $D:=\sigma_1+\cdots+\sigma_n$.

The construction of \S~\ref{sec:aremark}, carried out on $T$-valued points, is
compatible with base change. Hence tensor product defines morphisms $\LS_{X,D,e}\times_S\LS_{X,D,e'}\to\LS_{X,D,e+e'}$ and inversion defines
morphisms $\LS_{X,D,e}\to\LS_{X,D,-e}$, for $e,e'\in\mathbb Z$. In particular, $\LS_{X,D,0}$ is a commutative group scheme over $S$
(Theorem~\ref{th:relativestructure}). 
The neutral element of $\LS_{X,D,0}$ is
\(
\mathbf{1}
=
\bigl(\O_X,\gamma^{\mathrm{triv}}_1,\dots,\gamma^{\mathrm{triv}}_n\bigr),
\)
where, for each $i=1,\dots,n$, 
\(
\gamma^{\mathrm{triv}}_i:
\sigma_i^*\O_X
\xrightarrow{\sim}\O_S
\)
is the canonical  identification. Over a field $K$, $\LS_{X,D,0}(K)$ is the group of degree-zero line bundles on $X_K$ rigidified along $D_K$, modulo isomorphism. When $K$ is finite, its divisor
description is precisely the degree-zero part $\PicG^0(X_K,D_K)$ of the generalized Picard group used by Munuera and Pellikaan \cite[\S~7]{munuera-pellikaan}.

Let $T$ be an $S$-scheme. 
Let $\Res_i\colon\sigma_{iT}^{*}\omega_{X_T/T}(D_T)\simeq\O_T$
be the relative residue isomorphism. If $t$ is a local
equation of the divisor $\sigma_{iT}(T)$ and $a$ is a local section of $\O_{X_T}$, then $\Res_i(a\,dt/t)=a|_{\sigma_{iT}}$, which does not depend on $t$. These
isomorphisms are compatible with base change and involve no assumption on the
characteristic. They define the \emph{canonical level structure} $\gs_{\mathrm{can}}:=(\omega_{X/S}(D),\Res_1,\dots,\Res_n)\in
\LS_{X,D,2g-2+n}(S)$.

From now on we assume $n=2(1-g+d)$, so that $2g-2+n=2d$. For $\gs=(L,\gamma_1,\dots,\gamma_n)\in\LS_{X,D,d}(T)$ we set $\gs^{\bot}:=\gs_{\mathrm{can},T}\otimes\gs^{-1}\in\LS_{X,D,d}(T)$. Explicitly,
$\gs^{\bot}=(M,\eta_1,\dots,\eta_n)$, where $M=\omega_{X_T/T}(D_T)\otimes L^{-1}$ and $\eta_i=\Res_i\otimes\gamma_i^{-1}$,
with $\gamma_i^{-1}$ as in \S~\ref{sec:aremark}. This defines an $S$-morphism $\bot\colon\LS_{X,D,d}\to\LS_{X,D,d}$ with $\bot\circ\bot=\mathrm{id}$. For every
field-valued point $\Spec K\to S$ and every $\gs\in\LS_{X,D,d}(K)$, the
evaluation map of $\gs^{\bot}$ is the map $\Res_{\gs}$ of \S~\ref{subsec:dual-code}. Hence $C_{\gs^{\bot}}=C_{\gs}^{\bot}$.

\begin{definition}
A level structure $\gs\in\LS_{X,D,d}(T)$ is \emph{self-dual} if $\gs^{\bot}=\gs$. The \emph{self-dual locus} is the fixed-point subscheme $\fix(\bot):=(\bot,\mathrm{id})^{-1}(\Delta_{\LS_{X,D,d}/S})$. It is closed in $\LS_{X,D,d}$ because $\LS_{X,D,d}\to S$ is separated (Theorem~\ref{th:relativestructure}).
\end{definition}

\begin{lemma}\label{lm:char-dual}
For every $S$-scheme $T$ and every $\gs\in\LS_{X,D,d}(T)$, one has $\gs^{\bot}=\gs$ if and only if $\gs^{\otimes2}=\gs_{\mathrm{can},T}$.
\end{lemma}

\begin{proof}
Tensor the identity $\gs^{\bot}=\gs_{\mathrm{can},T}\otimes\gs^{-1}$ with $\gs$.
\end{proof}

\begin{lemma}\label{lem:selfdual-level-code}
Assume $d\geq2g+1$. Let $\Spec K\to S$ be a field-valued point and let $\gs\in\LS_{X,D,d}(K)$. Then $\gs$ is self-dual if and only if $C_{\gs}$ is self-dual.
\end{lemma}

\begin{proof}
Since $n-d=d-2g+2\geq3$ and $2d>3g$, one has $(n-d)(d-g)\geq3(d-g)>d$. Hence $\Phi_{X,D,d}$ is a monomorphism by Lemma~\ref{alvarez-porras}, and $\gs^{\bot}=\gs$ if and only if $C_{\gs^{\bot}}=C_{\gs}$, that is, $C_{\gs}^{\bot}=C_{\gs}$.
\end{proof}

\begin{lemma}\label{lm:LS0-structure}
\begin{enumerate}
\item The morphism $\iota\colon H_S\to\LS_{X,D,0}$, $\lambda\mapsto\lambda\cdot\mathbf 1$, and $o$ form an exact sequence of fppf sheaves of abelian groups $1\to H_S\to\LS_{X,D,0}\to\Pic_{X,D,0}\to1$.
\item The morphism $[2]$ of $\LS_{X,D,0}$ is an fppf epimorphism, and there is an exact sequence of fppf sheaves $1\to\mu_{2,S}^{n-1}\to\LS_{X,D,0}[2]\to\Pic_{X,D,0}[2]\to1$.
\item $\LS_{X,D,0}[2]$ is finite locally free of rank $2^{2g+n-1}$ over $S$, and \'etale if $\operatorname{char}\k\neq2$.
\end{enumerate}
\end{lemma}

\begin{proof}
(1) Since $\lambda\cdot\gs=\iota(\lambda)\otimes\gs$ for $\lambda\in H_S(T)$ and $\gs\in\LS_{X,D,0}(T)$, this is a reformulation of
Lemma~\ref{lm:preliminariesraynaud}(1) and (2).

(2) The $S$-scheme $\Pic_{X,D,0}$ is smooth and proper (see the discussion preceding Lemma~\ref{lm:preliminariesraynaud}) with connected geometric
fibers. Hence it is an abelian scheme, and its morphism $[2]$ is finite locally free of rank $2^{2g}$ \cite[Proposition~20.7]{milne-abelian-varieties}. The morphism $[2]$ of $H_S\simeq\mathbb G_{m,S}^{n-1}$ is finite locally free of rank $2^{n-1}$ with kernel $\mu_{2,S}^{n-1}$. Both are fppf epimorphisms, and (2) follows from the
snake lemma applied to $[2]$ on the sequence in (1).

(3) By (2), $\LS_{X,D,0}[2]\to\Pic_{X,D,0}[2]$ is an fppf torsor under $\mu_{2,S}^{n-1}$. Hence it is finite locally free of rank $2^{n-1}$ and \'etale if $\operatorname{char}\k\neq2$, these properties being local on the base for the fppf topology \cite[2.5.1, 2.7.1]{ega4-2}, \cite[17.7.3]{ega4-4}. The assertion follows by composing with $\Pic_{X,D,0}[2]\to S$ and using (2).
\end{proof}

\begin{theorem}\label{prop:selfdual-torsor}
Assume $n=2(1-g+d)$. Then $\fix(\bot)$ is an $\LS_{X,D,0}[2]$-torsor over $S$ for the fppf topology, the action being given by tensor product. In particular, $\fix(\bot)\to S$ is finite locally free of rank $2^{2d+1}$ and \'etale if $\operatorname{char}\k\neq2$, and every $\gs_0\in\fix(\bot)(S)$ induces an isomorphism $\LS_{X,D,0}[2]\to\fix(\bot)$, $\tau\mapsto\tau\otimes\gs_0$.
\end{theorem}

\begin{proof}
By Lemma~\ref{lm:char-dual}, $\fix(\bot)(T)$ is the set of $\gs\in\LS_{X,D,d}(T)$ with $\gs^{\otimes2}=\gs_{\mathrm{can},T}$. Hence $\LS_{X,D,0}[2]$ preserves $\fix(\bot)$. Moreover, for $\gs,\gs'\in\fix(\bot)(T)$ the element $\gs\otimes\gs'^{-1}$ is the unique $\tau\in\LS_{X,D,0}[2](T)$ with $\tau\otimes\gs'=\gs$. It remains to find
sections fppf-locally on $S$. Zariski-locally on $S$ the sheaves $\sigma_i^{*}\O_X(d\sigma_1)$ are trivial, which gives $\xi_0\in\LS_{X,D,d}(U)$ on each member $U$ of an open covering of $S$. By
Lemma~\ref{lm:LS0-structure}(2), after an fppf covering there is $\tau\in\LS_{X,D,0}$ with $\tau^{\otimes2}=\gs_{\mathrm{can}}\otimes(\xi_0^{\otimes2})^{-1}$, and then $\tau\otimes\xi_0\in\fix(\bot)$. The remaining assertions follow from Lemma~\ref{lm:LS0-structure}(3), since $\fix(\bot)$ is fppf-locally isomorphic
to $\LS_{X,D,0}[2]$, these properties are local on the base for the fppf topology, and $2g+n-1=2d+1$.
\end{proof}

\begin{remark}\label{rmk:selfdual}
Let $S=\Spec\k$ and assume $d\geq2g+1$.

(1) Let $G$ be a divisor of degree $d$ on $X$ whose support contains no $p_i$, let $v\in(\k^{*})^{n}$, and let $\gs$ be the level structure on $\O_X(G)$ with trivializations $f\mapsto v_if(p_i)$. Every element of $\LS_{X,D,d}(\k)$ is of this form, because invertible sheaves on $X$ are free on the semilocal ring of $X$ at $p_1,\dots,p_n$. By Lemma~\ref{lm:char-dual} and the identity $\Res_{p_i}(hw)=h(p_i)\Res_{p_i}(w)$ for differentials $w$ with at most simple
poles at $p_i$, the structure $\gs$ is self-dual if and only if there is a
differential $w$ with $(w)=2G-D$ and $\Res_{p_i}(w)=v_i^{2}$ for all $i$. By Lemma~\ref{lem:selfdual-level-code}, this is the level-structure form
of the classical differential criterion for the self-duality of $C_{\gs}$.
Over finite fields, in the present numerical range, a necessary-and-sufficient formulation is given in
\cite[Theorem~5.5]{munuera-pellikaan}. See also
\cite{stichtenoth-selfdual,driencourt-stichtenoth}.

(2) By Theorem~\ref{prop:selfdual-torsor} and
Lemma~\ref{lem:selfdual-level-code}, the set of self-dual codes $C_{\gs}$ with $\gs\in\LS_{X,D,d}(\k)$ is either empty or a torsor under the finite group $\LS_{X,D,0}[2](\k)$. Assume that $\k$ is perfect of characteristic $2$. Then $\mu_2(\k)=\{1\}$, $H_S(\k)$ is $2$-divisible, and $\LS_{X,D,0}(\k)\to\Pic_{X,D,0}(\k)$ is surjective because $X$ has rational
points \cite[Theorem~9.2.5]{kleiman-picard}. Hence $\LS_{X,D,0}[2](\k)\simeq\Pic_{X,D,0}[2](\k)$ by Lemma~\ref{lm:LS0-structure}(1). Indeed, any element of
$\Pic_{X,D,0}[2](\k)$ lifts to $b\in\LS_{X,D,0}(\k)$. Since $b^{\otimes2}\in H_S(\k)$ and $H_S(\k)$ is $2$-divisible, multiplying $b$
by a suitable element of $H_S(\k)$ gives a $2$-torsion lift, while
injectivity follows from $\mu_2(\k)=\{1\}$. Since every element of $\k^{*}$ is a square,
(1) shows that such codes exist if and only if $2G-D$ is a canonical divisor
for some divisor $G$ on $X$. Compare
\cite{stichtenoth-selfdual,munuera-pellikaan}.
\end{remark}

\section{The genus zero case}\label{sec:rational}

Let $\k$ be an arbitrary base field and let $n,d$ be integers with $n\geq 3$
and $n>d\geq 1$.

\medskip
The genus-zero case is exceptional in the theory of Goppa codes for two reasons. First, every smooth projective curve of genus zero over $\k$ with at least three rational points is isomorphic to $\mathbb{P}^1_\k$, and $\PGL_2$ acts simply transitively on ordered triples of distinct rational points.  This rigidity gives an explicit description of the moduli space, $\LS_{0,n,d}\simeq\Ht\times\M_{0,n}$, where $\Ht$ encodes the
trivialization data and $\M_{0,n}$ is the usual moduli space of $n$-pointed
rational curves. Second, the sheaf $\O(dp_1)$ on $\mathbb{P}^1_\k$ has an explicit basis of global sections, namely the polynomials of bounded
degree in the affine coordinate, so that all evaluation maps, generator matrices and parity-check matrices can be written in closed form in terms of the marked points and the trivialization parameters.

These two features make the genus-zero theory completely explicit. The main results of the section are Proposition~\ref{prop:LS0nd-explicit}, which
identifies $\LS_{0,n,d}$ with $\Ht\times\M_{0,n}$ as schemes, and the description of the Goppa morphism in coordinates (Subsection~\ref{subsec:space-and-goppa}). We specialize
Theorem~\ref{thm:fibers-goppa-projective-moduli} to the genus-zero setting, thereby describing the
fibers of the Goppa morphism in terms of rational normal curves, and we deduce
that in the range $2\leq d\leq n-3$ the Goppa morphism is itself an immersion (Theorem~\ref{thm:goppa-immersion-g0},
Subsection~\ref{subsec:goppa-immersion-g0}). We also characterize self-dual codes by an explicit residue condition (Proposition~\ref{prop:selfdual-g0}).

 \subsection{Fundamental calculations}\label{subsec:fund}
 \subsubsection{The pointed curve}
Consider the curve $\mathbb{P}^{1}_{\k}=\Proj(\k[x,y])$ together with $n$ rational points 
\begin{equation}
p_1=\infty:=[1:0], p_{2}:=[0:1], p_3:=[1:1], p_4=[\alpha_4:1],\dots, p_n=[\alpha_{n}:1],
\end{equation}
where 
\[
(\alpha_4,\dots,\alpha_n)\in (\mathbb{P}^{1}_{\k}\setminus\{\infty,0,1\})^{n-3}\setminus\Delta,
\]
where $\Delta$ denotes the union of the diagonals. Note that $p_i\ne p_j$ for $i\ne j$. Set $\alpha_2:=0$, $\alpha_3:=1$ and $D:=p_1+\cdots+p_n$.
 \subsubsection{The invertible sheaf}

Let us now consider the invertible sheaf of degree $d$ given by 
 $L=\O_{\mathbb{P}^{1}_{\k}}(d p_1)$.  Let us set $t=xy^{-1}$, so that $\k(t)$ is the field of rational functions of $\mathbb{P}^{1}_{\k}$. One has
\begin{equation}\label{eq:basis}
 H^{0}(\mathbb{P}^{1}_{\k},L)=\k[t]_{\leq d}=\langle 1,t,\dots,t^d \rangle.
 \end{equation}
 \subsubsection{The trivializations}\label{subsec:trivializations}
By \eqref{eq:canonicaltrivializations}, the canonical trivializations are
\begin{equation*}
\begin{split}
 \gamma^0_i:L|_{p_i}&\simeq \k, \ i=2,\dots,n\\
 f(t)&\mapsto f(\alpha_i).
\end{split}
\end{equation*}
Take $U:=\mathbb{P}^{1}_{\k}\setminus\{0\}$ and, for each $f(t)\in L(U)$, set $g(u):=u^{d}f(u^{-1})$, where $u=t^{-1}$ is the local parameter at $\infty$. The function $g(u)$ is regular on $U$, and $f(t)\mapsto g(u)$ defines an isomorphism
\begin{equation*}
\begin{split}
L(U)&\xrightarrow{\;\sim\;} \O_{\mathbb{P}^{1}_{\k}}(U)\\
f(t)&\mapsto g(u).
\end{split}
\end{equation*}
This isomorphism induces a trivialization
\begin{equation*}
\begin{split}
\gamma^0_{1}:L|_{\infty}&\simeq \k\\
f(t)&\mapsto g(0).
\end{split}
\end{equation*}
\subsubsection{The evaluation map and the code}
The data 
\begin{equation}
\gs^0=(\mathbb{P}^{1}_{\k},\infty,0,1,p_4,\dots,p_n, L, \gamma^0_1,\gamma^0_2,\dots,\gamma^0_n)
\end{equation}
determines a level structure. The corresponding evaluation map \eqref{eq:evaluationmapcurve} is
\begin{equation*}
\begin{split}
\ev_{\gs^0}:H^{0}(\mathbb{P}^{1}_{\k},L)&\rightarrow \k^{n}\\
f(t)&\mapsto (g(0), f(0), f(1),f(\alpha_4),\dots,f(\alpha_{n})).
\end{split}
\end{equation*}
When $\k=\Fq$, the image $C_{\gs^0}$ is a generalized Reed--Solomon code.
\subsubsection{The generator matrix}
With respect to the basis \eqref{eq:basis}, the matrix of this linear map is 
\begin{equation}
G_{\gs^{0}}=
\left(
\begin{array}{cccccc}
0 & 1 & 1 & 1 & \dots & 1\\
0 & 0 & 1 & \alpha_4 & \dots & \alpha_n\\
0 & 0 & 1 & \alpha_{4}^{2} & \dots & \alpha_{n}^{2}\\
\vdots & \vdots & \vdots & \vdots & \dots & \vdots \\
0 & 0 & 1 & \alpha_{4}^{d-1} & \dots & \alpha_{n}^{d-1}\\
1 & 0 & 1 & \alpha_{4}^{d} & \dots & \alpha_{n}^{d}
\end{array}
\right)\in\Mat_{1+d,n}(\k).
\end{equation}
Given a tuple $(l_1,\dots,l_{n-1})\in\mathbb{G}_m^{n-1}$, define new trivializations $\gamma_i$ by
$$
L|_{p_{i}}\overset{\gamma_{i}^{0}}{\simeq} \k\overset{l_{i}}{\simeq} \k, \ i=1,\dots, n-1.
$$
The new level structure
\begin{equation}\label{eq:composed_level_structure}
\gs=(\mathbb{P}^{1}_{\k},\infty,0,1,p_4,\dots,p_n, L, \gamma_1,\gamma_2,\dots,\gamma_n)
\end{equation}
determines a new code $C_{\gs}$ whose generator matrix is given by
\begin{equation}\label{eq:basic_generator}
G_{\gs}=
\left(
\begin{array}{cccccc}
0 & l_{2} & l_{3} & l_{4} & \dots & 1\\
0 & 0 & l_{3} & l_{4} \alpha_4 & \dots &  \alpha_n\\
0 & 0 & l_{3}  & l_{4} \alpha_{4}^{2} & \dots &  \alpha_{n}^{2}\\
\vdots & \vdots & \vdots & \vdots & \dots & \vdots \\
0 & 0 & l_{3}  & l_{4} \alpha_{4}^{d-1} & \dots &  \alpha_{n}^{d-1}\\
l_{1}  & 0 & l_{3}  & l_{4} \alpha_{4}^{d} & \dots &\alpha_{n}^{d}
\end{array}
\right)\in\Mat_{1+d,n}(\k).
\end{equation}

\subsubsection{The dual code and the parity-check matrix}
By \S~\ref{subsec:dual-code}, the dual code is the image of the following
residue map.
\begin{equation}
\begin{split}
\Res_{\gs^{0}}:H^{0}(\mathbb{P}^{1}_{\k},\omega_{\mathbb{P}^{1}_{\k}}(D)\otimes L^{-1})&\rightarrow \k^{n}\\
w&\mapsto (\Res_{\gs^{0},p_1}(w),\dots, \Res_{\gs^{0},p_n}(w)).
\end{split}
\end{equation}
A basis of $H^{0}(\mathbb{P}^{1}_{\k},\omega_{\mathbb{P}^{1}_{\k}}\otimes L^{-1}(D))$ is
$$
H^{0}(\mathbb{P}^{1}_{\k},\omega_{\mathbb{P}^{1}_{\k}}\otimes L^{-1}(D))=\left\langle \dfrac{1}{\prod_{i\neq 1} (t-\alpha_i)}dt,\dots, \dfrac{t^{n-d-2}}{\prod_{i\neq 1}(t-\alpha_i)}dt \right\rangle.
$$
Setting $w_{m}:=\frac{t^{m}}{\prod_{i\neq 1}(t-\alpha_i)}dt$, a direct residue computation gives
$$
\left\{
\begin{array}{ccccl}
\Res_{\gs^{0},p_j}(w_{m}) &=& \dfrac{\alpha_{j}^{m}}{\prod_{i\neq 1,j}(\alpha_{j}-\alpha_{i})}  & \text{if} & j=2,\dots,n\\
\Res_{\gs^{0},p_1}(w_{m}) &=& 0  &   \text{if}  & m<n-d-2\\
\Res_{\gs^{0},p_1}(w_{m}) &=& -1  &   \text{if}  & m=n-d-2
\end{array}
\right.
$$
\begin{remark}
To compute $\Res_{\gs^{0},p_1}(w_m)$, write $t=u^{-1}$. By the definition
of the induced residue trivialization in \S~\ref{subsec:dual-code}, at $p_1$
this amounts to multiplying by $u^{-d}$. The resulting local expression is
$$
-\dfrac{u^{-m-d+n-3}}{\prod_{i\neq 1}(1- u\alpha_i)}du,
$$
from which the stated values follow.
\end{remark}
The parity-check matrix of $C_{\gs^{0}}$ is
\begin{equation}
P_{\gs^{0}}=
\left(
\begin{array}{cccccc}
0                              & 0                            & \dots            &-1 \\
h_{2}       &  0 & \dots            & 0  \\
h_{3}       &  h_{3} & \dots            & h_{3}  \\
h_{4}       &  \alpha_{4} h_{4} & \dots            & \alpha_{4}^{n-d-2}h_{4}  \\
\vdots                      & \vdots                     &                      & \vdots\\
h_{n}       &  \alpha_{n}h_{n} & \dots            & \alpha_{n}^{n-d-2}h_{n}
\end{array}
\right)\in\Mat_{n,n-d-1}(\k),
\end{equation}
where $h_{j}=\frac{1}{\prod_{i\neq 1,j}(\alpha_{j}-\alpha_{i})}$. A direct computation gives $G_{\gs^{0}}P_{\gs^{0}}=0$.
For the level structure $\gs$ of \eqref{eq:composed_level_structure}, the parity-check matrix is
\begin{equation}\label{eq:basic_parity_check}
P_{\gs}=
\left(
\begin{array}{cccccc}
0                              & 0                            & \dots            &-l^{-1}_{1} \\
l^{-1}_{2}h_{2}       & 0& \dots            & 0 \\
l^{-1}_{3}h_{3}       & l^{-1}_{3} h_{3} & \dots            &  l^{-1}_{3} h_{3}  \\
l^{-1}_{4}h_{4}       & l^{-1}_{4} \alpha_{4}h_{4} & \dots            &  l^{-1}_{4} \alpha_{4}^{n-d-2}h_{4}  \\
\vdots                      & \vdots                     &                      & \vdots\\
 h_{n}       &\alpha_{n}h_{n} & \dots            &  \alpha_{n}^{n-d-2}h_{n}
\end{array}
\right)\in\Mat_{n,n-d-1}(\k).
\end{equation}

\subsection{The space of level structures and the Goppa morphism}\label{subsec:space-and-goppa}

Every family of $n$-pointed smooth projective curves of genus $0$ with $n\geq 3$, $X\rightarrow S$, parametrized by a $\k$-scheme $S$, is trivialized by the unique $S$-isomorphism $X\xrightarrow{\sim}\mathbb{P}^{1}_{S}$ sending any prescribed ordered triple of pairwise disjoint sections to the constant sections $\infty:=[1:0]$, $0:=[0:1]$, $1:=[1:1]$. Indeed, $\underline{\Isom}_S(X,\mathbb{P}^1_S)$ is a $\PGL_{2,S}$-torsor, and $\PGL_2$ acts simply transitively on ordered triples of pairwise distinct sections of $\mathbb{P}^1$. On the other hand,
the group of automorphisms of $\mathbb{P}^{1}_{\k}$ is $\PGL_2(\k)$ and for any three rational points $p_1,p_2,p_3\in\mathbb{P}^{1}_{\k}$ there is a unique automorphism $f:\mathbb{P}^{1}_{\k}\simeq \mathbb{P}^{1}_{\k}$ such that $f(p_1)=\infty:=[1:0]$, $f(p_2)=0:=[0:1]$, $f(p_3)=1:=[1:1]$. These two facts imply that, for $n\geq 3$, the stack $\M_{0,n}$ is representable by the scheme $(\mathbb{P}^{1}_{\k}\setminus\{\infty,0,1\})^{n-3} \setminus\Delta$, where $\Delta$ denotes the closed subscheme given by the diagonals. On $\k$-valued points the isomorphism reads as follows. The tuple $(\mathbb{P}^{1}_{\k},p_1,\dots,p_n)$ is equivalent to $(\mathbb{P}^{1}_{\k},\infty, 0, 1, f(p_4),\dots,f(p_n))$, where $f$ is the unique automorphism of $\mathbb{P}^{1}_{\k}$ with $f(p_1)=\infty$, $f(p_2)=0$ and $f(p_3)=1$. The tuple $(\mathbb{P}^{1}_{\k},p_1,\dots,p_n)$ is sent to $(f(p_4),\dots,f(p_n))$. We identify $\M_{0,n}$ with this scheme.

Define the morphism
\[
\vartheta:\Ht\times \bigl((\mathbb{P}^{1}_{\k}\setminus\{\infty,0,1\})^{n-3}\setminus\Delta\bigr)\longrightarrow \LS_{0,n,d}.
\]
Let $T$ be a $\k$-scheme. A $T$-valued point of
$\Ht\times \bigl((\mathbb{P}^{1}_{\k}\setminus\{\infty,0,1\})^{n-3}\setminus\Delta\bigr)$
is given by a pair
$\bigl((l_1,\dots,l_{n-1}),(p_4,\dots,p_n)\bigr)$,
where $(l_1,\dots,l_{n-1})\in \Ht(T)$ and
$(p_4,\dots,p_n)\in \bigl((\mathbb{P}^{1}_{\k}\setminus\{\infty,0,1\})^{n-3}\setminus\Delta\bigr)(T)$. The tuple $(p_4,\dots,p_n)$ determines sections
$\sigma_i:T\rightarrow \mathbb{P}^{1}_{T}:=T\times\mathbb{P}^{1}_{\k}$, for $i=4,\dots,n$,
which are pairwise disjoint and disjoint from the constant sections
$\sigma_1,\sigma_2,\sigma_3:T\rightarrow \mathbb{P}^{1}_{T}$
corresponding respectively to $p_1=\infty, p_2=0, p_3=1$. In this way we obtain a family of $n$-pointed smooth projective curves of genus $0$
$(\mathbb{P}^{1}_{T}\to T,\sigma_1,\dots,\sigma_n)$.
Let $L$ be the invertible sheaf on $\mathbb{P}^{1}_{T}$ defined by
$L:=\pr_2^{*}\O_{\mathbb{P}^{1}_{\k}}(dp_1)$,
where $\pr_2:\mathbb{P}^{1}_{T}\to \mathbb{P}^{1}_{\k}$ is the second projection. Let
$\gamma_1^0,\gamma_2^0,\gamma_3^0$
be the pullbacks to $T$ of the canonical trivializations of
$\O_{\mathbb{P}^{1}_{\k}}(dp_1)$ at the points $p_1,p_2,p_3$ introduced in \S~\ref{subsec:trivializations}. Since the sections $\sigma_4,\dots,\sigma_n$ do not meet the divisor $T\times\{p_1\}$, there are also canonical trivializations
$\gamma_i^0:\sigma_i^{*}L\xrightarrow{\sim} \O_T, \qquad i=4,\dots,n$. This gives a level
structure $\gs^0=(\mathbb{P}^{1}_{T},\sigma_1,\dots,\sigma_n,L,\gamma_1^0,\dots,\gamma_n^0)$.
Now we modify the trivializations by the parameters $(l_1,\dots,l_{n-1})$. More precisely, for $i=1,\dots,n-1$ we define
$\gamma_i:=l_i\circ \gamma_i^0:\sigma_i^{*}L\xrightarrow{\sim} \O_T$,
and we set
$\gamma_n:=\gamma_n^0$.
This determines a new level structure $\gs=(\mathbb{P}^{1}_{T},\sigma_1,\dots,\sigma_n,L,\gamma_1,\dots,\gamma_n)$.
We define
$\vartheta \bigl((l_1,\dots,l_{n-1}),(p_4,\dots,p_n)\bigr):=\gs$.
This construction is functorial in $T$, and so defines a morphism of functors
\[
\vartheta:\Ht\times \bigl((\mathbb{P}^{1}_{\k}\setminus\{\infty,0,1\})^{n-3}\setminus\Delta\bigr)\longrightarrow \LS_{0,n,d}.
\]

\begin{remark}\label{rmk:LS0nd-no-aut}
Since $g=0$ and $n\ge 3$, every automorphism of a $T$-valued level structure
$\gs\in\LS_{0,n,d}(T)$ is trivial. Indeed, an automorphism of the
underlying $n$-pointed smooth curve fixes the marked sections and so
induces, through the canonical normalization determined by the first three
sections, an automorphism of $\mathbb{P}^1_T$ fixing $\infty$, $0$ and $1$.
Such an automorphism is the identity. The stack $\LS_{0,n,d}$ has trivial
inertia, and is equivalent to the sheaf of isomorphism classes of its objects,
regarded as a category fibered in sets. In particular, a natural isomorphism on
$T$-valued points is in this case the same as an isomorphism of stacks.
\end{remark}

\begin{proposition}\label{prop:LS0nd-explicit}
The morphism
\[
\vartheta\colon \Ht\times \bigl((\mathbb{P}^{1}_{\k}\setminus\{\infty,0,1\})^{n-3}\setminus\Delta\bigr)
\longrightarrow \LS_{0,n,d}
\]
is an isomorphism of schemes. In particular, $\LS_{0,n,d}$ is the
$(2n-4)$-dimensional smooth and irreducible scheme
$\Ht\times \bigl((\mathbb{P}^{1}_{\k}\setminus\{\infty,0,1\})^{n-3}\setminus\Delta\bigr)$.
\end{proposition}

\begin{proof}
By Theorem~\ref{thm:goppa-representable}, $\LS_{0,n,d}$ is represented by a
scheme. It suffices to construct a natural inverse $\lambda$ to $\vartheta$ on
$T$-valued points for every $\k$-scheme $T$.

Let $\gs=(\pi:X\to T,\sigma_1,\dots,\sigma_n,L,\gamma_1,\dots,\gamma_n)\in \LS_{0,n,d}(T)$
be a $T$-valued point. Since we are in genus zero and $n\ge 3$, the ordered triple of sections
$\sigma_1,\sigma_2,\sigma_3$ determines a unique $T$-isomorphism
$\alpha:X\xrightarrow{\sim}\mathbb{P}^{1}_{T}$
such that $\alpha\circ \sigma_1,\alpha\circ \sigma_2,\alpha\circ \sigma_3$ are the constant sections corresponding to $p_1=\infty,\ p_2=0,\ p_3=1$.
Indeed, $\underline{\Isom}_{T}(X,\mathbb{P}^{1}_{T})$ is a $\PGL_{2,T}$-torsor, and $\PGL_2$ acts simply transitively on ordered triples of pairwise distinct points of $\mathbb{P}^1$.
Via $\alpha$, we may therefore regard $\gs$ as a level structure on the normalized family
\[
(\mathbb{P}^{1}_{T}\to T,\sigma_1^0,\sigma_2^0,\sigma_3^0,\tau_4,\dots,\tau_n),
\]
where $\sigma_1^0,\sigma_2^0,\sigma_3^0$ are the constant sections defined by $p_1,p_2,p_3$, and
\[
\tau_i:=\alpha\circ \sigma_i:T\to \mathbb{P}^{1}_{T}\qquad (4\le i\le n).
\]
Since the sections $\sigma_i$ are pairwise disjoint, the tuple $(\tau_4,\dots,\tau_n)$ determines a $T$-valued point of
$(\mathbb{P}^{1}_{\k}\setminus\{\infty,0,1\})^{n-3}\setminus\Delta$.

Next, by the standard description of the Picard group of a relative projective line, there exists a unique invertible sheaf $M$ on $T$ such that
$L\simeq \O_{\mathbb{P}^{1}_{T}}(dp_1)\otimes \pi^{*}M$.
Let $\gamma_i^0:\sigma_i^{*}\O_{\mathbb{P}^{1}_{T}}(dp_1)\xrightarrow{\sim} \O_T$ ($i=1,\dots,n$)
be the canonical trivializations of $\O_{\mathbb{P}^{1}_{T}}(dp_1)$ along the marked sections of the normalized family (the pullbacks of the trivializations of Subsection~\ref{subsec:trivializations}). The above decomposition of $L$ makes each given trivialization
$\gamma_i:\sigma_i^{*}L\xrightarrow{\sim} \O_T$
induce a trivialization
\[
\eta_i:M\xrightarrow{\sim} \O_T,
\qquad
\eta_i:=\gamma_i\circ (\gamma_i^0\otimes \mathrm{id}_M)^{-1}.
\]
In particular, $\eta_n$ trivializes $M$. Using this trivialization, we identify $M$ with $\O_T$, so $L\simeq \O_{\mathbb{P}^{1}_{T}}(dp_1)$.
Under this identification, each $\gamma_i$ becomes
$\gamma_i=l_i\circ \gamma_i^0$
for a unique unit $l_i\in \Gamma(T,\O_T)^{*}$, and by construction $l_n=1$.
Then the tuple $(l_1,\dots,l_{n-1})$ defines a $T$-valued point of $\Ht$.
Set
\[
\lambda_T(\gs):=\bigl((l_1,\dots,l_{n-1}),(\tau_4,\dots,\tau_n)\bigr).
\]
This construction is compatible with base change in $T$, and so defines a morphism of functors
\[
\lambda:\LS_{0,n,d}\longrightarrow \Ht\times \bigl((\mathbb{P}^{1}_{\k}\setminus\{\infty,0,1\})^{n-3}\setminus\Delta\bigr).
\]
The definitions give $\lambda_T\circ \vartheta_T=\mathrm{id}$ and $\vartheta_T\circ \lambda_T=\mathrm{id}$ for every $T$. Consequently, $\vartheta$ is an isomorphism of schemes, and $\LS_{0,n,d}$ is representable by $\Ht\times \bigl((\mathbb{P}^{1}_{\k}\setminus\{\infty,0,1\})^{n-3}\setminus\Delta\bigr)$.

Finally, since $\Ht$ is smooth and irreducible of dimension $n-1$, and
$(\mathbb{P}^{1}_{\k}\setminus\{\infty,0,1\})^{n-3}\setminus\Delta$
is a smooth irreducible scheme of dimension $n-3$, the product is smooth, irreducible, and of dimension $(n-1)+(n-3)=2n-4$.
\end{proof}
\begin{remark}
From now on we retain the explicit description
$\LS_{0,n,d}\simeq \Ht\times \M_{0,n}$
of Proposition~\ref{prop:LS0nd-explicit}. We shall use coordinates on the
standard affine charts of the Grassmannian, but all the constructions below are
to be understood functorially on $T$-valued points.
\end{remark}

On $T$-valued points, the Goppa morphism is given by
\begin{equation}
\begin{split}
\Goppa_{0,n,d}:\Ht\times \M_{0,n}&\rightarrow \Gr(k,n),\\
((l_1,\dots,l_{n-1}),(p_{4},\dots,p_{n}))&\mapsto \ker(P_\gs^{T})
\end{split}
\end{equation}
where $((l_1,\dots,l_{n-1}),(p_{4},\dots,p_{n}))=\vartheta^{-1}(\gs)$ and $P_\gs$ is the matrix \eqref{eq:basic_parity_check}. Likewise, the extended Goppa morphism is
\begin{equation}
\begin{split}
\Phi_{0,n,d}:\Ht\times \M_{0,n}&\rightarrow \Gr(k,n)\times \M_{0,n},\\
((l_1,\dots,l_{n-1}),(p_{4},\dots,p_{n}))&\mapsto (\ker(P_\gs^{T}),p_4,\dots,p_n).
\end{split}
\end{equation}

\subsection{Fibers of the Goppa morphism in genus zero}\label{subsec:fibers-g0}

We specialize Theorem~\ref{thm:fibers-goppa-projective-moduli} to $g=0$. Here
$k=d+1$ and $2g+1=1$, so the hypothesis $n>d\geq 2g+1$ of that theorem reads
$n>d\geq 1$. Below we assume moreover $2\leq d\leq n-2$. The case $d=n-1$ is
excluded because then $k=n$, so $C=\k^{n}$ and $\Gr(n,n)$ is a point.

Points $q_1,\dots,q_m\in\P^{r}$ are in \emph{linear general position} if any
$\min(m,r+1)$ of them are linearly independent.

We use the following two classical facts about rational normal curves, in a
form valid in arbitrary characteristic (compare \cite[Lecture~1]{harris-ag}).

\begin{lemma}\label{lem:rnc}
Let $\k$ be algebraically closed and let $d\geq 2$.
\begin{enumerate}
\item Any $d+1$ distinct points of a rational normal curve of degree $d$ in
$\P^{d}$ are linearly independent.
\item Two rational normal curves of degree $d$ in $\P^{d}$ through the same
$d+3$ points in linear general position coincide.
\end{enumerate}
\end{lemma}

\begin{proof}
For (1), after a projectivity the curve is the standard rational normal curve
\(
[s:t]\longmapsto[s^d:s^{d-1}t:\dots:t^d]
\)
(see \cite[\S~6A]{eisenbud-syzygies}). The matrix attached to $d+1$ distinct
points is a Vandermonde matrix, hence has nonzero determinant.

For (2), the reduced subscheme defined by the $d+3$ points is
Gorenstein and geometrically in linear general position. Hence it lies on a rational normal curve by \cite[Theorem~4.1]{eisenbud-popescu-gale}, and uniqueness is
Castelnuovo's $d+3$ theorem recalled in \cite[\S~4]{eisenbud-popescu-gale}.
\end{proof}

\begin{proposition}\label{prop:fibers-g0}
Assume $2\leq d\leq n-2$, and let $C\in\Gr(k,n)(\k)$ be a non-degenerate code,
with distinguished points $q_1,\dots,q_n\in\P_C$ as in
Lemma~\ref{lem:canonical-marked-projective}. Since $k=d+1$, the projective space $\P_C$ has dimension $d$. Let us consider the fiber $\mathcal{G}_C=\Goppa_{0,n,d}^{-1}(C)$.
The following statements hold.
\begin{enumerate}
\item If $d\leq n-3$, then $\mathcal{G}_C$ is either empty or
isomorphic to $\Spec\k$.
\item If $d=n-2$, then
\[
\mathcal{G}_C \simeq
\begin{cases}
\M_{0,n}, & \text{if } q_1,\dots,q_n \text{ are in linear general position},\\
\varnothing, & \text{otherwise,}
\end{cases}
\]
of dimension $n-3$ in the non-empty case.
\end{enumerate}
\end{proposition}

\begin{proof}
By Theorem~\ref{thm:fibers-goppa-projective-moduli} with $g=0$, we identify $\mathcal{G}_C$ with $\mathcal{P}_C$. On every
geometric fiber, the embedded curve is a rational normal curve of degree $d$
through $q_1,\dots,q_n$.
By Proposition~\ref{prop:LS0nd-explicit} the stack $\LS_{0,n,d}$ is
representable by a $\k$-scheme of finite type. The same holds for $\mathcal{G}_C=\LS_{0,n,d}\times_{\Gr(k,n)}\Spec\k$. In particular, $\mathcal{G}_C$ is Jacobson, and the residue fields at its closed points are
algebraic over $\k$ \cite[Tag~0478]{stacks}.

(1) Assume first that $d\leq n-3$, so that $n\geq d+3$.
Assume further that $\mathcal{G}_C$ is non-empty, and let $x$ be a closed point.
It determines a rational normal curve $R_x\subseteq\P_C\otimes_\k\kappa(x)$ of degree $d$ through $q_1,\dots,q_n$. 
By Lemma 7.16(1), after base change to an algebraic closure of $\kappa(x)$, the points $q_1,\ldots,q_{d+3}$ are in linear general position. Hence they are already in linear general position over $K$.
By \cite[Theorem~4.1]{eisenbud-popescu-gale}, there exists a rational normal curve $R\subseteq \P_C$ through these points which is unique. Indeed, after base change to an algebraic closure of $K$, uniqueness follows from
Lemma~7.16(2), and faithful flatness then gives uniqueness over $K$. After base change to $\kappa(x)$, the same uniqueness argument gives $R_{\kappa(x)}=R_x$. Hence, by faithful flatness, $R$ contains all $q_i$ and defines a $K$-rational point $[R]\in\mathcal G_C$.

If $y$ is any closed point of $\mathcal G_C$, then after base change
to an algebraic closure of $\kappa(y)$, Lemma~7.16(2) shows that the
corresponding rational normal curve is the base change of $R$.
Faithful flatness then gives the same equality over $\kappa(y)$.
Hence the associated $\kappa(y)$-point of $\mathcal G_C$ is the base
change of $[R]$, and therefore $y=[R]$.

Since $\mathcal{G}_C$ is of finite type over $\k$, it is Jacobson \cite[Tag~0478]{stacks}. Hence closed points are dense in every closed
subset. Since $[R]$ is the unique closed point, it follows that $\overline{\{z\}}=\{[R]\}$ for every $z\in\mathcal{G}_C$. Thus $\mathcal{G}_C$ has a single point.

Uniqueness of the closed point does not yet imply that $\mathcal{G}_C$ is
reduced. It remains to show that its tangent space vanishes. Let $\varepsilon^2=0$ and let $\xi\in\mathcal{G}_C(\k[\varepsilon])$ lie over $[R]$. Then $\xi$ determines a first-order deformation of $R$ as a closed subscheme
of $\P_C$, leaving $q_1,\dots,q_{d+3}$ fixed.
The space of such first-order deformations is $H^0\bigl(R,N_{R/\P_C}(-q_1-\cdots-q_{d+3})\bigr)$
(see \cite[\S~1.b, Proposition~1.2]{perrin-points}). 
The morphism which sends an object of $\mathcal{G}_C\simeq\mathcal{P}_C$ to the image of its embedding in $\P_C$ induces, on tangent spaces, a linear map \(T_{[R]}\mathcal{G}_C\to
H^0\bigl(R,N_{R/\P_C}(-q_1-\cdots-q_{d+3})\bigr)\).
This map is injective. Indeed, if the associated embedded deformation is
trivial, then its underlying embedded family is $R\times\Spec\k[\varepsilon]$, and its marked sections, being the preimages
of the $q_i$, are those of the trivial deformation of $[R]$. Hence $\xi$ is
trivial.

It remains to show that $H^0\bigl(R,N_{R/\P_C}(-q_1-\cdots-q_{d+3})\bigr)=0$.
After base change to $\bar{\k}$, one has $N_{R_{\bar{\k}}/\P_{C,\bar{\k}}}\simeq \O_{\mathbb{P}^1_{\bar{\k}}}(d+2)^{\oplus(d-1)}$ \cite[Corollary~2.6]{coskun-riedl-normal-bundles}. Hence the twist by
$q_1+\cdots+q_{d+3}$ is $\O_{\mathbb{P}^1_{\bar{\k}}}(-1)^{\oplus(d-1)}$, which has no nonzero
global sections. By compatibility of global sections with extension of the
base field, $H^0\bigl(R,N_{R/\P_C}(-q_1-\cdots-q_{d+3})\bigr)=0$, and therefore $T_{[R]}\mathcal{G}_C=0$.

Since $[R]$ is $\k$-rational, its residue field is $\k$. Since $T_{[R]}\mathcal{G}_C=0$, we have $\mathfrak m/\mathfrak m^2=0$, hence $\mathfrak m=\mathfrak m^2$.
As $\mathcal{G}_C$ is of finite type over $\k$, the maximal ideal $\mathfrak m$ of $\O_{\mathcal{G}_C,[R]}$ is finitely generated, and
Nakayama's lemma gives $\mathfrak m=0$. Since $\mathcal{G}_C$ has a single
point, it is affine, and therefore $\mathcal{G}_C=\Spec\O_{\mathcal{G}_C,[R]}=\Spec\k$.

(2) Assume now that $d=n-2$. If $\mathcal{G}_C$ is non-empty, then $q_1,\dots,q_n$
lie on a rational normal curve of degree $n-2$ in $\mathbb{P}^{n-2}$, on which
any $n-1$ distinct points are linearly independent by
Lemma~\ref{lem:rnc}(1), so $q_1,\dots,q_n$ are in linear general position.

Conversely, assume that $q_1,\dots,q_n$ are in linear general position, and let $(\pi:X\to T,\sigma_1,\dots,\sigma_n)$ be an object of $\M_{0,n}(T)$. The sheaf $\pi_{\ast}\O_X((n-2)\sigma_1)$ is locally free of rank $n-1$ and its formation
commutes with base change, since $R^{1}\pi_{\ast}\O_X((n-2)\sigma_1)=0$. Let $\nu$
be the closed immersion of $X$ into $\mathbb{P}(\pi_{\ast}\O_X((n-2)\sigma_1))$
defined by the relatively very ample sheaf $\O_X((n-2)\sigma_1)$ of relative degree $n-2$. 
Its geometric fibers are rational normal curves of degree $n-2$, so any $n-1$ of the sections $\nu\circ \sigma_i$ are fiberwise linearly independent by
Lemma~\ref{lem:rnc}(1). As $n=(n-2)+2$, these sections form a projective frame.
As in Subsection~\ref{subsec:space-and-goppa}, the sheaf of isomorphisms
between two projective bundles of the same rank is a $\PGL$-torsor, and $\PGL$
acts simply transitively on projective frames. There is a unique $T$-isomorphism $\varphi:\mathbb{P}(\pi_{\ast}\O_X((n-2)\sigma_1))\xrightarrow{\sim}\P_C\times T$ carrying $(\nu\circ \sigma_i)_i$ to $(q_i)_i$, and $\varphi\circ\nu$ defines an
object of $\mathcal{P}_C(T)$.

This construction is inverse to the functor which forgets the embedding $j$ and retains the underlying pointed curve. Indeed, for such an object the sheaf $L=j^{\ast}\O_{\P_C\times T}(1)$ has relative degree $n-2$, so 
$\mathbb{P}(\pi_{\ast}L)$ and $\mathbb{P}(\pi_{\ast}\O_X((n-2)\sigma_1))$ are canonically $T$-isomorphic.
Under this identification, both $j$ and $\varphi\circ\nu$ send each $\sigma_i$ to $q_i$. By the uniqueness of $\varphi$, we have $j=\varphi\circ\nu$.
Thus the two constructions are inverse to each other, and hence $\mathcal{P}_C\simeq\M_{0,n}$. By
Theorem~\ref{thm:fibers-goppa-projective-moduli}, it follows that $\mathcal{G}_C\simeq\M_{0,n}$, of dimension $n-3$.
\end{proof}

\begin{remark}
For $d=n-2$, Proposition~\ref{prop:fibers-g0}(2) recovers the open part of Kapranov's Hilbert-scheme realization of $\M_{0,n}$ as the space of rational
normal curves of degree $n-2$ through $n$ points in linear general position. Kapranov further proved that its closure in the Hilbert scheme is isomorphic to $\overline{\M}_{0,n}$ \cite[Theorem~0.1]{kapranov-veronese}.
\end{remark}

\subsection{The Goppa morphism as an immersion}\label{subsec:goppa-immersion-g0}

Theorem~\ref{th:immersionstacks} shows that the extended morphism $\Phi_{g,n,d}$ is an immersion whenever $n\geq d+g+2$. In genus zero this is the same as $d\leq n-2$. We prove here that in the smaller range $2\leq d\leq n-3$ the Goppa morphism itself is an immersion, so that the marked points need not be remembered since they are already encoded in the code. The bound is sharp in the
sense that for $d=n-2$ the fibers have dimension $n-3$ by
Proposition~\ref{prop:fibers-g0}(2).

\begin{theorem}\label{thm:goppa-immersion-g0}
Assume $2\leq d\leq n-3$ and set $k=d+1$. The Goppa morphism $\Goppa_{0,n,d}\colon \LS_{0,n,d}\longrightarrow \Gr(k,n)$
is an immersion, and it determines, for each $\lambda\in \Ht(\k)$, an immersion $\Goppa^{\lambda}_{0,n,d}:\M_{0,n}\hookrightarrow\Gr(k,n)$.
\end{theorem}

\begin{proof}
The second assertion follows directly from the first. We prove that the Goppa morphism is an immersion.

\emph{Step 1: reduction to the valuative criterion.} 
By Proposition~\ref{prop:LS0nd-explicit} and the projectivity of $\Gr(k,n)$, the morphism $\Goppa_{0,n,d}$ is of finite type. If $y\in\Gr(k,n)$ has non-empty fiber, choose a point $x$ in this fiber. Then $y_{\kappa(x)}$ is a genus-zero Goppa code, non-degenerate by \S~\ref{subsec:parameters}. Thus $y$ is non-degenerate, and
Proposition~\ref{prop:fibers-g0}(1), applied over $\kappa(y)$, gives $\Goppa_{0,n,d}^{-1}(y)\simeq\Spec\kappa(y)$. Therefore $\Goppa_{0,n,d}$ is a monomorphism \cite[Proposition~17.2.6]{ega4-4}. We verify the valuative condition of Lemma~\ref{lem:immersion-criterion}.

Let $A$ be a DVR with fraction field $K$, residue field $\kappa$, valuation $v$ and closed point $t_0$. Let $h:\Spec A\to\Gr(k,n)$ satisfy $h(t_0)\in\Goppa_{0,n,d}(\LS_{0,n,d})$, and let $\gs_\eta\in\LS_{0,n,d}(K)$ lift 
$h|_{\Spec K}$. An extension of $\gs_\eta$ to an $A$-point lifting $h$, if it exists, is unique since $\Goppa_{0,n,d}$ is a monomorphism.

\emph{Step 2: determinantal control of the extension.}
It remains to prove existence. By Proposition~\ref{prop:LS0nd-explicit},
after normalization it is enough to show that the marked points extend to
pairwise disjoint sections over $A$ and that the level parameters extend to units of $A$. We shall deduce both facts from the maximal minors of a generator matrix of $h$.

Normalize $\gs_\eta$ as in \S~\ref{subsec:space-and-goppa}, with $p_1=\infty$, $p_2=0$, $p_3=1$, parameters $l_1,\dots,l_{n-1}\in K^\ast$ and $l_n=1$. Properness of $\P^1$ gives
unique extensions $\widetilde p_i\in\P^1(A)$ of the marked points. Write $\widetilde p_i=[x_i:y_i]$, with $Ax_i+Ay_i=A$, and take $\widetilde p_1=[1:0]$, $\widetilde p_2=[0:1]$ and $\widetilde p_3=[1:1]$. Set
\((ab):=x_ay_b-x_by_a\). 
Then $(12),(13),(23)\in A^\ast$, $(1i)=y_i$ and $(2i)=-x_i$. Furthermore, $(ab)$ satisfies the following condition.
\begin{equation}\label{eq:disjoint}
(ab)\in A^\ast\Leftrightarrow \widetilde p_a \textrm{ is disjoint from } \widetilde p_b.
\end{equation}

We first write the generic generator matrix in homogeneous coordinates, so that its maximal minors record both the marked points and the level
parameters. Since $p_i\neq p_1$, one has $y_i\neq0$ in $K$ for $i\ge2$.
Set $w_i=(y_i^d,x_iy_i^{d-1},\dots,x_i^d)^T$, $W=(w_1|\cdots|w_n)$, $\mu_1=l_1$, and $\mu_i=l_i y_i^{-d}$ for $i\ge2$. On the generic fiber, $\alpha_i=x_i/y_i$ for $i\ge2$, so $w_i=y_i^d(1,\alpha_i,\dots,\alpha_i^d)^T$, while $w_1=(0,\dots,0,1)^T$. Hence, by \eqref{eq:basic_generator}, the generic
code has generator matrix $G=W\Diag(\mu_1,\dots,\mu_n)$.

Let $\mathcal E\subset A^n$ be the rank-$(d+1)$ subbundle defined by $h$.
Choose an $A$-basis of $\mathcal E$ and let $G'$ be the matrix whose rows
are the vectors of this basis. Since $G'$ and $G$ generate the same $K$-subspace of $K^n$, one has $G'=UG$ for some $U\in\GL_{d+1}(K)$. For $J\subseteq\{1,\dots,n\}$ with $|J|=d+1$, let
$G'_J$ denote the submatrix formed by the columns indexed by $J$, and set $\Delta_J=\det G'_J$. The homogeneous Vandermonde identity gives
\begin{equation}\label{eq:homogeneous-minors}
\Delta_J=\pm\det U\prod_{i\in J}\mu_i
\prod_{\substack{a<b\\a,b\in J}}(ab).
\end{equation}

We first show that every $\Delta_J$ is a unit. This will allow us to compare
valuations in \eqref{eq:homogeneous-minors}. Let $C_0$ be the special fiber
of $\mathcal E$. By the description of the fibers above, $C_0=C_{\gs_0}$ for some $\gs_0\in\LS_{0,n,d}(\kappa)$. The associated
points $q_1,\dots,q_n$ lie on a rational normal curve of degree $d$ and are
pairwise distinct. After base change to $\overline{\kappa}$,
Lemma~\ref{lem:rnc}(1) shows that any $d+1$ of them are linearly independent,
and hence the same holds over $\kappa$. Under the basis of $C_0$ given by $G'\otimes\kappa$, its columns represent the points $q_i$. Thus every $\Delta_J$ is nonzero modulo the maximal ideal of $A$, and therefore $\Delta_J\in A^\ast$.

\emph{Step 3: extension of the marked points.}
We now use these minors to prove that the extended marked points remain
pairwise disjoint. In the chosen homogeneous coordinates, this amounts to
showing that every $(ab)$ is a unit by \eqref{eq:disjoint}. Fix $i\ge4$ and choose $I\subseteq\{1,\dots,n\}$, with $|I|=d-1$, disjoint from $\{1,2,3,i\}$. This is possible since $n\ge d+3$. Write $\Delta_{Iuv}:=\Delta_{I\cup\{u,v\}}$. 
Cancellation in \eqref{eq:homogeneous-minors} gives
\[
\frac{\Delta_{I12}\Delta_{I3i}}
     {\Delta_{I13}\Delta_{I2i}}
 =\pm\frac{(12)(3i)}{(13)(2i)},
\qquad
\frac{\Delta_{I12}\Delta_{I3i}}
     {\Delta_{I1i}\Delta_{I23}}
 =\pm\frac{(12)(3i)}{(1i)(23)}.
\]
The left-hand sides are units, and so are $(12),(13),(23)$. Taking
valuations therefore gives $v(3i)=v(2i)$ and $v(3i)=v(1i)$. Since $(1i)=y_i$ and $(2i)=-x_i$, we obtain $v(x_i)=v(y_i)=v(3i)$. As $Ax_i+Ay_i=A$, at least one of $x_i,y_i$
is a unit. Hence their common valuation is zero. Thus $(1i),(2i),(3i)\in A^\ast$.
If $i,j\ge4$, choose $I$ of cardinality $d-1$ disjoint from $\{1,2,i,j\}$. The same cancellation gives
$$\frac{\Delta_{I12}\Delta_{Iij}}{\Delta_{I1i}\Delta_{I2j}}=\pm\frac{(12)(ij)}{(1i)(2j)}.$$ 
All factors except $(ij)$ are now units, so $(ij)\in A^\ast$. Thus every $(ab)$ is a unit, so the sections $\widetilde p_1,\dots,\widetilde p_n$ are pairwise disjoint. Moreover, $(1i)=y_i$ shows that $y_i\in A^\ast$ for every $i\ge2$.

\emph{Step 4: extension of the level parameters.}
It remains to show that the level parameters extend to units. Since all bracket factors in \eqref{eq:homogeneous-minors} are now units, quotients
of suitable maximal minors determine the ratios $\mu_i/\mu_n$ up to units. For $i<n$, choose a $(d+1)$-subset $J_i$ containing $n$ but not $i$, and
set $J_i'=(J_i\setminus\{n\})\cup\{i\}$. Then $\Delta_{J_i'}/\Delta_{J_i} =\pm(\mu_i/\mu_n)u_i$ for some $u_i\in A^\ast$. Since both minors are
units, $\mu_i/\mu_n\in A^\ast$. Moreover, $\mu_n=y_n^{-d}\in A^\ast$, because $l_n=1$ and $y_n$ is a unit. Consequently every $\mu_i$ is a unit. Since $l_1=\mu_1$ and $l_i=\mu_i y_i^d$ for $i\ge2$, all the $l_i$ are units.

We have therefore shown that the marked sections
are pairwise disjoint and the level parameters extend over $A$. By
Proposition~\ref{prop:LS0nd-explicit}, they define $\gs\in\LS_{0,n,d}(A)$ extending $\gs_\eta$. It remains only to check that
this extension lifts $h$. The morphisms $\Goppa_{0,n,d}(\gs)$ and $h$ agree over $K$. Since $A$ is a DVR and $\Gr(k,n)$ is separated,
they agree on all of $\Spec A$. Lemma~\ref{lem:immersion-criterion} now
shows that $\Goppa_{0,n,d}$ is an immersion.
\end{proof}

\subsection{Self-dual rational Goppa codes}\label{subsec:selfdual}

Fix $d>0$ and set $n=2(1+d)$. Every level structure $\gs\in\LS_{0,n,d}(\k)$ is isomorphic to one of the form
\[
\gs'=(\mathbb{P}^{1}_{\k},p_1,p_2,p_3,p_4,\dots,p_n,\O_{\mathbb{P}^{1}_{\k}}(d\infty),l_1\gamma_{1}^{0}, \dots, l_{n-1}\gamma_{n-1}^{0},l_n\gamma_{n}^{0})
\]
with $p_1=[1:0]$, $p_2=[0:1]$, $p_3=[1:1]$ and $l_n=1$. By
Lemma~\ref{lem:selfdual-level-code}, $\gs'$ is self-dual if and only if $C_{\gs'}$ is self-dual, so it suffices to analyze level structures of this
shape.

\begin{proposition}\label{prop:selfdual-g0}
Fix $d>0$, set $n=2(1+d)$, and let $\gs'$ be as above. For $i=2,\dots,n$, set $h_i:=\bigl(\prod_{j\neq1,i}(\alpha_i-\alpha_j)\bigr)^{-1}$. Then $\gs'$
(and equivalently $C_{\gs'}$) is self-dual if and only if
\[
l_1^2=-h_n^{-1},\qquad l_i^2=h_i h_n^{-1},\quad i=2,\dots,n-1.
\]
Equivalently, if $\eta_0:=dt/\prod_{j=2}^{n}(t-\alpha_j)$, then $\eta:=h_n^{-1}\eta_0$ satisfies $(\eta)=(2d-1)p_1-\sum_{i=2}^{n}p_i$, $\Res_{p_i}(\eta)=l_i^2$ for $i=2,\dots,n$, and $\Res_{p_1}(\eta/u^{2d})=l_1^2$, where $u=t^{-1}$.
\end{proposition}

\begin{proof}
Since $n=2d+2$, one has
$(\eta_0)=(2d-1)p_1-\sum_{i=2}^{n}p_i$. By
Lemma~\ref{lm:char-dual} and the trivializations of
Subsection~\ref{subsec:trivializations}, $\gs'$ is self-dual if and only if
there is a nonzero meromorphic differential $\eta$ with $(\eta)\geq(2d-1)p_1-\sum_{i=2}^{n}p_i$ such that $\Res_{p_i}(\eta)=l_i^2$ for $i=2,\dots,n$ and $\Res_{p_1}(\eta/u^{2d})=l_1^2$. Here $u^{-2d}$ is a local generator of
$\O_{\mathbb P^1}(2dp_1)$ at $p_1$. Since both divisors have degree $-2$,
the displayed inequality is an equality, and hence $\eta=c\eta_0$ for some $c\in\k^*$. For $i\geq2$, $\Res_{p_i}(\eta_0)=h_i$, while
\[
\frac{\eta_0}{u^{2d}}=-\frac{du}{u\prod_{j=2}^{n}(1-\alpha_j u)},
\]
so $\Res_{p_1}(\eta_0/u^{2d})=-1$. Thus the residue conditions are
$l_i^2=ch_i$ for $i=2,\dots,n$ and $l_1^2=-c$. Since $l_n=1$, necessarily
$c=h_n^{-1}$, and the stated conditions follow.
\end{proof}

\section*{Acknowledgement}
The author thanks F.~J. Plaza Mart\'in for his enlightening comments and suggestions.

\bibliographystyle{plain}
\bibliography{mibiblio7}

@book{oxley-matroid-theory,
  author    = {Oxley, J.},
  title     = {Matroid Theory},
  series    = {Oxford Graduate Texts in Mathematics},
  volume    = {21},
  edition   = {2nd},
  publisher = {Oxford University Press},
  address   = {Oxford},
  year      = {2011},
  isbn      = {978-0-19-960339-8},
}

@article{xing-equal,
  author  = {C.-P. Xing},
  title   = {When are two geometric {G}oppa codes equal?},
  journal = {IEEE Transactions on Information Theory},
  volume  = {38},
  number  = {3},
  pages   = {1140--1142},
  year    = {1992},
  doi     = {10.1109/18.135656},
}

@article{eisenbud-popescu-gale,
  author  = {Eisenbud, D. and Popescu, S.},
  title   = {The Projective Geometry of the Gale Transform},
  journal = {Journal of Algebra},
  volume  = {230},
  number  = {1},
  pages   = {127--173},
  year    = {2000},
  doi     = {10.1006/jabr.1999.7940}
}

@incollection{milne-abelian-varieties,
  author    = {Milne, J. S.},
  title     = {Abelian Varieties},
  booktitle = {Arithmetic Geometry},
  editor    = {Cornell, Gary and Silverman, Joseph H.},
  pages     = {103--150},
  publisher = {Springer},
  address   = {New York},
  year      = {1986},
  doi       = {10.1007/978-1-4613-8655-1_5}
}

@incollection{alvarez2,
  author    = {\'Alvarez-V\'azquez, A. and Mu\~noz Porras, J. M.},
  title     = {Canonical immersions of moduli schemes of vector bundles into Grassmannians},
  booktitle = {Workshop on Abelian Varieties and Theta Functions (Morelia, 1996)},
  series    = {Aportaciones Mat. Investig.},
  volume    = {13},
  pages     = {87--95},
  publisher = {Sociedad Matem\'atica Mexicana},
  address   = {M\'exico},
  year      = {1998},
}

@article{Gop77,
  author  = {Goppa, V. D.},
  title   = {Codes Associated with Divisors},
  journal = {Problemy Peredachi Informatsii},
  volume  = {13},
  number  = {1},
  pages   = {33--39},
  year    = {1977},
}

@article{feng-rao,
  author  = {Feng, G.-L. and Rao, T. R. N.},
  title   = {Decoding algebraic-geometric codes up to the designed minimum distance},
  journal = {IEEE Transactions on Information Theory},
  volume  = {39},
  number  = {1},
  pages   = {37--45},
  year    = {1993},
  doi     = {10.1109/18.179340},
}

@misc{savin-vector-bundles,
  author       = {Savin, V.},
  title        = {Algebraic-geometric codes from vector bundles and their decoding},
  howpublished = {arXiv:0803.1096},
  year         = {2008},
  primaryClass = {cs.IT},
}

@book{eisenbud-syzygies,
  author    = {Eisenbud, D.},
  title     = {The Geometry of Syzygies: A Second Course in Algebraic Geometry
               and Commutative Algebra},
  series    = {Graduate Texts in Mathematics},
  volume    = {229},
  publisher = {Springer},
  address   = {New York},
  year      = {2005},
  doi       = {10.1007/b137572}
}

@article{coskun-riedl-normal-bundles,
  author  = {Coskun, I. and Riedl, E.},
  title   = {Normal bundles of rational curves on complete intersections},
  journal = {Communications in Contemporary Mathematics},
  volume  = {21},
  number  = {2},
  year    = {2019},
  pages   = {1850011},
  doi     = {10.1142/S0219199718500116}
}

@article{ega4-4,
  author  = {Grothendieck, A.},
  title   = {\'{E}l\'{e}ments de g\'{e}om\'{e}trie alg\'{e}brique. {IV}. \'{E}tude locale des sch\'{e}mas et des morphismes de sch\'{e}mas. {IV}},
  journal = {Publications Math\'{e}matiques de l'IH\'{E}S},
  volume  = {32},
  pages   = {5--361},
  year    = {1967},
}

@article{saint-donat-equations,
   author  = {Saint-Donat, B.},
   title   = {Sur les \'equations d\'efinissant une courbe alg\'ebrique},
   journal = {C. R. Acad. Sci. Paris S\'er. A},
   volume  = {274},
   year    = {1972},
   pages   = {324--327; 487--489},
 }

@incollection{fujita-defining,
   author    = {Fujita, T.},
   title     = {Defining equations for certain types of polarized varieties},
   booktitle = {Complex analysis and algebraic geometry},
   pages     = {165--173},
   publisher = {Iwanami Shoten, Tokyo},
   year      = {1977},
 }

@incollection{fulton-pandharipande,
  author    = {Fulton, W. and Pandharipande, R.},
  title     = {Notes on stable maps and quantum cohomology},
  booktitle = {Algebraic Geometry---Santa Cruz 1995},
  series    = {Proc. Sympos. Pure Math.},
  volume    = {62, Part 2},
  pages     = {45--96},
  publisher = {Amer. Math. Soc., Providence, RI},
  year      = {1997},
}

@article{pellikaan-which,
  author  = {Pellikaan, R. and Shen, B.-Z. and van Wee, G. J. M.},
  title   = {Which Linear Codes Are Algebraic-Geometric?},
  journal = {IEEE Transactions on Information Theory},
  volume  = {37},
  number  = {3},
  pages   = {583--602},
  year    = {1991},
  doi     = {10.1109/18.79915},
}

@article{curto-convolutional,
  author  = {Iglesias Curto, J. I. and Mu\~noz Casta\~neda, \'A. L. and Mu\~noz Porras, J. M. and Serrano Sotelo, G.},
  title   = {Every Convolutional Code is a {Goppa} Code},
  journal = {IEEE Transactions on Information Theory},
  volume  = {59},
  number  = {10},
  pages   = {6628--6641},
  year    = {2013},
  doi     = {10.1109/TIT.2013.2271033},
}

@book{Tsfasman,
  author    = {Tsfasman, M. and Vl{\u{a}}du{\c{t}}, S. and Nogin, D.},
  title     = {Algebraic Geometric Codes: Basic Notions},
  series    = {Mathematical Surveys and Monographs},
  volume    = {139},
  publisher = {American Mathematical Society},
  address   = {Providence, RI},
  year      = {2007},
  isbn      = {978-0-8218-4306-2},
}

@book{hartshorne,
  author    = {Hartshorne, R.},
  title     = {Algebraic Geometry},
  series    = {Graduate Texts in Mathematics},
  volume    = {52},
  publisher = {Springer-Verlag},
  address   = {New York-Heidelberg},
  year      = {1977},
  pages     = {xvi+496},
  isbn      = {0-387-90244-9},
  doi       = {10.1007/978-1-4757-3849-0},
  mrnumber  = {0463157},
  zblnumber = {0367.14001},
}

@article{Gop82,
  author  = {Goppa, V. D.},
  title   = {Algebraico-geometric codes},
  journal = {Izvestiya Akademii Nauk SSSR. Seriya Matematicheskaya},
  volume  = {46},
  number  = {4},
  pages   = {762--781},
  year    = {1982},
}

@book{stichtenoth,
  author    = {Stichtenoth, H.},
  title     = {Algebraic Function Fields and Codes},
  series    = {Graduate Texts in Mathematics},
  volume    = {254},
  edition   = {2nd},
  publisher = {Springer-Verlag},
  address   = {Berlin},
  year      = {2009},
  isbn      = {978-3-540-76877-7},
}

@article{munuera-pellikaan,
  author  = {Munuera, C. and Pellikaan, R.},
  title   = {Equality of geometric {G}oppa codes and equivalence of divisors},
  journal = {Journal of Pure and Applied Algebra},
  volume  = {90},
  number  = {3},
  pages   = {229--252},
  year    = {1993},
  doi     = {10.1016/0022-4049(93)90043-S},
}

@article{raynaud,
  author  = {Raynaud, M.},
  title   = {Sp\'{e}cialisation du foncteur de {P}icard},
  journal = {Publications Math\'{e}matiques de l'IH\'{E}S},
  volume  = {38},
  pages   = {27--76},
  year    = {1970},
}

@book{laumon,
  author    = {Laumon, G. and Moret-Bailly, L.},
  title     = {Champs alg\'{e}briques},
  series    = {Ergebnisse der Mathematik und ihrer Grenzgebiete. 3. Folge},
  volume    = {39},
  publisher = {Springer-Verlag},
  address   = {Berlin},
  year      = {2000},
  isbn      = {3-540-65761-4},
}

@misc{stacks,
  author       = {{The Stacks Project Authors}},
  title        = {The Stacks Project},
  howpublished = {\url{https://stacks.math.columbia.edu}},
  year         = {2024},
}

@article{ega4-2,
  author  = {Grothendieck, A.},
  title   = {\'{E}l\'{e}ments de g\'{e}om\'{e}trie alg\'{e}brique. {IV}. \'{E}tude locale des sch\'{e}mas et des morphismes de sch\'{e}mas. {II}},
  journal = {Publications Math\'{e}matiques de l'IH\'{E}S},
  volume  = {24},
  pages   = {5--231},
  year    = {1965},
}

@article{ega4-3,
  author  = {Grothendieck, A.},
  title   = {\'{E}l\'{e}ments de g\'{e}om\'{e}trie alg\'{e}brique. {IV}. \'{E}tude locale des sch\'{e}mas et des morphismes de sch\'{e}mas. {III}},
  journal = {Publications Math\'{e}matiques de l'IH\'{E}S},
  volume  = {28},
  pages   = {5--255},
  year    = {1966},
}

@incollection{Artin,
  author    = {Artin, M.},
  title     = {Algebraization of formal moduli: {I}},
  booktitle = {Global Analysis (Papers in Honor of {K}.~Kodaira)},
  editor    = {Spencer, D.~C. and Iyanaga, S.},
  pages     = {21--71},
  publisher = {University of Tokyo Press and Princeton University Press},
  year      = {1969},
}

@article{butler,
  author  = {Butler, D. C.},
  title   = {Normal generation of vector bundles over a curve},
  journal = {Journal of Differential Geometry},
  volume  = {39},
  number  = {1},
  pages   = {1--34},
  year    = {1994},
  doi     = {10.4310/jdg/1214454673},
}

@book{alper,
  author    = {Alper, J.},
  title     = {Stacks and Moduli},
  publisher = {Unpublished lecture notes},
  year      = {2026},
}

@article{marquez4,
  author  = {M\'{a}rquez-Corbella, I. and Mart\'{\i}nez-Moro, E. and Pellikaan, R.},
  title   = {On the unique representation of very strong algebraic geometry codes},
  journal = {Designs, Codes and Cryptography},
  volume  = {70},
  number  = {1--2},
  pages   = {215--230},
  year    = {2014},
  doi     = {10.1007/s10623-012-9758-3}
}

@article{eisenbud-green-harris,
  author  = {Eisenbud, D. and Green, M. and Harris, J.},
  title   = {Cayley--{B}acharach theorems and conjectures},
  journal = {Bulletin of the American Mathematical Society (N.S.)},
  volume  = {33},
  number  = {3},
  pages   = {295--324},
  year    = {1996},
  doi     = {10.1090/S0273-0979-96-00666-0},
}

@book{harris-ag,
  author    = {Harris, J.},
  title     = {Algebraic Geometry: A First Course},
  series    = {Graduate Texts in Mathematics},
  volume    = {133},
  publisher = {Springer-Verlag},
  address   = {New York},
  year      = {1992},
  isbn      = {0-387-97716-3},
  doi       = {10.1007/978-1-4757-2189-8},
}

@article{deligne-mumford,
  author  = {Deligne, P. and Mumford, D.},
  title   = {The irreducibility of the space of curves of given genus},
  journal = {Publications Math\'{e}matiques de l'IH\'{E}S},
  volume  = {36},
  pages   = {75--109},
  year    = {1969},
  doi     = {10.1007/BF02684599},
}

@article{knudsen,
  author  = {Knudsen, F. F.},
  title   = {The projectivity of the moduli space of stable curves, {II}: {T}he stacks $M_{g,n}$},
  journal = {Mathematica Scandinavica},
  volume  = {52},
  number  = {2},
  pages   = {161--199},
  year    = {1983},
}

@incollection{kleiman-picard,
  author    = {Kleiman, S. L.},
  title     = {The {P}icard scheme},
  booktitle = {Fundamental Algebraic Geometry: {G}rothendieck's {FGA} Explained},
  series    = {Mathematical Surveys and Monographs},
  volume    = {123},
  pages     = {235--321},
  publisher = {American Mathematical Society},
  address   = {Providence, RI},
  year      = {2005},
}

@article{driencourt-stichtenoth,
  author  = {Driencourt, Y. and Stichtenoth, H.},
  title   = {A criterion for self-duality of geometric codes},
  journal = {Communications in Algebra},
  volume  = {17},
  number  = {4},
  pages   = {885--898},
  year    = {1989},
  doi     = {10.1080/00927878908823765},
}

@article{stichtenoth-selfdual,
  author  = {Stichtenoth, H.},
  title   = {Self-dual {G}oppa codes},
  journal = {Journal of Pure and Applied Algebra},
  volume  = {55},
  number  = {1--2},
  pages   = {199--211},
  year    = {1988},
  doi     = {10.1016/0022-4049(88)90046-1},
}

@article{couvreur-mcc-ag,
  author  = {Couvreur, A. and M\'{a}rquez-Corbella, I. and Pellikaan, R.},
  title   = {Cryptanalysis of {M}c{E}liece cryptosystem based on algebraic geometry codes and their subcodes},
  journal = {IEEE Transactions on Information Theory},
  volume  = {63},
  number  = {8},
  pages   = {5404--5418},
  year    = {2017},
}

@article{perrin-points,
  author  = {Perrin, D.},
  title   = {Courbes passant par $m$ points g\'{e}n\'{e}raux de $\mathbb{P}^3$},
  journal = {M\'{e}moires de la Soci\'{e}t\'{e} Math\'{e}matique de France (N.S.)},
  number  = {28--29},
  pages   = {1--137},
  year    = {1987},
  doi     = {10.24033/msmf.330},
}

@incollection{grothendieck-hilbert,
  author    = {Grothendieck, A.},
  title     = {Techniques de construction et th\'{e}or\`{e}mes d'existence en
               g\'{e}om\'{e}trie alg\'{e}brique. {IV}. {L}es sch\'{e}mas de {H}ilbert},
  booktitle = {S\'{e}minaire Bourbaki, Vol.~6 (1960/1961), Exp.~No.~221},
  pages     = {249--276},
  publisher = {Soci\'{e}t\'{e} Math\'{e}matique de France},
  address   = {Paris},
  year      = {1995},
}

@article{kapranov-veronese,
  author  = {Kapranov, M. M.},
  title   = {Veronese curves and the {G}rothendieck--{K}nudsen moduli space $\overline{M}_{0,n}$},
  journal = {Journal of Algebraic Geometry},
  volume  = {2},
  number  = {2},
  pages   = {239--262},
  year    = {1993},
}

@article{tsfasman-vladut-zink,
  author  = {Tsfasman, M. A. and Vl{\u{a}}du{\c{t}}, S. G. and Zink, T.},
  title   = {Modular curves, {S}himura curves, and {G}oppa codes, better than {V}arshamov--{G}ilbert bound},
  journal = {Mathematische Nachrichten},
  volume  = {109},
  pages   = {21--28},
  year    = {1982},
  doi     = {10.1002/mana.19821090103},
}

@article{janwa-moreno,
  author  = {Janwa, H. and Moreno, O.},
  title   = {{McEliece} public key cryptosystems using algebraic-geometric codes},
  journal = {Designs, Codes and Cryptography},
  volume  = {8},
  number  = {3},
  pages   = {293--307},
  year    = {1996},
}

@article{sidelnikov-shestakov,
  author  = {Sidelnikov, V. M. and Shestakov, S. O.},
  title   = {On insecurity of cryptosystems based on generalized {R}eed--{S}olomon codes},
  journal = {Discrete Mathematics and Applications},
  volume  = {2},
  number  = {4},
  pages   = {439--444},
  year    = {1992},
  doi     = {10.1515/dma.1992.2.4.439},
}

@book{demazure-gabriel,
  author    = {Demazure, M. and Gabriel, P.},
  title     = {Groupes alg\'ebriques. {T}ome~{I}: G\'eom\'etrie alg\'ebrique,
               g\'en\'eralit\'es, groupes commutatifs},
  publisher = {Masson \& Cie / North-Holland},
  address   = {Paris / Amsterdam},
  year      = {1970},
}

@article{conv-goppa-type,
  author  = {Dom\'{\i}nguez P\'erez, J. \'A. and Mu\~noz Porras, J. M. and Serrano Sotelo, G.},
  title   = {Convolutional Codes of {Goppa} Type},
  journal = {Applicable Algebra in Engineering, Communication and Computing},
  volume  = {15},
  number  = {1},
  pages   = {51--61},
  year    = {2004},
  doi     = {10.1007/s00200-004-0151-y},
}

@article{convolutional-goppa,
  author  = {Mu\~noz Porras, J. M. and Dom\'{\i}nguez P\'erez, J. \'A. and Iglesias Curto, J. I. and Serrano Sotelo, G.},
  title   = {Convolutional {Goppa} Codes},
  journal = {IEEE Transactions on Information Theory},
  volume  = {52},
  number  = {1},
  pages   = {340--344},
  year    = {2006},
  doi     = {10.1109/TIT.2005.860447},
}

@incollection{ag-convolutional-survey,
  author    = {Dom\'{\i}nguez P\'erez, J. \'A. and Mu\~noz Porras, J. M. and Serrano Sotelo, G.},
  title     = {Algebraic Geometry Constructions of Convolutional Codes},
  booktitle = {Advances in Algebraic Geometry Codes},
  editor    = {Mart\'{\i}nez-Moro, Edgar and Munuera, Carlos and Ruano, Diego},
  series    = {Series on Coding Theory and Cryptology},
  volume    = {5},
  pages     = {391--417},
  publisher = {World Scientific},
  address   = {Singapore},
  year      = {2008},
  doi       = {10.1142/9789812794017_0011},
}

@article{mds-convolutional-goppa,
  author  = {Plaza-Mart\'{\i}n, F. J. and Iglesias-Curto, J. I. and Serrano-Sotelo, G.},
  title   = {On the Construction of 1-D {MDS} Convolutional {Goppa} Codes},
  journal = {IEEE Transactions on Information Theory},
  volume  = {59},
  number  = {7},
  pages   = {4615--4625},
  year    = {2013},
  doi     = {10.1109/TIT.2013.2251926},
}

@article{keel-tevelev,
  author  = {Keel, S. and Tevelev, J.},
  title   = {Equations for $\overline{M}_{0,n}$},
  journal = {International Journal of Mathematics},
  volume  = {20},
  number  = {9},
  pages   = {1159--1184},
  year    = {2009},
}

@article{hmsv,
  author  = {Howard, B. and Millson, J. and Snowden, A. and Vakil, R.},
  title   = {The Equations for the Moduli Space of $n$ Points on the Line},
  journal = {Duke Mathematical Journal},
  volume  = {146},
  number  = {2},
  pages   = {175--226},
  year    = {2009},
  doi     = {10.1215/00127094-2008-063},
}

@book{schmitt-decorated-principal,
  author    = {Schmitt, A. H. W.},
  title     = {Geometric Invariant Theory and Decorated Principal Bundles},
  series    = {Zurich Lectures in Advanced Mathematics},
  publisher = {European Mathematical Society},
  address   = {Z\"urich},
  year      = {2008},
  isbn      = {978-3-03719-065-4},
}

@article{gomez-langer-schmitt-sols,
  author  = {G\'omez, T. L. and Langer, A. and Schmitt, A. H. W. and Sols, I.},
  title   = {Moduli Spaces for Principal Bundles in Arbitrary Characteristic},
  journal = {Advances in Mathematics},
  volume  = {219},
  number  = {4},
  pages   = {1177--1245},
  year    = {2008},
  doi     = {10.1016/j.aim.2008.05.015},
}

@article{ghorpade-lachaud-mds,
  author  = {Ghorpade, S. R. and Lachaud, G.},
  title   = {Hyperplane Sections of {Grassmannians} and the Number of {MDS} Linear Codes},
  journal = {Finite Fields and Their Applications},
  volume  = {7},
  number  = {4},
  pages   = {468--506},
  year    = {2001},
  doi     = {10.1006/ffta.2000.0299},
}

@article{ghorpade-lachaud-etale,
  author  = {Ghorpade, S. R. and Lachaud, G.},
  title   = {\'{E}tale Cohomology, {Lefschetz} Theorems and Number of Points of Singular Varieties over Finite Fields},
  journal = {Moscow Mathematical Journal},
  volume  = {2},
  number  = {3},
  pages   = {589--631},
  year    = {2002},
  doi     = {10.17323/1609-4514-2002-2-3-589-631},
}

@incollection{deligne-determinant,
  author    = {Deligne, P.},
  title     = {Le d\'{e}terminant de la cohomologie},
  booktitle = {Current Trends in Arithmetical Algebraic Geometry},
  series    = {Contemporary Mathematics},
  volume    = {67},
  pages     = {93--177},
  publisher = {American Mathematical Society},
  address   = {Providence, RI},
  year      = {1987},
}

@incollection{nitsure-hilbert,
  author    = {Nitsure, N.},
  title     = {Construction of {H}ilbert and {Q}uot Schemes},
  booktitle = {Fundamental Algebraic Geometry: Grothendieck's {FGA} Explained},
  editor    = {Fantechi, Barbara and G{\"o}ttsche, Lothar and Illusie, Luc and
               Kleiman, Steven L. and Nitsure, Nitin and Vistoli, Angelo},
  series    = {Mathematical Surveys and Monographs},
  volume    = {123},
  pages     = {105--137},
  publisher = {American Mathematical Society},
  address   = {Providence, RI},
  year      = {2005}
}

@article{henocq-rotillon,
  author  = {Henocq, T. and Rotillon, D.},
  title   = {The theta divisor of a {J}acobian variety and the decoding of geometric {G}oppa codes},
  journal = {Journal of Pure and Applied Algebra},
  volume  = {112},
  number  = {1},
  pages   = {13--28},
  year    = {1996},
  doi     = {10.1016/0022-4049(95)00008-9}
}

@article{beelen-grs,
  author  = {Beelen, P. and Glynn, D. and H{\o}holdt, T. and Kaipa, K.},
  title   = {Counting generalized {R}eed--{S}olomon codes},
  journal = {Advances in Mathematics of Communications},
  volume  = {11},
  number  = {4},
  pages   = {777--790},
  year    = {2017},
  doi     = {10.3934/amc.2017057}
}

\end{document}